\documentclass[10pt]{amsart}
\usepackage{amsmath,amssymb,amsthm,mathtools,mathrsfs}
\usepackage[left=2.5cm,right=2.5cm,top=2.5cm,bottom=2.5cm]{geometry}
\usepackage{enumitem}
\usepackage{hyperref}
\usepackage{verbatim}
\usepackage{cite}
\hypersetup{colorlinks=true,linkcolor=blue,citecolor=blue,filecolor=blue,urlcolor=blue,
pdftitle={Isolated yet highly unstable 2D Euler equilibria},
pdfauthor={Alberto Enciso, Antonio J. Fernandez, Marc Nualart}}
\renewcommand\leq\leqslant
\renewcommand\geq\geqslant
\setlist[enumerate]{label=\textup{(\roman*)},leftmargin=*}
\numberwithin{equation}{section}
\newtheorem{theorem}{Theorem}[section]
\newtheorem{introtheorem}[theorem]{Theorem}
\newtheorem{lemma}[theorem]{Lemma}
\newtheorem{proposition}[theorem]{Proposition}
\newtheorem{corollary}[theorem]{Corollary}
\newtheorem*{deconjecture}{Conjecture}
\theoremstyle{definition}

\newtheorem{remark}[theorem]{Remark}
\newtheorem*{remark*}{Remark}

\DeclareMathOperator{\dist}{dist}

\DeclareMathOperator{\osc}{osc}

\DeclareMathOperator{\curl}{curl}
\DeclareMathOperator{\diver}{div}
\newcommand{\R}{\mathbb R}

\newcommand{\Z}{\mathbb Z}
\newcommand{\Sone}{\mathbb S^1}
\newcommand{\eps}{\varepsilon}
\newcommand{\dd}{\,d}
\begin{document}
\title{Isolated yet highly unstable 2D Euler equilibria}

\author[A. Enciso]{Alberto Enciso}
\address{   
\newline
{\upshape\bfseries\small Alberto Enciso} 
\vspace{0.1cm}
\newline \indent Instituto de Ciencias Matem\'aticas, Consejo Superior de Investigaciones Cient\'\i ficas, 28049 Madrid, \newline \indent Spain}
\email{aenciso@icmat.es}

 \author[A. J. Fern\'andez]{Antonio J.\ Fern\'andez}
 \address{ \vspace{-0.5cm}
\newline 
{\upshape\bfseries\small Antonio J. Fern\'andez} 
\vspace{0.1cm}
\newline \indent Departamento de Matem\'aticas, Universidad Aut\'onoma de Madrid, 28049 Madrid, Spain}
 \email{antonioj.fernandez@uam.es}

 \author[M. Nualart]{Marc Nualart}
\address{ \vspace{-0.5cm}
\newline
{\upshape\bfseries\small Marc Nualart}
\vspace{0.1cm}
\newline \indent Instituto de Ciencias Matem\'aticas, Consejo Superior de Investigaciones Cient\'\i ficas, 28049 Madrid, \newline \indent Spain}
\email{marc.nualart@icmat.es}

\begin{abstract}
We construct steady solutions of the two-dimensional Euler equations on
smooth annular domains that are isolated, modulo scaling, in the $L^\infty$
vorticity topology. They are highly unstable: any Yudovich trajectory
remaining sufficiently close to such an equilibrium on an open time
interval is a stationary multiple of it, they admit wandering neighborhoods of initial data for all non-zero time and, moreover, at fixed energy,
enstrophy and circulations, distinct nearby trajectories cannot return
even weakly to the equilibrium. Further, on any non-circular annular domain, we show that this kind of isolated, highly unstable equilibria occurs arbitrarily close in the Yudovich norm to smooth equilibria which are themselves non-isolated.  
\end{abstract}
\maketitle
\section{Introduction}
\label{sec:introduction}

Consider the incompressible Euler equations on a bounded smooth connected
domain $\Omega\subset\mathbb R^2$:
\begin{equation}\label{intro:euler}
 \partial_tu+(u\cdot\nabla)u+\nabla p=0,\qquad
 \operatorname{div}u=0,\qquad u\cdot n|_{\partial\Omega}=0.
\end{equation}
Here $u$ is the velocity of the fluid, $p$ is the pressure, and $n$ denotes
the outward unit normal. The vorticity $\omega=\operatorname{curl}u$ is
transported by the flow, so its distribution function and all its $L^p$
norms are independent of time. The evolution also preserves the kinetic
energy $E(u)=\frac12\int_\Omega|u|^2\,dx$ and
the circulations along the boundary components.

A natural setting for the evolution problem is the Yudovich class
\[
 \mathcal Y(\Omega)=
 \{u\in L^2(\Omega;\mathbb R^2):
   \operatorname{div}u=0,\ u\cdot n|_{\partial\Omega}=0,\
   \operatorname{curl}u\in L^\infty(\Omega)\},
\]
with norm $\|u\|_{\mathcal Y}:=\|u\|_{L^2(\Omega)}
+\|\operatorname{curl}u\|_{L^\infty(\Omega)}$.
Every initial datum in this class gives rise to a unique global solution,
with $u\in C(\mathbb R;L^2)$ and
$\omega\in C_{\mathrm{w}^*}(\mathbb R;L^\infty)\cap C(\mathbb R;L^p)$ for
$1\leq p<\infty$, confer \cite{Yudovich1963,Sueur2011,LacaveMiotWang2014}.
Here $\mathrm{w}^*$ denotes the weak-$*$ topology. The equation is globally well-posed, but the long-time dynamics are still insufficiently understood. A major problem is that, although the vorticity distribution is conserved along the evolution, the weak limits of the vorticity orbits can be very different due to the formation of small scales.

Steady solutions are the most basic candidates for such limits.
Their possible role depends in part on the local structure of the set of equilibria.
Indeed, if a steady flow is accumulated by other equilibria which are not
multiples of it, its scaling line cannot attract a whole neighborhood of initial data.
As pointed out by Drivas and Elgindi
\cite[Problem 1]{DrivasElgindi2023}, this
observation makes the existence of isolated equilibria a natural question
in the study of relaxation. In two dimensions, the known smooth steady flows leave considerable freedom to construct nearby equilibria, which led them to
formulate the following conjecture:

\begin{deconjecture}[Drivas--Elgindi \cite{DrivasElgindi2023}]
Steady Euler flows with $C^1$ vorticity are not isolated in the
$C^1$ vorticity topology, modulo scaling and isometries of the domain.
\end{deconjecture}

Our objective in this paper is to construct isolated equilibria in the Yudovich class on every smooth annular domain whose boundary components are not concentric circles. Under suitable technical hypotheses, these isolated equilibria approximate in the Yudovich norm non-isolated smooth steady flows with the same energy. The isolated equilibria we find have discontinuous vorticity, so the smooth non-isolation conjecture of Drivas and Elgindi remains open.

The isolated equilibria we construct exhibit a strong dynamical rigidity:
a trajectory that stays sufficiently close to the constructed equilibrium
in $L^\infty$ vorticity on an open time interval must be a constant
multiple of it. Furthermore, they are strongly unstable: at fixed energy, enstrophy and circulations, distinct
nearby trajectories cannot tend weakly to these equilibria.

As a related result, we also provide an independent (and less stringent) obstruction to weak
attraction at fixed energy and
circulations which applies to every nonzero equilibrium on any bounded domain.

\subsection{Statement of the main results}

Let us first specify the meaning of isolation used here. We will say that a steady
flow $u_*$ is {\em isolated}\/ if every steady Yudovich flow $u$ with
$\|\curl u- \curl u_*\|_{L^\infty}$ sufficiently
small is a constant multiple of $u_*$ (we do not need to take a further quotient by
isometries of the domain; see Remark~\ref{rig:isometries}). Here and in what follows, we typically denote by $\omega:=\operatorname{curl} u$ the vorticity of an Euler flow with velocity~$u$, and $\omega_*$ denotes the vorticity of the reference equilibrium of velocity~$u_*$. When the Euler flow~$u$ is time-dependent, $u_0:=u(0)$ and $\omega_0:=\curl u_0$ respectively denote its initial velocity and vorticity. Unless otherwise specified, norms are henceforth taken over $\Omega$, and closeness of
vorticities is measured in $L^\infty$.

Our first theorem guarantees the existence of isolated equilibria with strong dynamical instability properties. To state it, write $Z(\omega)=\int_\Omega\omega^2\,dx$ for the enstrophy and
$\Gamma(u)$ for the vector of all boundary circulations. For a reference
flow $u_*$, let
$\mathcal M_Z(u_*)$ consist of the fields in $\mathcal Y(\Omega)$ with
the same energy, enstrophy and boundary circulations as $u_*$.

\begin{introtheorem}\label{intro:theorem-a}
There exist a smooth bounded annular domain $\Omega$,
a nonzero isolated steady Euler flow $u_*\in\mathcal Y(\Omega)$,
and a number $\delta_*>0$ with the following properties.
\begin{enumerate}
\renewcommand{\labelenumi}{\textup{(\roman{enumi})}}
\item If a Yudovich Euler flow satisfies
\[
 \sup_{t\in J}\|\omega(t)-\omega_*\|_{L^\infty}<\delta_*
\]
for some nonempty open interval $J\subset\mathbb R$, then
$u(t)=cu_*$ for all $t\in\mathbb R$ and some constant $c>0$.  \smallbreak
\item {Near the steady Euler flow $u_*$ there exist subsets of initial data that are wandering for all time. That is, any $\mathcal Y(\Omega)$-neighborhood of~$u_*$ contains {a nonempty open set}  $\mathcal U \subset \mathcal {Y}(\Omega)$ such that each trajectory $v(t)$ with initial datum $v_0\in\mathcal U$ satisfies the wandering condition $v(t)\not \in \mathcal U$ for all $t\neq 0$.} Moreover, the {corresponding} vorticities remain uniformly separated from $\omega_*$:
$$
\inf_{v_0 \in \mathcal{U}} \inf_{t \in \R} \|\omega(t)-\omega_*\|_{L^\infty} > 0.
$$

\item Every Yudovich Euler flow with initial data
$u_0\in\mathcal M_Z(u_*)\setminus\{u_*\}$ satisfying
$\|\omega_0-\omega_*\|_{L^\infty}<\delta_*$ obeys
\begin{equation}\label{intro:nonreturn}
 \omega_*\notin
 \overline{\{\omega(t):t\in\mathbb R\}}^{\,\mathrm{w}^*},
 \qquad
 \inf_{t\in\mathbb R}\|u(t)-u_*\|_{L^2}>0.
\end{equation}
\end{enumerate}
\end{introtheorem}

\begin{remark}\label{intro:immediate-departures}
Nontrivial data in
$\mathcal M_Z(u_*)$ arbitrarily close to $u_*$ do exist by Proposition~\ref{dyn:invariant-perturbations}. 
\end{remark}

For each initial datum~$u_0$ as 
above, Theorem~\ref{intro:theorem-a}\textup{(i)} forces departure from
any $L^\infty$-vorticity ball centered at $\omega_*$ of radius less than
$\delta_*$ at arbitrarily small positive and negative times. Part~\textup{(iii)} complements this immediate departure in $L^\infty$ with an
obstruction to return even in a weak topology, which would be compatible with mixing.
Note that departure in $L^\infty$ alone would not exclude such a return.
The positive separation in~\eqref{intro:nonreturn} depends on~$u_0$; it is not uniform over all nearby data. Actually, the proof of Theorem~\ref{intro:theorem-a}\textup{(iii)} shows that, for
$\|\omega_0-\omega_*\|_{L^\infty}<\delta_*$, weak-$*$ return to $\omega_*$
requires
\[
 Z(\omega_0)-Z(\omega_*)\geq\|\omega_0-\omega_*\|_{L^2}^2.
\]
At fixed energy and boundary circulations, it remains open whether such
return can occur for arbitrarily close initial data of higher enstrophy, which would also be compatible with mixing. 

{While $L^\infty$ non-return occurs during finite nonempty open time-intervals, the existence of discrete times where some $\omega(t)$ may return close to $\omega_*$ is nevertheless not excluded by (i). Part~\textup{(ii)} however shows the existence of an open set of initial data close to $\omega_*$ that remains uniformly bounded away from $\omega_*$ and that does not intersect itself for all non-zero time; it is wandering. This is the first all-time wandering result we are aware of and, in fact, the proof illustrates that non-stagnant Arnold stable steady states also admit all (non-zero) time wandering neighborhoods of initial data arbitrarily close, see Proposition~\ref{prop:wandering} and Corollary~\ref{cor:wandering-arnold-stable}. So far, Nadirashvili in \cite{nadirashvili} and later Drivas, Elgindi and Jeong in \cite{DrivasElgindiJeong2024} found neighborhoods of initial vorticities that wander on phase-space after some finite time.}

We now want to show that this kind of isolated,
highly unstable equilibria can appear arbitrarily close in the Yudovich norm to non-isolated smooth equilibria. To explain
this approximation result, suppose that
$u=\nabla^\perp\psi$ and $\Delta\psi=F(\psi)$, where the stream function
is constant on each boundary component of $\Omega$ and
$\nabla^\perp\psi=(-\partial_2\psi,\partial_1\psi)$.
Under the technical hypotheses below, we will show that $u$ is a non-isolated equilibrium. Indeed, if the linearized operator $-\Delta + F'(\psi)$ with Dirichlet boundary conditions is invertible, the implicit function theorem allows us to
replace $F$ by $F+\varepsilon$ while keeping the two boundary values of
$\psi$ fixed, and the resulting smooth stream functions converge to $\psi$
in $C^{2,\alpha}$, for $0<\alpha<1$. Since the boundary values are
distinct, the corresponding velocities are not multiples of $u$, and
rescaling them then fixes the energy. Thus $u$ is not isolated, even on
its energy surface, both in this topology of stream functions and in the
norm $\|\cdot\|_{\mathcal Y}$.

However, the next theorem shows that such a non-isolated equilibrium is a
limit of isolated ones in the Yudovich norm and on the same energy
surface. Here, the geometry of the annular domain plays a key role.

\begin{introtheorem}\label{intro:approximation}
Let $\Omega$ be a smooth bounded annular domain and let
$u=\nabla^\perp\psi$ be a steady Euler flow with
$\psi\in C^\infty(\overline\Omega)$ and $\Delta\psi=F(\psi)$,
where $F$ is {smooth} on the range of $\psi$.
Assume that $\psi$ has no critical points on $\overline\Omega$ and
that $-\Delta+F'(\psi)$ has trivial Dirichlet kernel.
There exist isolated steady Euler flows $u_n\in\mathcal Y(\Omega)$ with
\[
 E(u_n)=E(u),\qquad
 \|u_n-u\|_{\mathcal Y}\longrightarrow0,
\]
if and only if $\Omega$ is not a concentric circular annulus.
\end{introtheorem}

Heuristically, what makes the circular case different is a rigidity result of Hamel and Nadirashvili
\cite{HamelNadirashvili2023} ensuring that every $C^2$ steady flow on a
concentric circular annulus, without stagnation points up to the boundary,
is circular. However, the local version needed here must also allow competitors
with merely bounded vorticity. This case is handled in
Proposition~\ref{app:cir-local-rigidity} by a substantially shorter
argument showing that a sufficiently
small neighborhood of the smooth reference flow contains only circular
steady Yudovich flows. Of course, those cannot be isolated.

The approximating flows in the theorem have discontinuous (in fact, finitely valued) vorticity, so
their stream functions are not $C^{2,\alpha}$. {However, the lack of smoothness is not the main underlying reason behind the isolation, but rather the special features of the $L^\infty$ topology. To illustrate this, in Section \ref{app:constraints} we prove the existence of steady vorticities in $W^{k,\infty} \setminus C^k$ that are isolated in the $W^{k,\infty}$ topology, for any fixed~$k$,} {and show that the equilibria we construct can be chosen not to be isolated in the weak-$*$ $L^{\infty}$ topology. While this contrast involves only equilibria, it is morally consistent with the fact that the long-time behavior of the 2D Euler equations is strongly dependent on the choice of topology \cite{BedrossianMasmoudi2015,LinZeng2011}}.

Since any smooth non-circular annular domain admits functions $\psi$ as in Theorem~\ref{intro:approximation}, (for example, given by harmonic functions with distinct constant boundary values), we conclude that all non-circular annular domains exhibit isolated equilibria.

\begin{corollary}\label{intro:harmonic-reference}
Every smooth bounded annular domain $\Omega$ whose boundary components
are not concentric circles admits a nonzero isolated steady Euler flow
$u\in\mathcal Y(\Omega)$ with finitely valued vorticity, and with the same instability properties as in Theorem~\ref{intro:theorem-a}.
\end{corollary}

We conclude with a simple observation showing that a much weaker analogue of the failure of attraction exhibited by our isolated equilibria will in fact occur with any arbitrary nonzero steady state (and, in fact, for any reference Euler flow with precompact vorticity orbit, see Remark~\ref{app:precompact-motions}). Specifically, one can find arbitrarily close data
with the nonreturn property, keeping the energy and all boundary
circulations fixed but leaving the enstrophy free, thereby excluding attraction of a full neighborhood to any fixed
equilibrium or its scaling line. 

\begin{introtheorem}\label{intro:theorem-c}
Let $\Omega\subset\mathbb R^2$ be a bounded connected domain with
boundary of class $C^{1,1}$ and let
$u_*\in\mathcal Y(\Omega)\setminus\{0\}$ be a steady Euler flow.
There are initial data $u_0\in\mathcal Y(\Omega)$ arbitrarily close to
$u_*$, with
$E(u_0)=E(u_*)$ and $\Gamma(u_0)=\Gamma(u_*)$, whose solutions
satisfy~\eqref{intro:nonreturn}. In fact, they satisfy
\[
 \inf_{t,c\in\mathbb R}
       \|u(t)-cu_*\|_{L^2}>0.
\]
\end{introtheorem}

\begin{remark}\label{intro:regularity}
If $\partial\Omega$ is smooth, the initial data~$u_0$ in
Theorem~\ref{intro:theorem-c} can be chosen so that the difference $u_0-u_*$ is in $C^\infty(\overline\Omega;\mathbb R^2)$ and arbitrarily small in any $C^m$~norm. In particular, $u_0$ is as smooth as the initial equilibrium~$u_*$. See Section~\ref{app:constraints} for more details.
\end{remark}

\begin{remark}
    {While Theorem~\ref{intro:theorem-c} above shows the failure of asymptotic stability of $u_*$ in $L^2$ among initial data within the same energy and {boundary circulations}, no statement is made about its (in)stability.}
\end{remark}

\subsection{Previous results and related questions}
\label{intro:previous-results}

The local structure of the set of smooth steady flows has been studied
in several settings. On annular domains, Choffrut and \v{S}ver\'ak
\cite{ChoffrutSverak2012} described the neighborhood of a nondegenerate
steady state without stagnation points by means of vorticity distribution
functions, under suitable additional hypotheses. In particular, their
local uniqueness within a rearrangement class is consistent with an
infinite-dimensional family of nearby equilibria when the distribution
is allowed to vary. Constantin, Drivas and Ginsberg
\cite{ConstantinDrivasGinsberg2021} then proved rigidity results in domains with symmetries and that Arnold stable steady flows without stagnation points persist under perturbations of the domain and the vorticity law.

There are two difficulties in extending the perturbative description to
all smooth equilibria. The Dirichlet linearization may have a kernel,
and the vorticity need not be a globally defined function of the stream
function. The second issue was addressed by Elgindi, Huang, Said and Xie in
\cite{ElgindiHuangSaidXie2024}, where they proved that analytic steady flows on bounded
analytic simply connected domains are either radial or admit a global
semilinear relation. In contrast, Elgindi and Huang \cite{ElgindiHuang2026} showed that
such a relation can fail for smooth equilibria, even under additional
Morse or stability assumptions. 

The local structure of steady states also depends on the topology and the background flow. The topological dependence is clearly illustrated by the Couette flow: Lin and Zeng \cite{LinZeng2011} constructed
nonparallel steady flows arbitrarily close to Couette in vorticity
$H^s$, for $s<3/2$, and also showed that steady and traveling waves sufficiently close in
vorticity $H^s$,  with $s>3/2$, must be shear flows. The existence of nearby
nonparallel states close to Couette precludes inviscid damping, {where
the velocity converges to a nearby shear flow}, for all perturbations
in these lower-regularity neighborhoods. The sensitivity to
the reference flow is illustrated by Coti Zelati, Elgindi and Widmayer
\cite{CotiZelatiElgindiWidmayer2023}: nontrivial stationary structures
occur near Kolmogorov flow on the square torus even in analytic regularity. Dynamically, both results show these steady states cannot be attractors of whole neighborhoods of initial data, as they are not isolated from nearby stationary solutions.

The conservation laws of the Euler equation also constrain the dynamics near equilibria. For example, Arnold's energy–Casimir method~\cite{Arnold1965} and Burton's variational approach on rearrangement classes~\cite{Burton2005} used these laws to establish Lyapunov stability for appropriate steady flows. For equilibria with nonconstant
vorticity, our argument uses these invariants (specifically, a decrease of a convex Casimir at fixed
energy and circulations) to separate an orbit from a prescribed weak limit. A more direct antecedent is the work of Dolce and Drivas
\cite{DolceDrivas2022} on maximally mixed equilibria. Among other
results, they constructed vorticities arbitrarily close in $L^1$ to a
given bounded vorticity in a periodic channel whose weak rearrangement
closures, at the energy and momentum of the perturbed state, contain no shear
flows. Those perturbations need not be small in stronger vorticity
norms and do not fix the reference energy. Our perturbations are small
in $L^\infty$, preserve the reference energy and all circulations, and
apply on arbitrary bounded domains. However, the conclusion excludes the
reference scaling line rather than an entire class of shear equilibria.

Instability in strong norms and weak relaxation are closely
related, but distinct, questions. Drivas, Elgindi and Jeong
\cite{DrivasElgindiJeong2024} showed that sustained twisting near stable
flows leads to filamentation and strong-norm instability, while the recent
work of Alazard and Said \cite{AlazardSaid2026} establishes generic
creation of small scales for the Euler equation on the two-torus.
Such results concern the evolution of spatial regularity; they do not
identify all weak limit motions. The relaxation conjectures of
\v{S}ver\'ak and Shnirelman \cite{Sverak2011,Shnirelman2013}, in the
form stated in \cite[Conjecture~4.5]{Elgindi2026}, predict that generic
vorticity orbits are not precompact in $L^2$, while weak limiting states
generate $L^2$-precompact orbits. The latter need not be stationary; our nonreturn statements only concern one prescribed steady state within
this broader relaxation problem. 

The behavior predicted by these relaxation conjectures has been
established in perturbative regimes through inviscid damping. In periodic channels, these results require sufficiently small perturbations of strictly monotone shear flows in suitable Gevrey classes. Bedrossian and Masmoudi \cite{BedrossianMasmoudi2015} proved nonlinear inviscid damping near Couette flow in the periodic strip. Subsequently, Ionescu and Jia
\cite{IonescuJia2023} and Masmoudi and Zhao \cite{MasmoudiZhao2024}
established nonlinear damping near monotone shear flows in a finite
channel under suitable regularity and spectral hypotheses. In these
results, mixing produces decay of the non-shear component of the
velocity, while the limiting shear depends on the initial datum,
consistently with Theorem~\ref{intro:theorem-c}. This convergence can
coexist with growth of strong vorticity norms. Guo and Luo
\cite{guo2026asymptotic} also established asymptotic stability at
Yudovich regularity near Couette flow in the nonperiodic infinite
channel $\mathbb R\times[0,1]$, for compactly supported vorticity
perturbations sufficiently small in $L^\infty$. In this setting,
transport to spatial infinity yields local $L^\infty$ velocity decay and weak
$L^p$ convergence to zero of the vorticity perturbation for every
$1<p<\infty$, a mechanism that is unavailable in bounded domains.

A related nonperturbative relaxation result was
obtained by Elgindi, Murray and Said \cite{ElgindiMuraySaid2025}
for scale-invariant solutions in $\mathbb R^2$. Under $m$-fold symmetry with $m\geq4$, they proved that solutions whose
initial angular vorticity has bounded variation relax to piecewise
constant angular profiles with finitely many jumps.

{We finish the discussion noting that the first author, Kepplinger and Peralta-Salas \cite{EncisoKepplingerPeraltaSalas2025} constructed isolated steady flows in three dimensions through a mechanism based on robust chaotic dynamics of Anosov type.} On suitable closed Riemannian three-manifolds, the combination of Anosov particle dynamics and a simple curl eigenvalue gives smooth
equilibria that are isolated,
modulo scaling, in the $C^1$~velocity topology.
On the flat torus $\mathbb T^3$, the same argument confines nearby analytic equilibria
to a six-dimensional space, a weaker conclusion than isolation modulo
scaling. This idea is not applicable to the two-dimensional case: any nonzero steady Euler flow has its stream function as a nonconstant continuous first integral, so its flow lines cannot exhibit Anosov dynamics. Instead, as we will see below, our argument uses a finitely valued vorticity and a rough
vorticity law as a source of rigidity.  

\subsection{Ideas of the proofs}
Let us first sketch the main idea of the construction. To obtain an isolated equilibrium $u_*=\nabla^\perp\psi_*$, we seek a
stream function whose level curves foliate the annular domain and whose vorticity
satisfies $\omega_*=\Delta\psi_*=F_*(\psi_*)$. The key idea is to choose a Borel
function $F_*$ that takes finitely many values, with at least two distinct
values occurring on sets of positive measure in every open interval of the stream-function values.

The first step is to show that any nearby steady flow preserves \(\omega_*\). To see this, let $v$ be another steady Euler flow, with vorticity
$\widetilde\omega=\operatorname{curl}v$, and let $\delta$ be half the
smallest gap between the values of $F_*$. If
$\|\widetilde\omega-\omega_*\|_{L^\infty}$ is less than $\delta$, then rounding
$\widetilde\omega$ to the nearest of these values recovers
$\omega_*$ almost everywhere as a function of~$\widetilde\omega$. Since $v$ transports its own vorticity,
it must also transport $\omega_*$. 

Suppose now that a streamline of $v$ crosses
a level curve of $\psi_*$ transversely. Near the crossing, $s:=\psi_*$
and a local stream function for $v$ form a system of local coordinates.
In these coordinates, transport of $\omega_*$ gives
$\partial_s F_*(s)=0$ in the sense of distributions, contrary to the
choice of $F_*$. Thus $v$ is tangent to the level curves of $\psi_*$, and
incompressibility gives
$v=A(\psi_*)u_*$ for some absolutely continuous scalar function $A$. To show that $A$ is constant, we arrange that the speed $|u_*|^2$ varies
along almost every level curve. On $\{\psi_*=s\}$, the vorticity of $v$
is
\[
 \widetilde\omega=A(s)F_*(s)+A'(s)|u_*|^2.
\]
For almost every $s$ with $A(s)\ne0$, stationarity makes
$\widetilde\omega$ constant along the curve. Since $|u_*|^2$ varies
there, $A'(s)=0$, so $v=cu_*$ for some constant~$c$. 

For a time-dependent solution that
stays in the same vorticity neighborhood on an open time interval,
rounding again recovers the fixed function $\omega_*$ at every time, and similar arguments give the
rigidity assertion of Theorem~\ref{intro:theorem-a}. 

To solve the equation $\Delta\psi_*=F_*(\psi_*)$, we
start from a smooth nondegenerate equilibrium for which the oscillation
of the speed has a positive lower bound over all level curves. We
approximate its smooth vorticity law uniformly by finitely valued laws
that oscillate on every interval. Although the usual implicit function
argument no longer applies, the foliation gives a uniform area bound
$|\{x:\psi(x)\in B\}|\leq C|B|$ for measurable sets $B$ and nearby
stream functions $\psi$ with the same boundary data. By approximation
with continuous functions, this bound makes $\psi\mapsto F_*(\psi)$
continuous into $L^p$ for finite $p$. The inverse of the smooth
Dirichlet linearization then permits a Schauder fixed-point argument.
The resulting solution is close enough in $C^1$ to preserve the
foliation and the variation of speed along its leaves, as shown in
Proposition~\ref{con:rough-approximation}.

The wandering proof of Theorem~\ref{intro:theorem-a} relies on the essential twisting ideas already present in \cite{nadirashvili}: if the smooth reference non-degenerate equilibrium is stable in $L^2$ and has distinct boundary travel times, the differential rotation on $\partial\Omega$ is still present for nearby solutions, in particular those with finitely valued vorticity laws. Consider one such steady state $\omega_*$ and let $Q$ denote the set of attained vorticity values. Since it is finite, $Q$ and $Q+\eps$ are disjoint for $\eps>0$ small enough so that $\tilde\omega_\eps = \omega_* + \eps\mathbf 1_\mathcal{N}$ is {an initial vorticity} with finitely many values, for some $\mathcal{N}\subset \Omega$ with $\overline{\mathcal N}\cap \partial\Omega\neq \emptyset$. For $0<\rho < \frac12 \dist(Q,Q+\eps)$, the Borel measurable map defined by $\mathsf b(s) =\mathbf{1}_{\lbrace \dist(s,Q+\eps) \leq \rho\rbrace}$ distinguishes $\mathcal{N}$, namely $\mathsf b(\curl v(t))  = \mathbf{1}_{{\Xi}_t^{{v_0}}(\mathcal N)}$ almost everywhere, for any $v_0$ sufficiently close to $\tilde u_\eps ={u_*+\eps\nabla^\perp z}$, {with $\Delta z=\mathbf1_{\mathcal N}$ in $\Omega$ and $z|_{\partial\Omega}=0$, and} where $\Xi_t^{v_0}$ denotes the flow-map of $\mathcal S_tv_0$. Any finite time return of $v(t)$ would yield $\mathcal N = \Xi_t^{v_0}(\mathcal N)$, which in turn contradicts the differential rotation rates of $\partial\Omega$ under $\Xi_t^{v_0}$.

The two weak nonreturn arguments use different conserved quantities. For the third statement in Theorem~\ref{intro:theorem-a}, let $Q_*$ be the finite essential range of the
constructed equilibrium $\omega_*$. Initial closeness gives
$\operatorname{dist}(\omega_0,Q_*)=|\omega_0-\omega_*|$ almost
everywhere. Transport preserves the $L^2$ norm of the left-hand side.
Since $\omega_*(x)\in Q_*$ almost everywhere, it follows that
\[
 \|\omega(t)-\omega_*\|_{L^2}\geq
 \|\omega_0-\omega_*\|_{L^2}\qquad\text{for all }t.
\]
When the enstrophies agree, weak-$*$ convergence to $\omega_*$ would
imply strong $L^2$ convergence, because the $L^2$ norms agree. The
inequality therefore excludes weak return whenever $\omega_0\ne\omega_*$.

To prove Theorem~\ref{intro:theorem-c}, for an arbitrary equilibrium with nonconstant vorticity, we instead
construct small perturbations at fixed energy and circulations that
decrease a convex Casimir $C_H(f)=\int_\Omega H(f)\,dx$. Once the inner
circulations are fixed, the remaining circulation and energy constraints
require the vorticity perturbation to be orthogonal, to first order, to
$1$ and $\psi_*$. For $H(s)=s^2$ and $H(s)=s^4$, the Casimir gradients
are proportional to $\omega_*$ and $\omega_*^3$. These cannot both be
affine functions of $\psi_*$, since the cube of a nonconstant affine
function is not affine on an interval. Thus one of the two Casimirs has
a decreasing admissible direction. A second-order correction restores
the energy exactly and preserves the strict decrease. Conservation and
weak lower semicontinuity then give
{$C_H(\omega_\infty)\leq C_H(\omega_0)<C_H(\omega_*)$ for every weak-$*$ limit
$\omega_\infty$ of the orbit.} Since $H$ is even, this excludes both $\omega_*$
and $-\omega_*$. Energy conservation then separates the velocity orbit
from the whole scaling line in $L^2$. For constant-vorticity equilibria,
we obtain the same separation by constructing nearby stationary states
with nonconstant vorticity and the same energy and circulations.

For the approximation theorem, Theorem \ref{intro:approximation}, the smooth reference equilibrium need
not have the required variation of speed along its streamlines. On a
nonconcentric annular domain there are only finitely many level curves of
constant speed, but removing one by perturbing the vorticity law may
simply create another at a neighboring level. We therefore need a
response that cannot be reproduced by changing the level value.
For interior curves, oscillatory perturbations of the smooth law and
the corresponding elliptic response provide such a direction unless
the curvature is constant. That exceptional case would force concentric
circular boundaries by unique continuation. The boundary curves require a separate argument with fixed Dirichlet data. The resulting
independence of the two variations excludes constant speed on both the original and the neighboring levels, and allows us to prepare the geometry on the fixed domain. We then apply the Borel construction just described. A diagonal choice of approximations gives convergence in the Yudovich norm, and a final constant rescaling restores the prescribed energy. It will be clear from the proofs that sufficiently high finite regularity suffices wherever we assume smoothness, but throughout the paper we have chosen to present the results in the smooth setting to simplify the statements.
 
\subsection{Organization of the paper} Section~\ref{sec:rigidity} establishes
the rigidity mechanism for steady flows and for trajectories confined to
a small vorticity neighborhood. Section~\ref{sec:construction} constructs
the isolated equilibria used in Theorem~\ref{intro:theorem-a}, whose
remaining dynamical assertions are proved in Section~\ref{sec:dynamics}. Section~\ref{sec:approximation}
proves Theorem~\ref{intro:approximation} and
Corollary~\ref{intro:harmonic-reference}, and
Section~\ref{sec:no-attraction} proves
Theorem~\ref{intro:theorem-c}.
Finally, Section~\ref{app:constraints} justifies the constrained perturbations
in Remarks~\ref{intro:immediate-departures} and~\ref{intro:regularity}, {and the higher regularity claims following Theorem \ref{intro:approximation}.}
Appendix~\ref{app:coefficient-calculation} records the coefficient calculations
used in the approximation argument.

\section{Isolation Lemma}
\label{sec:rigidity}

We first prove the rigidity statement that will be used in both
constructions. The hypotheses concern the vorticity law and the
geometry of the level curves, and impose no restriction on the
circulations of a competing Euler flow.

\begin{lemma}\label{rig:rigidity}
Let $\Omega$ be a bounded annular domain with boundary of class $C^{1,1}$,
and let $\psi\in W^{2,p}(\Omega)$ for some $p>2$. Suppose that $\psi$ takes
the constant values $a<b$ on the two boundary components,
$|\nabla\psi|\geq m>0$ on $\overline\Omega$, and
\[
 \Delta\psi=F(\psi)\quad\text{a.e. in }\Omega.
\]
Assume that there is a $C^1$ diffeomorphism
\[
 X:[a,b]\times\mathbb S^1\longrightarrow\overline\Omega,
 \qquad \psi(X(s,\theta))=s,
\]
and that the following properties hold:
\begin{enumerate}
\item The function $F$ is Borel measurable, has finite range, and is
not essentially constant on any nonempty open subinterval of $(a,b)$. \smallbreak
\item The function $|\nabla\psi|^2$ is nonconstant on $\{\psi=s\}$
for a.e. $s\in(a,b)$.
\end{enumerate}
Let $Q=\{q_1,\ldots,q_N\}$ be the essential range of $F$ on $(a,b)$,
and set
\begin{equation}\label{rig:gap}
 \delta:=\frac12\min_{i\ne j}|q_i-q_j|>0,
 \qquad u:=\nabla^\perp\psi,\qquad \omega:=F(\psi).
\end{equation}
Then, $u$ is a steady Euler flow and, if a Yudovich solution $v(t)$ satisfies
\[
 \|\operatorname{curl}v(t)-\omega\|_{L^\infty}<\delta
 \quad\text{for every }t\text{ in a nonempty open interval }J,
\]
then $v(t)=cu$ for every $t\in\mathbb R$, for a single constant $c>0$.
If in addition $E(v(0))=E(u)$, then $c=1$.
\end{lemma}

The essential range contains at least two values because $F$ is not
essentially constant. Stationary isolation follows by applying the
lemma to a time-independent solution. In the proof, we establish
this special case first and use it after showing that a trajectory
confined to the stated neighborhood is stationary. The hypothesis
does not require a uniform margin below $\delta$ on $J$, so it implies
the version stated in
Theorem~\ref{intro:theorem-a}.

\begin{remark}\label{rig:isometries}
The lemma imposes no restriction on the isometries of $\Omega$.
If $g(x)=Ox+b$ preserves $\Omega$, with $O$ orthogonal, then
$u_g(x)=O u(g^{-1}x)$ is a steady Euler flow with $E(u_g)=E(u)$ and
\[
 \curl u_g=(\det O)\,\omega\circ g^{-1}.
\]
If $\|\curl u_g-\omega\|_{L^\infty}<\delta$, applying the lemma to the
constant trajectory $v(t)=u_g$ gives $u_g=u$ by equality of energy.
Thus every image under a domain isometry either coincides with $u$
or lies outside this vorticity neighborhood. No quotient by
isometries is required for the isolation asserted here.
\end{remark}

\subsection{Stream functions and level coordinates}

Throughout this section, $\Omega$ is as in Lemma~\ref{rig:rigidity}.
We record the facts needed to work with an arbitrary bounded-vorticity
competitor. Weak incompressibility and impermeability are understood as
\[
 \int_\Omega v\cdot\nabla\zeta\,dx=0
 \qquad\text{for all }\zeta\in H^1(\Omega).
\]

\begin{lemma}\label{rig:streamfunction}
Every $v\in\mathcal Y(\Omega)$ has a stream function
$\phi\in W^{2,r}(\Omega)$ for every $1<r<\infty$, constant on each
boundary component, such that $v=\nabla^\perp\phi$. In particular,
$v\in W^{1,r}(\Omega)$ and $\phi\in C^{1,\alpha}(\overline\Omega)$
for every $0<\alpha<1$. The boundary constants are not prescribed.
Moreover, for every $2\leq r<\infty$,
\begin{equation}\label{rig:div-curl-estimate}
 \|v\|_{W^{1,r}(\Omega)}
 \leq C_\Omega r\|v\|_{\mathcal Y},
\end{equation}
where $C_\Omega$ is independent of $r$.
Vorticity and boundary circulations together determine $v$ uniquely.
\end{lemma}

\begin{proof}
The stream-function representation and uniqueness from vorticity and
boundary circulations follow from the classical planar Hodge decomposition. See for instance \cite[Section~2.1]{LequeurreMunnier2020} and
\cite[Section~4]{Sueur2011}.
Estimate~\eqref{rig:div-curl-estimate} follows from the quantitative
Dirichlet estimate on $C^{1,1}$ domains in
\cite[Section~1]{LacaveMiotWang2014}, applied to the zero-Dirichlet
part and, after subtracting a smooth lift of its boundary values,
to the harmonic part. The coefficient of the latter is bounded by
$C_\Omega\|v\|_{L^2}$ by $L^2$ orthogonality.
The remaining regularity follows from Dirichlet regularity and
Sobolev embedding.
\end{proof}

\begin{remark}\label{rig:stationary-transport}
Let $v\in\mathcal Y(\Omega)$ be a stationary weak solution to the Euler
equations and set $\widetilde\omega:=\operatorname{curl}v\in L^\infty(\Omega)$.
Then, it follows that
\begin{equation}\label{rig:transport}
 \operatorname{div}(\widetilde\omega v)=0, \quad \textup{that is} \quad
 \int_\Omega\widetilde\omega v\cdot\nabla\zeta\,dx=0
 \quad\text{for all }\zeta\in C_c^1(\Omega).
\end{equation}
Indeed, Lemma~\ref{rig:streamfunction} gives
$v\in W^{1,2}(\Omega)\cap L^\infty(\Omega)$, so local mollification
justifies taking the curl of the stationary momentum equation in the
form $\operatorname{div}(\widetilde\omega v)=0$.
Since $\widetilde\omega v\in L^\infty(\Omega)$, the testing identity
extends from $C_c^\infty(\Omega)$ to $C_c^1(\Omega)$ by approximation.
\end{remark}

We will repeatedly use the following elementary consequence of
distributional differentiation.

\begin{lemma}\label{rig:zero-derivative}
If $W\in L^\infty(I_0\times J_0)$, where $I_0,J_0$ are intervals,
and $\partial_sW=0$ in the sense of distributions, then $W(s,t)=w(t)$ a.e. for a bounded measurable function $w$. Likewise, if
$H\in L^\infty((a,b)\times\mathbb S^1)$ and
$\partial_\theta H=0$ in the sense of distributions, then $H(s,\theta)=h(s)$ a.e.
for a bounded measurable function $h$.
\end{lemma}

\begin{proof}
By \cite[Theorem~A.1]{BoykoKunzingerPopovych2021}, the distribution
induced by $W$ is independent of $s$. Its representative is
$w(t)=\int_{I_0}W(s,t)\chi(s)\,ds$, where
$0\leq\chi\in C_c^\infty(I_0)$ and $\int_{I_0}\chi=1$,
so $w$ is bounded and measurable.
For the second assertion, apply the same theorem after identifying
$\mathbb S^1$ minus one point with $(0,2\pi)$. Averaging in $\theta$
yields a bounded measurable $h$.
\end{proof}

Write $I=(a,b)$ and orient the angular coordinate so that
\[
 D(s,\theta):=\det(\partial_sX,\partial_\theta X)>0.
\]
Since $X$ is a $C^1$ diffeomorphism, $D$ is continuous and strictly
positive on the compact cylinder $[a,b]\times\mathbb S^1$.
It therefore has a positive minimum and a finite maximum, so both
$D$ and $D^{-1}$ are bounded there. Change of variables therefore
yields, for any measurable $B\subset I$,
\begin{equation}\label{rig:null-transfer}
 c|B|\leq|\{x:\psi(x)\in B\}|\leq C|B|.
\end{equation}
In particular, a null set of values has a null spatial preimage,
and each value in $Q$ occurs on a spatial set of positive measure.
Differentiating $\psi\circ X=s$ also gives
\begin{equation}\label{rig:chart-identities}
 \partial_\theta X=D\,u(X),
 \qquad
 (u\cdot\nabla z)(X)=D^{-1}\partial_\theta(z\circ X).
\end{equation}
Indeed, $\partial_\theta X$ is parallel to $u(X)$, and taking its
determinant with $\partial_sX$ identifies the coefficient as $D$.

These observations also verify that the reference field is steady.
The Lipschitz function $G(s)=\int_a^s F(r)\,dr$ satisfies
$\nabla G(\psi)=F(\psi)\nabla\psi$ a.e. in $\Omega$. The exceptional
values in the one-dimensional chain rule have null spatial preimage
by~\eqref{rig:null-transfer}. Since
\[
 (u\cdot\nabla)u
 =\nabla\frac{|\nabla\psi|^2}{2}-F(\psi)\nabla\psi,
\]
the pressure associated with $\psi$ can be chosen as
$P[\psi]:=G(\psi)-|\nabla\psi|^2/2$. The boundary
condition follows from the constant boundary values of $\psi$.

\subsection{Rigidity of stationary neighbors}

\begin{lemma}\label{rig:stationary-rigidity}
Under the hypotheses of Lemma~\ref{rig:rigidity}, every steady Euler flow
$v\in\mathcal Y(\Omega)$ satisfying
$\|\operatorname{curl}v-\omega\|_{L^\infty}<\delta$ is of the form $v=cu$
for a constant $c>0$.
\end{lemma}

\begin{proof}
Let $v=\nabla^\perp\phi$ be a {nearby} steady Euler flow and choose
$\rho<\delta$ such that
\begin{equation}\label{rig:pointwise-closeness}
 |\widetilde\omega-F(\psi)|\leq\rho
 \quad\text{a.e. in }\Omega,
 \qquad \widetilde\omega:=\operatorname{curl}v.
\end{equation}
We first show that the two stream functions cannot cross
transversely. We argue by contradiction and assume that, at some
interior point $x_0$,
$K(x_0)=\det(\nabla\psi(x_0),\nabla\phi(x_0))\ne0$,
where $K:=\det(\nabla\psi,\nabla\phi)$. The inverse function theorem
provides local coordinates
\[
 T(x)=(\psi(x),\phi(x))=(s,r)
\]
on a neighborhood of $x_0$ mapping onto a rectangle $I_0\times J_0$, with
$I_0\subset I$, on which $K$ has constant sign and is bounded away
from zero. The chart and its inverse preserve null sets. For a
smooth compactly supported test $\xi$ in this rectangle, use
$\zeta=\xi\circ T$ in~\eqref{rig:transport}. The identities
\[
 v\cdot\nabla\psi=-K,\qquad
 v\cdot\nabla\phi=0,\qquad dx=|K|^{-1}\,ds\,dr
\]
then show that
$v\cdot \nabla \zeta(x) = (\partial_s \xi)\circ T(x)  (v\cdot \nabla\psi)(x) = -K(x) (\partial_s \xi)\circ T(x)$ and thus
\[
 0=\int_{T^{-1}(I_0\times J_0)} \widetilde\omega v\cdot\nabla \zeta\, d x = -\int_{T^{-1}(I_0\times J_0)}\widetilde\omega K ((\partial_s \xi)\circ T) d x = -\text{sgn}(K) \int_{I_0\times J_0}(\widetilde\omega\circ T^{-1})
                  \partial_s\xi\,ds\,dr.
\]
Lemma~\ref{rig:zero-derivative} implies that
$\widetilde\omega\circ T^{-1}=w(r)$ a.e. in $I_0\times J_0$.
Consequently, for a.e. $r\in J_0$,
\[
 |w(r)-F(s)|\leq\rho\quad\text{for a.e. }s\in I_0.
\]
At least two distinct values $q_i,q_j$ occur on positive-measure
subsets of $I_0$, since $F$ is finite valued and not essentially
constant there. The preceding inequalities then imply
$2\delta\leq|q_i-q_j|\leq2\rho<2\delta$, giving a contradiction.
We have therefore proved that
\begin{equation}\label{rig:parallel}
 \det(\nabla\psi,\nabla\phi)=0
 \quad\text{on }\overline\Omega.
\end{equation}

Having this identity at hand, we infer from~\eqref{rig:chart-identities} that 
\begin{align*}
    \partial_\theta(\phi\circ X)(s,\theta)= D(s,\theta)(u\cdot \nabla \phi)(X(s,\theta)) = 0.
\end{align*}
Since the level curves are connected, there is a $C^1$ function $h$ on $[a,b]$ such that
\begin{equation}\label{rig:level-multiplier}
 \phi=h(\psi),\qquad v=A(\psi)u,\qquad A=h'\in C([a,b]).
\end{equation}
Before differentiating this identity, observe that
\[
 B(x):=\frac{\nabla\phi\cdot\nabla\psi}{|\nabla\psi|^2}
       =A(\psi(x))
\]
belongs to $W^{1,p}(\Omega)$. This follows from the Sobolev product
rule, boundedness of the gradients, and the lower bound for
$|\nabla\psi|$. Pullback by the $C^1$ diffeomorphism $X$ preserves
$W^{1,p}$-regularity, and averaging $B\circ X=A(s)$ in $\theta$ shows that
$A\in W^{1,p}(I)$. Smooth approximation in this one-dimensional
space, together with~\eqref{rig:null-transfer}, now justifies
$\nabla(A(\psi))=A'(\psi)\nabla\psi$. Thus, it follows that
\begin{equation}\label{rig:curl-multiplier}
 \widetilde\omega=A(\psi)F(\psi)+A'(\psi)|\nabla\psi|^2.
\end{equation}
In fact $A\in W^{1,\infty}(I)$, because this identity and
$|\nabla\psi|\geq m$ give
\[
 |A'(\psi)|\leq
 m^{-2}\bigl(\|\widetilde\omega\|_{L^\infty}
       +\|A\|_{L^\infty(I)}\|F\|_{L^\infty(I)}\bigr),
\]
and~\eqref{rig:null-transfer} transfers the estimate to a.e.
$s\in I$.

Finally, let $S=|\nabla\psi|^2$, and write
$\widehat\omega=\widetilde\omega\circ X$ and $\widehat S=S\circ X$. Also, let $\xi \in C^\infty_0(I\times \mathbb{S}^1)$ be a smooth test function and $\zeta:= \xi \circ X^{-1}$. Then, $\xi = \zeta\circ X$ and
\begin{align*}
    \partial_\theta\xi(s,\theta) = \nabla\zeta(X)\cdot \partial_\theta X = D(s,\theta) (u\cdot \nabla \zeta)(X(s,\theta)).
\end{align*}
Hence, testing~\eqref{rig:transport} with $\zeta$,
\begin{align*}
0 & = \int_\Omega \widetilde\omega v\cdot \nabla \zeta d x = \int_\Omega \widetilde\omega(x) A(\psi(x)) (u \cdot \nabla \zeta)(x) d x \\
& = \int_\Omega \widetilde\omega(x) A(\psi(x)) \frac{(\partial_\theta\xi)(X^{-1}(x))}{(D\circ X^{-1})(x)} d x = 
 \int_I\int_{\mathbb S^1}A(s)\widehat\omega(s,\theta)
                      \partial_\theta\xi(s,\theta)\,d\theta\,d s,
\end{align*}
since $d x = D(s,\theta)\,ds \, d\theta$. Hence, $A\widehat\omega$ is independent of $\theta$ for a.e. $s\in I$. On the other hand,~\eqref{rig:curl-multiplier} reads
\[
 A(s)\widehat\omega(s,\theta)
 =A(s)^2F(s)+A(s)A'(s)\widehat S(s,\theta).
\]
The function $\widehat S(s,\cdot)$ is continuous and nonconstant
for a.e. $s\in I$. It cannot be a.e. constant on
any such circle. Therefore $AA'=0$ a.e. in $I$.
The absolutely continuous function $A^2$ is constant, and continuity
of $A$ on the connected interval implies that $A\equiv c\in\mathbb R$.

We are just missing to show that $c>0$. Note that each value of $Q$
occurs on a positive-measure spatial set,
so $\|\omega\|_{L^\infty}=\max_j|q_j|\geq \delta$. If $c\leq0$, then
\[
 \|\widetilde\omega-\omega\|_{L^\infty}
 =|c-1|\|\omega\|_{L^\infty}\geq \delta,
\]
contrary to the hypothesis. Thus $c>0$. Since $E(u)>0$, equality
of energies forces $c=1$.
\end{proof}

\subsection{Rigidity of trajectories}

The linear growth in $r$ in~\eqref{rig:div-curl-estimate} yields
log-Lipschitz regularity of a Yudovich velocity, locally uniformly
in time. Thus the usual uniqueness and Lagrangian transport
properties hold on the $C^{1,1}$ annular domain considered here
\cite{Yudovich1963,LacaveMiotWang2014}.
For $s,t\in\mathbb R$, let $\Xi_{t,s}(x)$ denote the position at time $t$ of the particle
located at $x\in\Omega$ at time $s$. Its trajectory is the unique
absolutely continuous solution to
\[
 \Xi_{t,s}(x)
 =x+\int_s^t v(\tau,\Xi_{\tau,s}(x))\,d\tau.
\]
Thus the second index is the initial time. These maps preserve area
and satisfy, for all $r,s,t\in\mathbb R$,
\[
 \Xi_{s,s}=\mathrm{id},\qquad
 \Xi_{t,s}\circ\Xi_{s,r}=\Xi_{t,r},\qquad
 \Xi_{t,s}^{-1}=\Xi_{s,t}.
\]
Vorticity transport takes the form
\begin{equation}\label{rig:lagrangian-transport}
 \widetilde\omega(t,\Xi_{t,s}(x))
 =\widetilde\omega(s,x)\quad\text{for a.e. }x\in\Omega,
 \qquad \widetilde\omega=\operatorname{curl}v.
\end{equation}
Every bounded Borel function $b(\widetilde\omega)$ is transported as
well. This assertion does not require a chain rule for $b$: compose
\eqref{rig:lagrangian-transport} with $b$, change variables by the
flow, and differentiate the smooth test function along the
absolutely continuous trajectories. The result is
\begin{equation}\label{rig:borel-transport}
 \partial_t b(\widetilde\omega)
       +\operatorname{div}\bigl(vb(\widetilde\omega)\bigr)=0
\end{equation}
for all bounded Borel measurable functions $b:\mathbb{R}\rightarrow \mathbb{R}$.
We also use the continuity $v\in C(\mathbb R;L^2)$ supplied by
Yudovich theory. Conservation of the vorticity bound and energy,
together with~\eqref{rig:div-curl-estimate}, gives uniform bounds
on the spatial $W^{1,r}$ norms for each finite $r$.

\begin{proof}[Proof of Lemma~\ref{rig:rigidity}]
Choose a Borel map $\mathfrak q:\mathbb R\to Q$ assigning to each
number a nearest element of $Q$, resolving ties in a fixed way.
The separation of $Q$ and the assumed closeness imply
\[
 \mathfrak q(\widetilde\omega(t,x))=\omega(x)
 \quad\text{for all }t\in J,\quad\text{for a.e. }x\in\Omega.
\]
Since $\mathfrak{q} \circ \widetilde \omega=\omega$ satisfies \eqref{rig:borel-transport} and $\omega$ is independent of time,~\eqref{rig:borel-transport} gives
$\operatorname{div}(v\omega)=0$ on $J\times\Omega$.
Testing with $\vartheta(t)\zeta(x)$, where $\vartheta\in C_c^\infty(J)$,
and using $v\in C(J;L^2)$ \leavevmode{\ignorespaces and $\omega = F(\psi)$}, we obtain
\begin{equation}\label{rig:frozen-labels}
 \int_\Omega F(\psi)v(t)\cdot\nabla\zeta\,dx=0,
 \qquad\text{for all }\zeta\in C_c^1(\Omega),
\end{equation}
for every $t\in J$.

Fix $t\in J$ and let $\phi_t$ be a stream function for $v(t)$.
If $\det(\nabla\psi,\nabla\phi_t)$ were nonzero at some point,
the coordinates $(s,r)=(\psi,\phi_t)$ used in Lemma~\ref{rig:stationary-rigidity} would transform~\eqref{rig:frozen-labels} into
\[
\int_{I_0\times J_0}
F(s)\,\partial_s\xi(s,r)\,ds\,dr=0
\qquad
\text{for all }\xi\in C_c^\infty(I_0\times J_0).
\]
By Lemma~\ref{rig:zero-derivative}, $F$ would be essentially constant
on $I_0$, giving a contradiction.
Thus $v(t)\cdot\nabla\psi=0$ on $\overline\Omega$ for every $t\in J$.
At each time the fields are continuous in space, so the argument
leading to~\eqref{rig:level-multiplier} gives
\begin{equation}\label{rig:time-multiplier}
 v(t,x)=A(t,\psi(x))u(x),\qquad A(t,\cdot)\in C([a,b]).
\end{equation}

We next use the first equation in~\eqref{intro:euler}, applied to $v$,
to show that $A(t,\cdot)$ is independent of time. Set
\[
 C_\psi(s):=\int_{\mathbb S^1}|u(X(s,\theta))|^2
                               D(s,\theta)\,d\theta.
\]
This is a positive continuous weight bounded away from zero, and
\[
 \|v(t)-v(t')\|_{L^2}^2
 =\int_I|A(t,s)-A(t',s)|^2C_\psi(s)\,ds
\]
together with $v\in C(\mathbb{R}; L^2)$ shows that $A$ is continuous in time with values in the weighted $L^2$
space. Moreover, it is bounded on compact time intervals, since $|u|\geq m$
and the velocities are uniformly bounded there.

Let $\beta\in C_c^\infty(I)$ and $\vartheta\in C_c^\infty(J)$ be fixed
but arbitrary, and set $w_\beta:=\beta(\psi)u$. This field belongs to $W^{1,p}$, is divergence
free, and is compactly supported in $\Omega$, since $\beta$ vanishes
near the boundary values $a$ and $b$ of $\psi$.
We test the first equation in~\eqref{intro:euler}, applied to $v$,
against $\vartheta(t)w_\beta(x)$. This is justified by extending $w_\beta$
by zero and mollifying in space: the smooth divergence-free spatial
approximations converge in $W^{1,p}$ and have compact support in
$\Omega$ for sufficiently small mollification radii.
The boundedness of $v$ on $\operatorname{supp}\vartheta\times\Omega$
allows passage to the limit, and the pressure term vanishes. Thus, it follows that
\begin{equation}\label{rig:tested-euler}
 \int_J\int_\Omega\Big[\vartheta'(t)v\cdot w_\beta
   +\vartheta(t)v\cdot\bigl((v\cdot\nabla)w_\beta\bigr)\Big]dx\,dt=0.
\end{equation}
The spatial integral involving $(v\cdot\nabla)w_\beta$ is zero.
Indeed, setting $S:=|u|^2$ and using~\eqref{rig:time-multiplier},
$u\cdot\nabla\psi=0$, and~\eqref{rig:chart-identities}, we obtain
\begin{align*}
 \int_\Omega v\cdot\bigl((v\cdot\nabla)w_\beta\bigr)\,dx
 &=\frac12\int_\Omega A(t,\psi)^2\beta(\psi)
                                  u\cdot\nabla S\,dx\\
 &=\frac12\int_I A(t,s)^2\beta(s)
       \left(\int_{\mathbb S^1}\partial_\theta(S\circ X)(s,\theta) \,d\theta\right)ds =0.
\end{align*}
Here $S\circ X\in W^{1,p}$, so for a.e. $s\in I$ its angular slice
is periodic and the integral of its derivative is zero.
Identity~\eqref{rig:tested-euler} therefore reduces to
\[
 \int_J\vartheta'(t)
       \left(\int_I A(t,s)\beta(s)C_\psi(s)\,ds\right)dt=0.
\]
By Lemma~\ref{rig:zero-derivative}, the inner expression is constant
a.e. in $J$, and hence everywhere by continuity in $t$.
Since $\beta\in C_c^\infty(I)$ is arbitrary and $C_\psi>0$, we have
$A(t,\cdot)=A(t_0,\cdot)$ a.e. in $I$ for any $t,t_0\in J$.
Continuity in $s$ extends this equality to $[a,b]$, so $v(t)=V$ on
$J$ for a fixed field $V$. Let $p$ be the pressure associated with
$v$, choose $\vartheta_0\in C_c^\infty(J)$ with $\int_J\vartheta_0(t)\,dt=1$,
and define $p_0:=\int_J\vartheta_0(t)p(t,\cdot)\,dt$ in the sense of distributions.
Multiplying the first equation in~\eqref{intro:euler} by $\vartheta_0(t)$
and integrating over $J$, we obtain $(V\cdot\nabla)V+\nabla p_0=0$
in the sense of distributions in $\Omega$, since $\partial_t v=0$. Thus $V$ is a stationary weak Euler
flow, and Lemma~\ref{rig:stationary-rigidity} gives $V=cu$, with $c>0$.
Yudovich uniqueness then gives $v(t)=cu$ for all $t\in\mathbb R$.
Finally, if $E(v(0))=E(u)$, then $c^2=1$, and hence $c=1$.
\end{proof}

\section{Construction of isolated equilibria}
\label{sec:construction}

We develop a construction of isolated steady Euler flows that will be
used in Theorems~\ref{intro:theorem-a} and~\ref{intro:approximation}.
We first construct finitely valued Borel approximations of a continuous
vorticity law and establish the composition estimates needed to solve
the corresponding semilinear equation. We then give a general
perturbation result and apply Lemma~\ref{rig:rigidity} to obtain
isolation. Finally, we construct a smooth reference Euler flow to which this
result applies, producing the equilibrium in
Theorem~\ref{intro:theorem-a}.

We begin by constructing the finitely valued Borel vorticity law on an arbitrary
compact interval $I=[a,b]$, $a<b$. By \cite[Chapter~2, Exercise~8]{Rudin1987} (see also \cite{Rudin1983}),
there is a Borel set $B\subset I$ such that $0<|B\cap U|<|U|$ for every
nonempty open interval $U\subset(a,b)$.  

Fix such a set $B$ once and for all. Given a continuous function
$F_0:I\to\R$ and $\eta>0$, define
\begin{equation}\label{con:profile}
 F_\eta(s)=\eta\left\lfloor\frac{F_0(s)}\eta\right\rfloor
                 +\eta\mathbf1_B(s),\qquad s\in I.
\end{equation}
The floor records the interval of length $\eta$ containing
$F_0(s)$; the second term chooses one of its two endpoints. It is
this choice on sets of positive measure at every scale that enforces
the rigidity hypothesis.

\begin{lemma}\label{con:profile-properties}
The function $F_\eta$ is Borel, has finite range in $\eta\Z$,
and satisfies
\begin{equation}\label{con:profile-error}
 |F_\eta(s)-F_0(s)|\leq\eta,\qquad\textup{ for all }s\in I.
\end{equation}
On every nonempty open subinterval of $(a,b)$, two values separated
by exactly $\eta$ occur on sets of positive measure. In particular,
$F_\eta$ is not essentially constant on any such interval.
\end{lemma}
\begin{proof}
The sets $\{s\in I:k\eta\leq F_0(s)<(k+1)\eta\}$ are Borel, and only
finitely many integers $k$ occur because $F_0$ is bounded. This
proves measurability and finiteness of the range. If
$F_0(s)=k\eta+\tau$ with $0\leq\tau<\eta$, then
$F_\eta(s)-F_0(s)$ is either $-\tau$ or $\eta-\tau$, proving
\eqref{con:profile-error}.

Let $U\subset(a,b)$ be a nonempty open interval. If $F_0$ is
constant there, its floor is constant there as well. Otherwise, since
$F_0$ is continuous, $F_0(U)$ contains an interval of positive length,
so there is $s_0\in U$ with $F_0(s_0)\notin\eta\Z$. By continuity, a nonempty open
$V\subset U$ satisfies $k\eta<F_0(V)<(k+1)\eta$ for some $k$.
In either case, the floor is a fixed integer $k$ on a nonempty open
$V\subset U$. On this interval, it follows that
\[
 F_\eta=(k+1)\eta\quad\hbox{on }B\cap V,
 \qquad F_\eta=k\eta\quad\hbox{on }V\setminus B.
\]
Both sets have positive measure by the choice of $B$.
\end{proof}

We next establish the continuity of composition with an $L^p$ vorticity
law under small $C^1$ perturbations of the stream function.

{Let $\psi_0 \in C^2(\overline{\Omega})$ have constant boundary values $a<b$
on a bounded annular domain with boundary of class $C^{1,1}$.} Assume that there is a $C^1$ diffeomorphism
\begin{equation}\label{con:parametrization}
 \Phi:I\times\Sone\longrightarrow\overline\Omega,
 \qquad \psi_0(\Phi(s,\theta))=s.
\end{equation}
Then, for all $s\in \psi_0(\overline{\Omega})$, the closed curve $\lbrace \Phi(s,\theta) : \theta \in \mathbb{S}^1 \rbrace$ corresponds to the streamline $\lbrace \psi_0 =s \rbrace$. We orient $\theta$ so that $J_0:=\det D\Phi>0$ and write
\[
 \mathcal{X}=\{v\in C^1(\overline\Omega):v|_{\partial\Omega}=0\},
 \qquad \|v\|_{\mathcal{X}}=\|v\|_{L^\infty}+\|\nabla v\|_{L^\infty}.
\]
Also, we choose $r_0>0$ sufficiently small that
$r_0\|\partial_s\Phi\|_{L^\infty(I\times\Sone)}\leq1/2$. For $\|v\|_{\mathcal{X}}\leq r_0$,
put $\psi_v=\psi_0+v$. Then
\begin{equation}\label{con:monotonicity}
 \frac12\leq\partial_s(\psi_v\circ\Phi)
 =1+\nabla v(\Phi)\cdot\partial_s\Phi\leq\frac32.
\end{equation}
For fixed $\theta$, this is an increasing map of $I$ onto itself,
since its endpoint values are $a,b$. Its inverse in the first
variable, denoted $\sigma_v(t,\theta)$, is jointly $C^1$. Therefore, the map
\begin{equation}\label{con:perturbed-chart}
 X_v(t,\theta)=\Phi(\sigma_v(t,\theta),\theta),
 \qquad \psi_v(X_v(t,\theta))=t,
\end{equation}
is a $C^1$ diffeomorphism from $I\times\Sone$ onto $\overline\Omega$.
Its Jacobian and inverse Jacobian are bounded uniformly on the ball. In particular, the perturbed functions have no critical points
and have connected level curves. As before, the closed curve $\lbrace X_v(t,\theta) : \theta \in \mathbb{S}^1 \rbrace$ corresponds to the streamline $\lbrace \psi_v =t \rbrace$, for all $t\in \psi_v(\overline{\Omega})$.

\begin{lemma}\label{con:pushforward}
There is a constant $C_\rho$, independent of $v\in\mathcal{X}$
satisfying $\|v\|_{\mathcal{X}}\leq r_0$, such that
\begin{equation}\label{con:pushforward-bound}
 \int_\Omega q(\psi_v(x))\dd x
 \leq C_\rho\int_a^b q(t)\dd t
\end{equation}
for every nonnegative measurable $q$. Consequently,
$\|g\circ\psi_v\|_{L^p(\Omega)}\leq
C_\rho^{1/p}\|g\|_{L^p(I)}$ for $1\leq p<\infty$, and composition
is well defined on almost-everywhere equivalence classes.
\end{lemma}
\begin{proof}
We change variables first by $x=\Phi(s,\theta)$ and then by
$t=\psi_v(\Phi(s,\theta))$. Tonelli's theorem yields
\[
 \int_\Omega q(\psi_v(x))\dd x = \int_a^b \int_{\mathbb{S}^1} q(\psi_v(\Phi(s,\theta))) J_0(s,\theta) \dd \theta \dd s=\int_a^b q(t)\rho_v(t)\dd t,
\]
with 
\begin{align*}
   \rho_v(t)=\int_{\mathbb S^1}
 \frac{J_0(\sigma_v(t,\theta),\theta)}
 {\partial_s(\psi_v\circ\Phi)(\sigma_v(t,\theta),\theta)}\dd\theta. 
\end{align*}
By \eqref{con:monotonicity}, one may take
$C_\rho=4\pi\max J_0$. Taking $q$ to be the indicator of a null
set shows that its preimage is null. This also extends the argument
from Borel to Lebesgue measurable functions and proves the assertion
about equivalence classes.
\end{proof}

\begin{lemma}\label{con:composition}
For every $F\in L^p(I)$, $1\leq p<\infty$, the map
$v\mapsto F(\psi_0+v)$ is continuous from the closed $\mathcal{X}$-ball of
radius $r_0$ into $L^p(\Omega)$.
\end{lemma}
\begin{proof}
Let $v_n\to v$ in $\mathcal{X}$, and approximate $F$ in $L^p(I)$ by continuous
functions $G_j\in C([a,b])$, since $(\psi_0+v_n)(\overline{\Omega})\subset [a,b]$ due to \eqref{con:monotonicity} and $v|_{\partial\Omega}=0$. Lemma~\ref{con:pushforward} gives
\[
 \begin{split}
 \|F(\psi_0+v_n)-F(\psi_0+v)\|_{L^p}
 &\leq2C_\rho^{1/p}\|F-G_j\|_{L^p(I)}+\|G_j(\psi_0+v_n)-G_j(\psi_0+v)\|_{L^p}.
 \end{split}
\]
For fixed $j$, the last term tends to zero by uniform continuity
of $G_j$. Hence,
\begin{align*}
    \lim_{n\rightarrow \infty} \|F(\psi_0+v_n)-F(\psi_0+v)\|_{L^p} \leq 2C_\rho^{1/p}\Vert F - G_j \Vert_{L^p(I)} \leq \varepsilon
\end{align*}
for all $j\geq j_0$ large enough.
This proves continuity without
any pointwise continuity assumption on $F$.
\end{proof}

We now solve the semilinear equation for a bounded Borel law sufficiently
close to the smooth reference law, and obtain isolation under the
additional assumptions on its essential values and variation.

\begin{proposition}\label{con:rough-approximation}
{Let $\psi_0\in C^2(\overline\Omega)$ satisfy
$\Delta\psi_0=F_0(\psi_0)$ on an annular domain with boundary of class $C^{1,1}$, with constant
boundary values $a<b$, $F_0\in C^2(I)$, and a $C^1$ parametrization
$\Phi$ as in \eqref{con:parametrization}.} Assume that
$-\Delta+F_0'(\psi_0)$ has trivial Dirichlet kernel and that
\[
\mathbf{o}:=\inf_{s\in I}\osc_{\{\psi_0=s\}}|\nabla\psi_0|^2>0.
\]
Set $u_0:=\nabla^\perp\psi_0$. There is $\eta_0>0$ such that, for every
$0<\eta<\eta_0$ and every bounded Borel function $F:I\to\R$ satisfying
$\|F-F_0\|_{L^\infty(I)}\leq\eta$, there is
$\psi_F\in\bigcap_{1<p<\infty}W^{2,p}(\Omega)$ with the following
properties:
\begin{enumerate}
\item\label{con:rough-steady} The function $\psi_F$ satisfies
\[
 \Delta\psi_F=F(\psi_F),\qquad
 \psi_F|_{\partial\Omega}=\psi_0|_{\partial\Omega},
\]
and the velocity $u_F=\nabla^{\perp}\psi_F\in\mathcal{Y}(\Omega)$
is a steady Euler flow.
\smallbreak
\item\label{con:rough-geometry} There is a $C^1$ diffeomorphism
\begin{equation}\label{con:rough-parametrization}
 \Phi_F:I\times\Sone\longrightarrow\overline\Omega,
 \qquad \psi_F(\Phi_F(s,\theta))=s.
\end{equation}
The level curves of $\psi_F$ are connected and satisfy
\[
 \inf_{s\in I}\osc_{\{\psi_F=s\}}|\nabla\psi_F|^2\geq\mathbf{o}/2.
\]
\smallbreak
\item\label{con:rough-bounds} For every $1<p<\infty$,
\begin{equation}\label{con:rough-estimates}
 \|\psi_F-\psi_0\|_{W^{2,p}}\leq C_p\eta,
 \qquad
 \|u_F-u_0\|_{\mathcal Y}\leq C\eta,
\end{equation}
with constants independent of $F$ and $\eta$.
\smallbreak
\item\label{con:rough-isolation} If, in addition, $F$ has finite essential range $Q$ and is not
essentially constant on any nonempty open subinterval of $(a,b)$, set
\[
 2\delta:=\min\{|q-q'|:q,q'\in Q,\ q\ne q'\}>0.
\]
Then every steady $v\in\mathcal Y(\Omega)$ with
$\|\curl v-\curl u_F\|_{L^\infty}<\delta$ is a positive constant
multiple of $u_F$. In particular, for $F=F_\eta$ defined by
\eqref{con:profile}, the resulting Euler flow is isolated modulo positive
scaling with vorticity radius $\eta/2$.  
\end{enumerate}
\end{proposition}

\begin{proof}
We first prove~\ref{con:rough-steady}. Fix $p_0>2$. Triviality of the
Dirichlet kernel and elliptic Fredholm theory give a bounded inverse
\[
 L^{-1}:L^{p_0}(\Omega)\longrightarrow
 W^{2,p_0}(\Omega)\cap W^{1,p_0}_0(\Omega),
 \qquad L=\Delta-F_0'(\psi_0).
\]
Since $p_0>2$, the latter space embeds compactly into $\mathcal{X}$. On the closed $\mathcal{X}$-ball of
radius $r_0$ define
\begin{equation}\label{con:fixed-point-map}
 T_F(v)=L^{-1}
 \bigl[F(\psi_0+v)-F_0(\psi_0)-F_0'(\psi_0)v\bigr].
\end{equation}
We shall find a fixed point $v=T_F(v)$. Note first that all arguments of the functions $F$ and $F_0$ above belong to $I$ by \eqref{con:monotonicity}. Since $\Vert F - F_0 \Vert_{L^\infty(I)}\leq \eta$, there exists $N\subset I$ with $|N|=0$ such that $|F(t) - F_0(t)|\leq  \eta$ for all $t\in I\setminus N$. Appealing to Lemma~\ref{con:pushforward}, we see that
\begin{align}\label{eq:L_infty_diff}
    |\lbrace x\in \Omega : \psi_0(x) + v(x) \in N \rbrace| = \int_\Omega \mathbf{1}_N(\psi_0(x) + v(x)) \dd x \leq C_\rho \int_a^b \mathbf{1}_N(t) \dd t = 0, 
\end{align}
from which we deduce that 
\begin{align}\label{eq:comp_diff}
    \Vert F\circ(\psi_0+v) - F_0\circ(\psi_0+v)\Vert_{L^\infty}\leq \eta,
\end{align}
for all $v\in \mathcal{X}$ with $\|v\|_{\mathcal X} \leq r_0$.
With $C_0:=\tfrac12\|F_0''\|_{L^\infty(I)}$, Taylor's theorem
therefore gives
\[
 |F(\psi_0+v)-F_0(\psi_0)-F_0'(\psi_0)v|
 \leq\eta+C_0|v|^2
 \quad\text{a.e. in }\Omega.
\]
Hence, for a constant $A>0$ independent of $F$ and $\eta$,
\begin{equation}\label{con:map-bound}
 \|T_F(v)\|_{\mathcal{X}}\leq A(\eta+C_0\|v\|_{\mathcal{X}}^2).
\end{equation}
For sufficiently small $\eta$, the closed $\mathcal{X}$-ball of radius
$2A\eta$ is contained in the radius-$r_0$ ball and is mapped
into itself: it suffices that $2A\eta\leq r_0$ and
$4C_0A^2\eta\leq1$.

Lemma~\ref{con:composition} makes $T_F$ continuous into
$W^{2,p_0}$, and the same estimates bound its image there. The
compact embedding into $\mathcal{X}$ therefore permits the application of
Schauder's fixed-point theorem \cite{GT,Schauder}. Let $v_F$
be a fixed point and set $\psi_F=\psi_0+v_F$. Thus, $v_F$ is the
scalar perturbation of the stream function. Cancelling
the linear terms in \eqref{con:fixed-point-map} yields the exact
equation $\Delta\psi_F=F(\psi_F)$, and
$\|v_F\|_{\mathcal{X}}\leq2A\eta$.

To complete the proof of~\ref{con:rough-steady}, we show that the {divergence free and tangent to the boundary} velocity field $u_F:= \nabla^\perp\psi_F\in W^{1,p_0}(\Omega)$ is a steady weak Euler solution. We define
\[
 \mathcal P_F(s)=\int_a^sF(t)\dd t,
 \qquad p_F=\mathcal P_F(\psi_F)
                    -\tfrac12|\nabla\psi_F|^2.
\]
Approximating $F$ in $L^{p_0}(I)$ by smooth functions and their
primitives uniformly on $I$, Lemma~\ref{con:pushforward} allows
passage to the limit in their chain rules, giving
\[
 \nabla(\mathcal P_F(\psi_F))
 =F(\psi_F)\nabla\psi_F.
\]
The Sobolev product rule and our convention for $\nabla^\perp$ give
\[
 (u_F\cdot\nabla)u_F
 =\nabla\bigl(\tfrac12|u_F|^2\bigr)
       -F(\psi_F)\nabla\psi_F
 =-\nabla p_F.
\]
so that the momentum equation \eqref{intro:euler} is satisfied.

We next prove~\ref{con:rough-geometry}. The parametrization
$\Phi_F:=X_{v_F}$ from \eqref{con:perturbed-chart} gives the claimed level geometry. To compare speeds on corresponding
level curves, observe that its inverse coordinate satisfies
\begin{align*}
    \psi_F(\Phi(s,\theta)) = s + v_F(\Phi(s,\theta))
\end{align*}
for all $(s,\theta)\in[a,b]\times \mathbb{S}^1$. For $s=\sigma_{v_F}(t,\theta)$ and $X_{v_F}(t,\theta) = \Phi(\sigma_{v_F}(t,\theta),\theta)$ we have
\begin{align*}
    t = \psi_F(X_{v_F}(t,\theta)) = \sigma_{v_F}(t,\theta) + v_F(X_{v_F}(t,\theta))
\end{align*}
for all $(t,\theta)\in [a,b]\times \mathbb{S}^1$. Hence,
\[
 \sigma_{v_F}(t,\theta)-t
 =-v_F(\Phi(\sigma_{v_F}(t,\theta),\theta)),
\]
so $\|\sigma_{v_F}-t\|_{L^\infty(I\times\Sone)}\leq\|v_F\|_{L^\infty}$.
Since $\psi_0\in C^2(\overline{\Omega})$ and $\Phi\in C^1([a,b]\times\mathbb{S}^1)$, we have that $|\nabla\psi_0|^2\circ\Phi$ is Lipschitz in $s$, and it follows that
\[
 \sup_{t,\theta}\left|
 |\nabla\psi_F|^2(X_{v_F}(t,\theta))
       -|\nabla\psi_0|^2(\Phi(t,\theta))\right|
 \leq C\|v_F\|_{\mathcal{X}}\leq C\eta.
\]
Oscillations of two functions differ by at most twice their uniform
distance. Decreasing $\eta_0$ once more proves the bound $\mathbf{o}/2$,
including on the boundary level curves.

To prove~\ref{con:rough-bounds}, note that 
\[
 \Delta v_F=F(\psi_F)-F_0(\psi_0),\qquad
 v_F|_{\partial\Omega}=0,
\]
with
\begin{align*}
    \Vert F(\psi_F) - F_0(\psi_0) \Vert_{L^\infty} \leq \Vert F(\psi_F) - F_0(\psi_F) \Vert_{L^\infty} + \Vert F_0(\psi_F) - F_0(\psi_0) \Vert_{L^\infty} &\leq \eta + \Vert F'_0\Vert_{L^\infty(I)} \Vert v_F\Vert_{L^\infty}\leq C\eta,
\end{align*}
due to \eqref{eq:comp_diff}.
Dirichlet estimates for this equation give the first bound in
\eqref{con:rough-estimates} for every finite $p>1$. Combining the case
$p=2$ with the preceding vorticity bound gives the estimate in
$\mathcal Y$. The smallness condition on $\eta$ was chosen only once,
using $p_0$.

To prove~\ref{con:rough-isolation}, suppose that $F$ has finite
essential range $Q$ and is not essentially constant on any nonempty open subinterval of $(a,b)$.
By changing $F$ on the null set where its values lie outside $Q$,
which leaves $F(\psi_F)$ unchanged a.e. in $\Omega$ by
Lemma~\ref{con:pushforward}, we may assume that $F$ has finite range.
All the hypotheses of Lemma~\ref{rig:rigidity} are then satisfied,
so that lemma gives isolation with radius $\delta$.
For the law $F_\eta$ defined by \eqref{con:profile},
Lemma~\ref{con:profile-properties} gives both the required essential
variation and $\delta\geq\eta/2$, proving the last assertion.
\end{proof}

As a by-product of the construction, we show that the steady states $\omega_F$ constructed in Proposition \ref{con:rough-approximation}  have wandering sets of initial data arbitrarily close.

\begin{proposition}\label{prop:wandering}
{Let $\psi_0\in C^2(\overline\Omega)$ satisfy
$\Delta\psi_0=F_0(\psi_0)$ on an annular domain with boundary of class $C^{1,1}$, with constant
boundary values $a<b$, $F_0\in C^2(I)$, and a $C^1$ parametrization
$\Phi$ as in \eqref{con:parametrization}.} Assume that
$-\Delta+F_0'(\psi_0)$ has trivial Dirichlet kernel. 
Set $u_0:=\nabla^\perp\psi_0$ and assume that it has distinct boundary travel times and satisfies
\begin{align}\label{eq:u0-L2-stability}
    \forall\ \tilde\delta>0, \, \exists\ \tilde\varepsilon>0 \, : \, \Vert v_0 - u_0 \Vert_{\mathcal Y}< \tilde\varepsilon\, \Longrightarrow\, \sup_{t\geq 0} \Vert \mathcal S_tv_0 - u_0 \Vert_{L^2}< \tilde\delta,
\end{align}
for $v_0\in \mathcal Y(\Omega)$, and where $\mathcal S_t$ denotes the Yudovich solution operator. 
Then, there are $C>0$ and $\eta_0>0$ such that, for every
$0<\eta<\eta_0$ and every finitely valued bounded Borel function $F:I\to\R$ satisfying
$\|F-F_0\|_{L^\infty(I)}\leq\eta$, there is $\psi_F\in\bigcap_{1<p<\infty}W^{2,p}(\Omega)$ such that 
\begin{align*}
    \begin{cases}
        \Delta\psi_F = F(\psi_F), & \text { in }\Omega, \\
        \psi_F = \psi_0, & \text{ on }\partial\Omega.
    \end{cases}
\end{align*}
The field $u_F:=\nabla^\perp\psi_F$ is a steady Euler flow and there holds $\Vert u_F - u_0 \Vert_{\mathcal Y}\leq C\eta$. Additionally, for all $r>0$ there is a non-empty set $\mathcal U$ open in $\mathcal Y(\Omega)$ such that $\mathcal U \subset \lbrace v \in \mathcal Y(\Omega) : \Vert v - u_F \Vert_{\mathcal Y} < r \rbrace$ and
\begin{align*}
    \mathcal S_t (\mathcal U) \cap \mathcal U = \emptyset, \quad \textup{for all } t \neq 0.
\end{align*}
Furthermore, for $\omega_F :=\curl u_F$, there holds
\begin{align*}
    \inf_{v_0\in \mathcal U}\inf_{t\in \R} \Vert \curl \mathcal S_t v_0 - \omega_F \Vert_{L^\infty} > 0.
\end{align*}
In particular, $u_0$ also admits wandering sets arbitrarily close to it and uniformly separated in vorticity from $\curl u_0$.
\end{proposition}

\begin{proof}

We follow the idea of \cite{nadirashvili}, albeit with technical differences due to the non-circular geometry of $\Omega$. We first introduce coordinates that identify the annular domain $\overline\Omega$ with $I\times \Sone$. We use the $C^1$ diffeomorphism $\Phi:=\Phi_{F_0}:I\times\Sone \rightarrow \overline{\Omega}$, which is such that $\psi_0(\Phi(s,\theta))=s$. For $c\in \lbrace a, b\rbrace$, let $\Gamma_c = \lbrace \psi_0 = c \rbrace$. We note that $\gamma_c(\theta):\Sone\rightarrow\Gamma_c$ defined by $\gamma_c(\theta)= \Phi(c,\theta)$ parametrizes $\Gamma_c$. Since $\Phi$ is a $C^1$ diffeomorphism, we observe that $\gamma'_c(\theta) = \partial_\theta\Phi(c,\theta)$ is non-zero and we define $\tau_c(\theta):= \frac{\gamma'_c(\theta)}{|\gamma'_c(\theta)|}$, for $\theta\in \Sone$. Since $u_0$ is tangent to $\partial\Omega$, both $\frac{u_0(\gamma_c(\theta))}{|u_0(\gamma_c(\theta))|}$ and $\tau_c(\theta)$ are unit vectors tangent to $\Gamma_c$. Therefore, there exists $\sigma_c\in \lbrace 1, -1\rbrace$ such that
\begin{align*}
    u_0(\gamma_c(\theta)) = \sigma_c {|u_0(\gamma_c(\theta))|}\tau_c(\theta), \quad \textup{for all } \theta \in \Sone.
\end{align*}
 
The travel time $\mathsf{T}_c(u_0)$ of particles on $\Gamma_c$ under the velocity field $u_0$ is defined as $\mathsf{T}_c(u_0) :=\int_{\Gamma_c}\frac{\dd \ell}{|u_0|}$ for $c\in \lbrace a, b \rbrace$. Using the parametrization given by $\gamma_c$, we write
\begin{align*}
    \mathsf T_c(u_0) = \int_0^{2\pi} \frac{|\gamma'_c(\theta)|}{|u_0(\gamma_c(\theta))|} \dd \theta
\end{align*}
and we define $\nu_c := 2\pi \sigma_c /\mathsf{T}_c(u_0)$, which denotes the signed normalized angle variation per unit time of particles along $\Gamma_c$. Since $\mathsf T_a(u_0) \neq \mathsf T_b(u_0)$, we have $|\nu_a|\neq |\nu_b|$ and thus $|\nu_a - \nu_b|>0$, i.e. there is a differential rate of rotation of the two boundary components.

To measure the actual angle displacement, since $\theta\mapsto \gamma_c(\theta)$ may not be an arc-length parametrization of $\Gamma_c$ and the speed along $\Gamma_c$ need not be constant, we define
\begin{align*}
    \vartheta_c(\theta) := \frac{2\pi}{\mathsf{T}_c(u_0)} \int_0^\theta \frac{|\gamma_c'(\xi)|}{|u_0(\gamma_c(\xi))|}\dd \xi, \quad \textup{for all } \theta \in \R.
\end{align*}
We observe that the integrand is positive and $2\pi$-periodic. Thus, $\vartheta_c(\theta+2\pi k) = \vartheta_c(\theta) +2\pi k$ for all $k\in \Z$.

{Now, let $v_0\in \mathcal Y(\Omega)$ and $v(t):= \mathcal{S}_t v_0$, denote its associated flow map by $\Xi_t^{v_0}$, and define $\gamma_c^{v}(t) = \Xi_t^{v_0}(\gamma_c(\alpha))$, for some $\alpha\in \Sone$.} Then, $\gamma_c^{v}(t) = \gamma_c(\Theta_c^{v}(t))$ for some locally absolutely continuous $\Theta_c^{v}:\R\to \R$ with $\Theta_c^v(0)=\alpha$. By definition of $\gamma_c^{v}$ and the chain rule, almost everywhere in time, 
\begin{align*}
    v(t,\gamma_c^{v}(t)) = \frac{\dd}{\dd t} \gamma_c^{v}(t) = \gamma'_c(\Theta_c^{v}(t)) \dot \Theta_c^{v}(t).
\end{align*}
Taking the scalar product with $\tau_c(\Theta_c^{v}(t))$ yields
\begin{align*}
    \dot \Theta_c^v(t) = \frac{v(t, \gamma_c^v(t))\cdot \tau_c(\Theta_c^v(t))}{|\gamma'_c(\Theta_c^v(t))|}.
\end{align*}
Therefore, the variation of angle displacement along $\Theta_c^v$ is
\begin{align*}
    \frac{\dd}{\dd t} \vartheta_c(\Theta_c^v (t) ) = \frac{2\pi}{\mathsf{T}_c(u_0)}\frac{v(t, \gamma_c^v(t))\cdot \tau_c(\Theta_c^v(t))}{|u_0(\gamma_c^v(t))|}.
\end{align*}
In particular, for $K_c:= \frac{2\pi}{\mathsf{T}_c(u_0)}\frac{1}{\min_{\Gamma_c}|u_0|}$,
\begin{align*}
    \left| \frac{\dd}{\dd t} \vartheta_c(\Theta_c^v (t) ) - \nu_c \right| &= \frac{2\pi}{\mathsf{T}_c(u_0)}\frac{1}{|u_0(\gamma_c^v(t))|} \left| \left(v(t,\gamma_c^v(t)) - u_0(\gamma_c^v(t)) \right) \cdot \tau_c(\Theta_c^v(t)) \right| \leq K_c \Vert v(t) - u_0 \Vert_{{L^\infty}}.
\end{align*}
Hence, if $\Vert v(t) - u_0 { \Vert_{L^\infty} \leq \varepsilon_0}$ for all $t> 0$ and $K_c {\varepsilon_0} \leq \frac{|\nu_b - \nu_a|}{4}$ then, for the angular lift displacement
\begin{align*}
    \tilde\vartheta_c(t) := \vartheta_c(\Theta_c^v(t)) - \vartheta_c(\Theta_c^v(0)),
\end{align*}
we obtain $|\tilde\vartheta_c(t) - \nu_c t| \leq \frac{|\nu_b - \nu_a|}{4}t$ for all $t\geq0$ and $c\in \lbrace a, b \rbrace$. As a result, $|\tilde\vartheta_b(t) - \tilde\vartheta_a(t) | \geq \frac{|\nu_b - \nu_a|}{2}t$ for all $t\geq 0$. In other words, there is at least linear separation of angular lift displacements between the inner and outer boundaries, a testament to the twisting of $u_0$ inherited by $v(t)$ for all positive time provided $v(t)$ remains {$L^\infty$} close to $u_0$.  

{Before fixing $F$ and $r$, choose the tolerance $\varepsilon_0$ above, fix $p>2$, and restrict to initial data with $\Vert v_0-u_0\Vert_{\mathcal Y}\leq1$. By conservation of energy and vorticity,}
the div-curl estimate yields
\begin{align*}
    \Vert v(t) - u_0 \Vert_{W^{1,p}}\leq C_{\Omega,p}\left( \Vert \curl v{(t)} - {\curl u}_0\Vert_{L^p} + \Vert v(t) - u_0 \Vert_{L^2} \right) \leq C,
\end{align*}
for some $C>0$. Then, by the compact embedding $W^{1,p}(\Omega) \hookrightarrow C(\overline{\Omega})$, for all $\varepsilon_{{0}}>0$ there exists $\tilde\delta>0$ such that $ \Vert v(t) - u_0 \Vert_{C(\overline{\Omega})}\leq \varepsilon_{{0}}$ if $\Vert v(t) - u_0 \Vert_{L^2(\Omega)}< \tilde\delta$. {For this $\tilde\delta$, choose $0<\tilde{\varepsilon}<1$ in \eqref{eq:u0-L2-stability}. The existence and approximation arguments in the proof of Proposition~\ref{con:rough-approximation}, which do not use the speed-oscillation or essential-variation assumptions, give the stated $\psi_F$ and bound $\Vert u_F-u_0\Vert_{\mathcal Y}\leq C_0\eta$. Reduce their uniform threshold $\eta_0$ so that $C_0\eta_0<\tilde{\varepsilon}/3$.}

{Now fix $0<\eta<\eta_0$, an admissible $F$, and $r>0$.} In what follows, we find a suitable subset of initial data that will wander and stay uniformly away from $\curl u_F$. To do so, define $\mathcal R:= (a,b)\times(\alpha,\beta)\subset I \times \Sone$. Let $\mathcal N = \Phi(\mathcal R)\subset \Omega$ and let $\omega_F:= \curl u_F$ where $u_F$ is {any such} steady state. {Set $Q:=\operatorname{ess\,ran}\omega_F$}. For $0<\varepsilon$ {sufficiently small} we set $\tilde\omega_\varepsilon = \omega_F + \varepsilon \mathbf 1_{\mathcal N}$ and we observe that $Q$ and $Q+\varepsilon$ have no intersection{, since $Q$ is finite; this also holds if $Q$ is a singleton}. {Let $z$ solve $\Delta z=\mathbf1_{\mathcal N}$ in $\Omega$, with $z=0$ on $\partial\Omega$.} Let $\tilde v_\varepsilon={u_F+\varepsilon\nabla^\perp z}$ be the velocity field associated to $\tilde\omega_\varepsilon$, with $\tilde v_\varepsilon\rightarrow u_F$ in $\mathcal{Y}$. {More precisely, $\Vert\tilde v_\varepsilon-u_F\Vert_{\mathcal Y}\leq C_{\mathcal N}\varepsilon$; choose $\varepsilon$ so that $C_{\mathcal N}\varepsilon<\min\lbrace r/2,\tilde{\varepsilon}/3\rbrace$.} We shall see that a neighborhood of $\tilde v_\varepsilon$ does not return to itself. 

Let $\varepsilon_1>0$ be chosen below and define $\mathcal U_{\varepsilon_1} := \lbrace v\in \mathcal Y(\Omega) : \Vert v - \tilde v_\eps \Vert_{\mathcal Y} < \varepsilon_1 \rbrace$. If $M(\omega) = \frac{1}{|\Omega|}\int_\Omega \omega(x) \dd x$ denotes the mean of the vorticity $\omega$, we note that $M(\tilde \omega_\eps ) - M(\omega_F) = \eps |\mathcal N | / |\Omega| >0$. Since $|M(\curl v_0) - M(\tilde\omega_\eps)|\leq \Vert \curl v_0 - \tilde\omega_\eps\Vert_{L^\infty}$, we have that
\begin{align*}
    \frac{\eps}{2|\Omega|}|\mathcal N | \leq |M(\curl v_0) - M(\omega_F)| = |M(\curl v(t)) - M(\omega_F)|  \leq \Vert \curl v(t) - \omega_F\Vert_{L^\infty}
\end{align*}
for all $v_0\in \mathcal U_{\varepsilon_1}$, for $\varepsilon_1>0$ small enough, as $M(\curl v_0)$ is conserved along the Euler solution. Hence, $\curl v(t)$ remains uniformly bounded away from $\omega_F$ for all $t\in \R$.

To take advantage of the twisting of $u_0$, we shall see that $ \Vert v(t) - u_0 {\Vert_{L^{\infty}}\leq \varepsilon_0}$ for all $t\geq 0$ and all $v_0\in \mathcal U_{\varepsilon_1}$. First, we note that for $\varepsilon_1>0$ small enough, we have $\Vert v - u_F \Vert_{\mathcal Y} < r $ for all $v\in \mathcal U_{\varepsilon_1}$. {By the choice of $\tilde{\varepsilon}$, it suffices to check that $\Vert v_0-u_0\Vert_{\mathcal Y}<\tilde{\varepsilon}$.} Additionally, we record
\begin{align}\label{eq:v0-to-u0-ineq}
    \Vert v_0 - u_0\Vert_{\mathcal Y}\leq \Vert v_0 - \tilde v_\eps \Vert_{\mathcal Y} + \Vert \tilde v_\eps - u_F \Vert_{\mathcal Y} + \Vert u_F - u_0 \Vert_{\mathcal Y} \leq \varepsilon_1 + C_\mathcal N\eps + C\eta.
\end{align}
{Further choosing $\varepsilon_1<\tilde{\varepsilon}/3$, the bounds $C_0\eta<\tilde{\varepsilon}/3$ and $C_{\mathcal N}\varepsilon<\tilde{\varepsilon}/3$ place \eqref{eq:v0-to-u0-ineq} within the stability neighborhood of $u_0$. Hence,} we obtain
\begin{align*}
    |\tilde\vartheta_b(t) - \tilde\vartheta_a(t) | \geq \frac{|\nu_b - \nu_a|}{2}\,t
\end{align*}
for all $t\geq 0$ and all $v_0\in \mathcal U_{{\varepsilon_1}}$.

Once a uniform differential rotation rate of the boundaries is ensured for all Euler evolutions of initial data in $\mathcal U_{\varepsilon_1}$, we next define a Borel map that distinguishes $\mathcal N$. For $0<\rho < \frac12 \dist(Q,Q+\eps)$, we define the Borel measurable map $\mathsf b:\R \to \lbrace 0,1 \rbrace$ by
\begin{align*}
    \mathsf b(s) =\mathbf{1}_{\lbrace \dist(s,Q+\eps) \leq \rho\rbrace}.
\end{align*}
Therefore, $\mathsf b(\curl v) = \mathbf{1}_\mathcal N$ for all $\Vert \curl v - \tilde \omega_\eps \Vert_{L^\infty}\leq \rho$. Further assuming ${\varepsilon_1} < \rho$, for $v_0\in \mathcal U_{{\varepsilon_1}}$
\begin{align*}
    \mathsf b(\curl v(t)) = \mathsf b (\curl v_0) \circ ({\Xi}_t^{v_{{0}}})^{-1} = \mathbf{1}_{{\Xi}_t^{v_{{0}}}(\mathcal N)}
\end{align*}
almost everywhere in $\Omega$, since $\curl v(t)$ is transported by ${\Xi}_t^{v_{{0}}}$. In particular, if $v(t)\in \mathcal U_{{\varepsilon_1}}$ for some $t\mathrel{{>}}0$ then $\mathcal N = {\Xi}_t^{v_{{0}}}(\mathcal N)$ almost everywhere. Since we are interested in the behavior of the boundaries, which have zero measure, we next show $\mathcal N = {\Xi}_t^{v_{{0}}}(\mathcal N)$: as $\mathcal N$ coincides with the interior of its closure relative to $\Omega$, the same is true for ${\Xi}_t^{v_{{0}}} (\mathcal N)$ because homeomorphisms map interiors to interiors and closure to closures and ${\Xi}_t^{v_{{0}}}$ is a homeomorphism. Additionally, any two open subsets that agree almost everywhere in $\Omega$ have the same relative closure (otherwise, if a point of the first set lies outside the closure of the second, so does a small ball centered at the point and properly included in the first set; the difference of sets has positive measure). Consequently, $\mathcal N = {\Xi}_t^{v_{{0}}}(\mathcal N)$.

For the remainder of the proof, fix $v_0\in \mathcal U_{\varepsilon_1}$ such that $v(t)\in \mathcal{U}_{\varepsilon_1}$ for some $t\neq 0$. The particle map in the cylinder $I\times \Sone$ of the evolution of $v_0$ is given by $H_t := \Phi^{-1}\circ \Xi_t^{v_0} \circ \Phi:I \times \Sone \rightarrow I\times \Sone$, with $H_0 = \operatorname{Id}$. Let $\tilde H_t: I\times \R \rightarrow I \times \R$ denote its lift to the universal covering, chosen continuously with $\tilde H_0 = \operatorname{Id}$, given by the map $\mathsf p:I\times \R\rightarrow I\times\Sone$ defined by $\mathsf{p}(s,\theta) = (s, \theta \mod 2\pi)$. Then, $\mathsf p \circ \tilde H_t = H_t \circ \mathsf p$ and we observe that
\begin{align*}
    \mathsf p^{-1} (\mathcal R) = \bigsqcup_{k\in \Z} \mathcal R_k, \quad \mathcal R_k := (a, b ) \times (\alpha + 2\pi k, \beta + 2\pi k).
\end{align*}
Since $\mathcal N = \Xi_t^{v_0}(\mathcal N)$ we have
\begin{align*}
    \mathsf p \circ \tilde H_t (\mathcal R_k) = H_t\circ \mathsf{p}(\mathcal R_k) = \Phi^{-1} \circ \Xi_t^{v_0} (\mathcal N) = \Phi^{-1}(\mathcal N) = \mathcal R.
\end{align*}
and thus $\tilde H_t (\mathsf p^{-1}(\mathcal R)) = \mathsf p^{-1}(\mathcal R)= \bigsqcup_{j\in \Z} \mathcal R_j$ for all $k\in \Z$, namely $\tilde H_t$ maps the union of the $\mathcal R_k$ onto itself. Since $\tilde H_t$ is a homeomorphism and $\mathcal R_0$ is connected, $\tilde H_t(\mathcal R_0)$ is connected as well. As the $\lbrace \mathcal R_k \rbrace_{k\in \Z}$ are the connected components of $\mathsf{p}^{-1}(\mathcal R)$, we must have $\tilde H_t(\mathcal R_0) = \mathcal R_k$, for some $k\in \Z$.

In particular, $\tilde H_t(\overline{\mathcal R}_0) = \overline{\mathcal R}_k$ and $\tilde H_t(\lbrace c \rbrace \times[\alpha,\beta]) = \lbrace c \rbrace \times[\alpha + 2\pi k, \beta + 2\pi k]$ because $\tilde H_t$ preserves each connected component of $\partial\Omega$. Furthermore, $\tilde H_t(a,\alpha) = (a, \alpha + 2\pi k)$ and similarly $\tilde H_t(b,\alpha) = (b, \alpha + 2\pi k)$ since ordered particles cannot exchange order along the lifted line, as they would at some point intersect and $\tilde H_t$ is a homeomorphism. In particular, we must have
\begin{align*}
    \tilde\vartheta_a(t) = \tilde\vartheta_b(t) = 2\pi k,
\end{align*}
thus reaching a contradiction with $|\tilde\vartheta_b(t) - \tilde\vartheta_a(t)| \geq \frac{|\nu_b - \nu_a|}{2} \,t >0$ for all $t>0$. To obtain wandering for all non-zero times, if $v_0\in \mathcal U_{\varepsilon_1}$ and $v(-T)\in \mathcal{U}_{\varepsilon_1}$ for some $T>0$, then the Euler solution starting at $v(-T)\in \mathcal{U}_{\varepsilon_1}$ returns at time $T>0$ to $\mathcal{U}_{\varepsilon_1}$, a contradiction with the positive-time conclusion.
\end{proof}

\begin{remark}\label{rem:wandering-stability}
Several observations are in order:
\begin{enumerate}
    \item The steady state $\omega_F=\curl u_F$ need not be isolated. Its only fundamental property is that it has finitely many vorticity values and its velocity field is close to an $L^2$-Lyapunov stable velocity field (it satisfies \eqref{eq:u0-L2-stability}) with different boundary travel times. A close inspection of the proof shows that $\mathcal U_{\varepsilon_1}$ is also a wandering neighborhood close to $u_0$, with arbitrarily small distance obtained by letting $F\rightarrow F_0$. \smallbreak
    \item The background steady field $u_0$ is $L^2$-Lyapunov stable if it is Arnold stable, namely $F_0'(\psi_0)>0$ in $\overline{\Omega}$ or $F_0'(\psi_0) < 0$ in $\overline{\Omega}$ and, for some $\kappa>0$, 
    \begin{equation}\label{eq:arnold-coercivity}
    \int_\Omega\bigl(|\nabla z|^2+F_0'(\psi_0)z^2\bigr)\dd x
    \geq\kappa\int_\Omega|\nabla z|^2\dd x,
    \qquad z\in V,
    \end{equation}
     where $V:=\{z\in H^1(\Omega):z|_{\Gamma_{\rm out}}=0,\ z|_{\Gamma_{\rm in}}\text{ is constant}\}$, see \cite{Arnold1965} and \cite[Theorem~1.5 and Remarks~1.6--1.7]{WangZuo2022}. 
     Satisfying either strict criterion is an open condition in $C^1$, so it is preserved among $C^1$ perturbations of $F_0$ and corresponding $C^1$ perturbations of $\psi_0$.
    \end{enumerate}
\end{remark}

\begin{corollary}\label{cor:wandering-arnold-stable}
    Non-stagnant Arnold stable steady states in $\Omega$ admit non-empty $\mathcal{Y}$-open sets of initial data arbitrarily close to them that wander for all non-zero time. Here, the regularity and parametrization assumptions of Proposition~\ref{prop:wandering} and one of the strict stability criteria in Remark~\ref{rem:wandering-stability} are assumed.
\end{corollary}

\begin{proof}
If $u_0$ is such a state with distinct travel times on $\partial\Omega$, then Proposition~\ref{prop:wandering} applies immediately. Otherwise, if $\mathsf T_a(u_0) = \mathsf T_b(u_0)$ we find $F_\lambda = F_0 + \lambda G$ for some smooth $G:I\rightarrow \R$ such that $\widehat\psi_\lambda$ is obtained by the implicit function theorem (confer Proposition~\ref{con:rough-approximation}\ref{con:rough-steady}). Indeed, the strict criterion in Remark~\ref{rem:wandering-stability} yields $\Delta - F_0'(\psi_0)$ Dirichlet invertible and the implicit function theorem gives a $C^1$ branch in $W^{2,p}(\Omega)$ for $p>2$ with $\widehat\psi_0 = \psi_0$ and
\begin{align*}
    \begin{cases}
        \Delta \widehat\psi_\lambda = F_\lambda(\widehat \psi_\lambda), & \text{ in } \Omega, \\
        \widehat \psi_\lambda = \psi_0, & \text{ on } \partial\Omega.
    \end{cases}
\end{align*}
Elliptic Schauder regularity then yields $\widehat \psi_\lambda\in C^2(\overline{\Omega})$.
Moreover $\widehat u_\lambda = \nabla^\perp\widehat\psi_\lambda$ satisfies the conditions of Proposition~\ref{prop:wandering}: note that
$\widehat u_\lambda$ is non-stagnant Arnold stable for $\lambda>0$ small enough and we shall just see that $\mathsf{T}_a(\widehat u_\lambda) \neq \mathsf{T}_b(\widehat u_\lambda)$ for some small non-zero $\lambda$. The derivative $\widehat v := (\partial_\lambda \widehat\psi_\lambda)|_{\lambda=0}$ satisfies
    \begin{align*}
        (\Delta - F'_0(\psi_0))\widehat v = G(\psi_0)
    \end{align*}
    in $\Omega$, with $\widehat v = 0$ on $\partial \Omega$. Let $\varrho(\lambda) := \mathsf{T}_b(\widehat u_\lambda) - \mathsf{T}_a(\widehat u_\lambda)$ and note that the outward normal derivative of $\psi_0$ is $|\nabla\psi_0|$ on $\Gamma_b$ and $-|\nabla\psi_0|$ on $\Gamma_a$. Then,
    \begin{align*}
        \varrho'(0) = - \int_{\partial\Omega} \frac{\partial_n \widehat v}{|\nabla\psi_0|^2}\dd \ell.
    \end{align*}
    As $|\nabla\psi_0|>0$ in $\overline{\Omega}$ and $\Delta - F'_0(\psi_0)$ has trivial Dirichlet kernel, let $w$ solve the adjoint Dirichlet problem
    \begin{align*}
        \begin{cases}
            (\Delta - F'_0(\psi_0))w = 0, & \text{ in } \Omega,\\
            w = -|\nabla\psi_0|^{-2}, & \text{ on } \partial\Omega.
        \end{cases}
    \end{align*}
    Then, the divergence theorem and co-area formula show
    \begin{align*}
        \varrho'(0) = \int_{\partial\Omega}w \partial_n \widehat v \dd \ell = \int_\Omega w G(\psi_0) \dd x = \int_a^b G(s) \left(\int_{\lbrace \psi_0 = s\rbrace} \frac{w}{|\nabla\psi_0|} \dd \ell\right) \dd s
    \end{align*}
    Since $w < 0$ on $\Gamma_a$, it stays negative in a small tubular neighborhood of $\Gamma_a$ in $\Omega$, so that $\left(\int_{\lbrace \psi_0 = s\rbrace} \frac{w}{|\nabla\psi_0|} \dd \ell\right)$ is negative for $s\in [a,a+\delta_a)$ for some $\delta_a>0$. Hence, for any non-zero non-negative smooth bump function $G$ supported in $(a,a+\delta_a)$ we observe that $\varrho'(0)<0$ and thus $\mathsf{T}_b(\widehat u_\lambda) \neq \mathsf{T}_a(\widehat u_\lambda)$ for $\lambda>0$ small enough. 
\end{proof}

We finish the section by constructing a smooth reference Euler flow satisfying the hypotheses of Proposition~\ref{con:rough-approximation}.

\begin{proposition}\label{con:seed}
There are a smooth bounded annular domain $\Omega$ and a function
$\bar\psi\in C^\infty(\overline\Omega)$ satisfying
\[
 \Delta\bar\psi=\bar\psi^{-1},\qquad
 \bar\psi|_{\Gamma_{\rm in}}=1,\qquad
 \bar\psi|_{\Gamma_{\rm out}}=R,
 \qquad R=e^{\pi/4},
\]
such that $\partial_r\bar\psi\geq1/2$, where $r=|x|$, and every
level curve is connected and closed. Moreover,
\begin{equation}\label{con:seed-oscillation}
 \mathbf{o}:=\inf_{s\in[1,R]}
 \osc_{\{\bar\psi=s\}}|\nabla\bar\psi|^2>0,
\end{equation}
$\Delta+\bar\psi^{-2}\operatorname{Id}$ has trivial Dirichlet kernel,
and there is a smooth diffeomorphism
\begin{equation}\label{con:seed-parametrization}
 \bar\Phi:[1,R]\times\Sone\longrightarrow\overline\Omega,
 \qquad \bar\psi(\bar\Phi(s,\theta))=s.
\end{equation}
{Moreover, the reference velocity
$\bar u:=\nabla^\perp\bar\psi$ satisfies the stability condition
\eqref{eq:u0-L2-stability} and has distinct travel times on the two
boundary components.}
\end{proposition}

\begin{proof}
Set $\sigma=\pi/8$, $r_-=e^{-\sigma}$ and $r_+=Re^\sigma$.
We work first on the fixed annular domain
$\widehat\Omega=\{r_-<|x|<r_+\}$. The radial function $r=|x|$
satisfies $\Delta r=1/r$, and the linearized operator is
$L_0=\Delta+r^{-2}\operatorname{Id}$. To check its invertibility, let
$\mathcal A(c,d)=\{c<|x|<d\}$ and $\ell=\log(d/c)$. The change of
variables $r=ce^t$ and the one-dimensional Dirichlet Poincar\'e
inequality give
\begin{equation}\label{con:coercivity}
 \int_{\mathcal A(c,d)}\bigl(|\nabla z|^2-r^{-2}z^2\bigr)\dd x
 \geq \left(1-\frac{\ell^2}{\pi^2}\right)
       \int_{\mathcal A(c,d)}|\nabla z|^2\dd x,
 \qquad z\in H^1_0(\mathcal A(c,d)).
\end{equation}
Indeed, the two integrals before subtraction become, respectively,
$\int(|z_t|^2+|z_\theta|^2)\dd t\dd\theta$ and
$\int z^2\dd t\dd\theta$ on $(0,\ell)\times\Sone$.
For $\widehat\Omega$ one has $\ell=\pi/2$, so the constant in
\eqref{con:coercivity} is $3/4$. Lax--Milgram and elliptic regularity
therefore make $L_0:C^{2,\alpha}_0(\overline{\widehat\Omega})\to C^{0,\alpha}(\overline{\widehat\Omega})$ an isomorphism, for $0<\alpha<1$.

Now, let $\beta=\sqrt{3}$, set $w(r,\theta)=r^\beta\cos(2\theta)$ and
consider the smooth map
\[
 \begin{aligned}
 \mathcal G:\mathcal U\subset\mathbb R\times
 C^{2,\alpha}_0(\overline{\widehat\Omega})
 &\longrightarrow C^{0,\alpha}(\overline{\widehat\Omega}),\\
 (\varepsilon,v)&\longmapsto
 \Delta(r+\varepsilon w+v)-(r+\varepsilon w+v)^{-1},
 \end{aligned}
\]
where $\mathcal U$ is an open neighborhood of $(0,0)$ chosen so that
$r+\varepsilon w+v>0$ on $\overline{\widehat\Omega}$.
Since $\mathcal G(0,0)=0$ and $D_v\mathcal G(0,0)=L_0$ is an
isomorphism, the implicit function theorem gives a smooth branch
$\varepsilon\mapsto v_\varepsilon$ with $v_0=0$ and
$\mathcal G(\varepsilon,v_\varepsilon)=0$. Moreover,
\[
 D_\varepsilon\mathcal G(0,0)=L_0w
 = (\beta^2-3)r^{\beta-2}\cos(2\theta)=0.
\]
Differentiating $\mathcal G(\varepsilon,v_\varepsilon)=0$ at
$\varepsilon=0$ and using the invertibility of $L_0$ gives
$\left.\partial_\varepsilon v_\varepsilon\right|_{\varepsilon=0}=0$.
By Taylor's theorem, we therefore have
\begin{equation}\label{con:seed-expansion}
 \psi_\varepsilon=r+\varepsilon w+v_\varepsilon,
 \qquad \Delta\psi_\varepsilon=\psi_\varepsilon^{-1},
 \qquad \|v_\varepsilon\|_{C^{2,\alpha}(\overline{\widehat\Omega})}\leq C\varepsilon^2.
\end{equation}
Smooth boundary values and elliptic bootstrapping give
$\psi_\varepsilon\in C^\infty(\overline{\widehat\Omega})$.

For small $\varepsilon$, we have $\partial_r\psi_\varepsilon\geq1/2$,
$\psi_\varepsilon(r_-,\theta)<1$ and
$\psi_\varepsilon(r_+,\theta)>R$. Consequently, for each
$(s,\theta)\in[1,R]\times\Sone$ there is a unique
$\rho_\varepsilon(s,\theta)$ with
\[
 \psi_\varepsilon(\rho_\varepsilon(s,\theta),\theta)=s.
\]
Define
\[
 \Omega_\varepsilon=\{1<\psi_\varepsilon<R\},\qquad
 \Phi_\varepsilon(s,\theta)
 =\rho_\varepsilon(s,\theta)(\cos\theta,\sin\theta).
\]
The map $\Phi_\varepsilon:[1,R]\times\Sone\to\overline{\Omega_\varepsilon}$
is a smooth diffeomorphism, since
$\partial_s\rho_\varepsilon=(\partial_r\psi_\varepsilon)^{-1}>0$.
Moreover, implicit differentiation in \eqref{con:seed-expansion} gives
\begin{equation}\label{con:level-expansion}
 \begin{aligned}
 \rho_\varepsilon(s,\theta)
 &=s-\varepsilon w(s,\theta)+O(\varepsilon^2)
 &&\hbox{in }C^1([1,R]\times\Sone),\\
 \rho_\varepsilon-s&=O(\varepsilon)
 &&\hbox{in }C^2([1,R]\times\Sone).
 \end{aligned}
\end{equation}
All remainder estimates in this proof are uniform on the indicated
compact sets.

The squared speed on these level curves satisfies
\[
 \begin{split}
 S_\varepsilon(s,\theta)
 &:=|\nabla\psi_\varepsilon|^2(\Phi_\varepsilon(s,\theta))\\
 &=1+2\varepsilon\beta s^{\beta-1}\cos(2\theta)
     +O(\varepsilon^2).
 \end{split}
\]
Thus its second Fourier coefficient is
\[
 \frac1\pi\int_0^{2\pi}S_\varepsilon(s,\theta)
                       \cos(2\theta)\dd\theta
 =2\varepsilon\beta s^{\beta-1}+O(\varepsilon^2)
 \geq\beta\varepsilon
\]
for all $s\in[1,R]$ and all sufficiently small $\varepsilon>0$.
On the other hand, for any continuous periodic function $S$, subtracting the midpoint of
its range in this integral shows that
\[
 \left|\frac1\pi\int_0^{2\pi}S(\theta)\cos(2\theta)\dd\theta\right|
 \leq\frac2\pi\osc S.
\]
It then follows that the oscillation on every level curve, including the two
boundary level curves, is at least $\pi\beta\varepsilon/2$.

Finally, the Dirichlet inverse must persist on the new domain. Let
$\Omega_0=\{1<r<R\}$ and define
$T_\varepsilon(r,\theta)=\Phi_\varepsilon(r,\theta)$. This map is
$O(\varepsilon)$ close to the identity in $C^1(\overline{\Omega}_0)$, and
$\psi_\varepsilon\circ T_\varepsilon=r$. For
$z\in H^1_0(\Omega_\varepsilon)$ and $v=z\circ T_\varepsilon$,
change of variables gives
\[
 \int_{\Omega_\varepsilon}(|\nabla z|^2-\psi_\varepsilon^{-2}z^2)\dd x
 =\int_{\Omega_0}(|\nabla v|^2-r^{-2}v^2)\dd x
       +O(\varepsilon)\|\nabla v\|_{L^2(\Omega_0)}^2.
\]
Here the potential error is controlled by the same Poincar\'e
inequality used in \eqref{con:coercivity}. On $\Omega_0$ the
coercivity constant is $15/16$, so, after decreasing $\varepsilon$,
the left-hand side is bounded below by
$\frac3{16}\|\nabla z\|_{L^2(\Omega_\varepsilon)}^2$. In particular, the Dirichlet kernel
is trivial. 

To verify the stability condition, let $\mathcal V_\varepsilon$ be the space
$V$ of Remark~\ref{rem:wandering-stability} on $\Omega_\varepsilon$.
For $v\in\mathcal V_0$, the one-dimensional Poincar\'e inequality with
one zero endpoint, applied in $t=\log r$, gives
\[
 \int_{\Omega_0}r^{-2}v^2\dd x
 \leq\frac{4(\log R)^2}{\pi^2}
       \int_{\Omega_0}|\nabla v|^2\dd x
 =\frac14\int_{\Omega_0}|\nabla v|^2\dd x.
\]
Thus the quadratic form
$\int_{\Omega_0}(|\nabla v|^2-r^{-2}v^2)\dd x$ has coercivity
margin $3/4$ on $\mathcal V_0$. The pullback $v=z\circ T_\varepsilon$
maps $\mathcal V_\varepsilon$ onto $\mathcal V_0$, and the same
change-of-variables estimate as above, now using this one-endpoint
Poincar\'e inequality, yields, for all sufficiently small $\varepsilon$,
\[
 \int_{\Omega_\varepsilon}
       (|\nabla z|^2-\psi_\varepsilon^{-2}z^2)\dd x
 \geq c\|\nabla z\|_{L^2(\Omega_\varepsilon)}^2,
 \qquad z\in\mathcal V_\varepsilon,
\]
with $c>0$ independent of $\varepsilon$. Since
$F_0'(\psi_\varepsilon)=-\psi_\varepsilon^{-2}<0$,
Remark~\ref{rem:wandering-stability} gives
\eqref{eq:u0-L2-stability} for
$u_\varepsilon:=\nabla^\perp\psi_\varepsilon$.
Finally, the $C^1$ expansions of the level parametrizations and velocities give
\[
 \mathsf T_1(u_\varepsilon)=2\pi+O(\varepsilon),\qquad
 \mathsf T_R(u_\varepsilon)=2\pi R+O(\varepsilon).
\]
As $R>1$, these travel times are distinct after decreasing
$\varepsilon$ once more.

Fix one such positive $\varepsilon$ and set
$\Omega=\Omega_\varepsilon$, $\bar\psi=\psi_\varepsilon|_\Omega$
and $\bar\Phi=\Phi_\varepsilon$.
\end{proof}

\begin{remark}\label{con:reference-profile}
For the reference Euler flow in Proposition~\ref{con:seed}, the formula
\eqref{con:profile} for $F_\eta$ is particularly simple:
\begin{equation}\label{con:specific-profile}
 F_\eta(s)=\eta\left\lfloor\frac1{\eta s}\right\rfloor
                 +\eta\mathbf1_B(s),\qquad 1\leq s\leq R.
\end{equation}
If $\eta<1/(2R)$, then $F_\eta(s)>1/(2R)$ on $[1,R]$.
Although this law is uniformly close to $1/s$, its essential
variation on every interval prevents the usual smooth perturbation
argument from producing nearby independent equilibria.
\end{remark}

Applying Proposition~\ref{con:rough-approximation} with this reference
Euler flow and the Borel law in Remark~\ref{con:reference-profile} gives the
isolated equilibrium required in Theorem~\ref{intro:theorem-a}.

\begin{corollary}\label{con:constructed}
There exist a smooth bounded annular domain $\Omega$ and a nonzero isolated steady
Euler flow $u_*=\nabla^\perp\psi_*\in\mathcal Y(\Omega)$ with positive,
finitely valued vorticity $\omega_*=F_*(\psi_*)$. For some $R>1$, the
stream function $\psi_*$ admits a $C^1$ parametrization
\begin{equation}\label{con:constructed-parametrization}
 X_*:[1,R]\times\Sone\longrightarrow\overline\Omega,
 \qquad \psi_*(X_*(s,\theta))=s,
\end{equation}
where $X_*$ is a diffeomorphism. The function $\psi_*$ has no critical
points, and its level curves are connected, with a positive uniform
lower bound for the oscillation of $|\nabla\psi_*|^2$.

The Borel law $F_*:[1,R]\to\R$ is not essentially constant on any
nonempty open subinterval of $(1,R)$. If $Q_*$ denotes its set of
essential values, then
\begin{equation}\label{con:gap}
 2\delta_*:=\min\{|q-q'|:q,q'\in Q_*,\ q\ne q'\}>0.
\end{equation}
\end{corollary}
\begin{proof}
Apply Proposition~\ref{con:rough-approximation} to the reference Euler flow of
Proposition~\ref{con:seed}, with $F_0(s)=1/s$ and $F=F_\eta$
as in Remark~\ref{con:reference-profile}, choosing one sufficiently
small $\eta<1/(2R)$, also small enough for Proposition \ref{prop:wandering}. Set $\psi_*=\psi_{F_\eta}$ and
$F_*=F_\eta$. Positivity follows from
\eqref{con:specific-profile}. The profile has at least two essential
values, and its finite range gives \eqref{con:gap}.
Lemma~\ref{rig:rigidity} applies with this actual separation.
\end{proof}

\section{Dynamical consequences of the construction. Proof of Theorem~\ref{intro:theorem-a}.}
\label{sec:dynamics}

Let $u_*=\nabla^\perp\psi_*$ be the steady Euler flow constructed in
Corollary~\ref{con:constructed}, and write
\[
 \omega_*=F_*(\psi_*),\qquad
 Q_*:=\operatorname{ess\,ran}F_*=\{q_1,\ldots,q_N\},\qquad
 2\delta_*:=\min_{i\ne j}|q_i-q_j|.
\]
Lemma~\ref{rig:rigidity} already supplies the dynamical rigidity
in Theorem~\ref{intro:theorem-a}\textup{(i)}, { and Proposition \ref{prop:wandering} the wandering in Theorem~\ref{intro:theorem-a}(ii).} To prove the third
assertion, we use a conserved quantity associated with the finite
set $Q_*$.

\begin{proposition}\label{dyn:nonreturn}
Let $v(t)$ be a Yudovich solution, with vorticity
$\omega(t)$ and initial vorticity $\omega_0$, such that
\[
 \|\omega_0-\omega_*\|_{L^\infty}<\delta_*,\qquad
 Z(\omega_0)=Z(\omega_*),\qquad \omega_0\ne\omega_*.
\]
Set $r:=\|\omega_0-\omega_*\|_{L^2}>0$ and $Z_*:=Z(\omega_*)$.
Then
\begin{equation}\label{dyn:correlation-gap}
 \int_\Omega\omega(t)\omega_*\,dx\leq Z_*-\frac{r^2}{2}
 \qquad\text{for all }t\in\mathbb R.
\end{equation}
Consequently,
\[
 \omega_*\notin
 \overline{\{\omega(t):t\in\mathbb R\}}^{\,\mathrm{w}^*},
 \qquad
 \inf_{t\in\mathbb R}\|v(t)-u_*\|_{L^2}>0.
\]
\end{proposition}

\begin{proof}
The initial closeness implies that $\omega_*(x)$ is the unique
nearest element of $Q_*$ to $\omega_0(x)$ for a.e. $x\in\Omega$.
Thus, it follows that
\[
 \operatorname{dist}(\omega_0(x),Q_*)
 =|\omega_0(x)-\omega_*(x)|.
\]
Transport and preservation of volume give the conserved quantity
\begin{equation}\label{dyn:distance-invariant}
 \int_\Omega\operatorname{dist}(\omega(t),Q_*)^2\,dx
 =\int_\Omega\operatorname{dist}(\omega_0,Q_*)^2\,dx=r^2.
\end{equation}
Since $\omega_*(x)\in Q_*$ a.e. in $\Omega$,
\begin{equation}\label{dyn:distance-lower-bound}
 \|\omega(t)-\omega_*\|_{L^2}^2\geq \int_\Omega\operatorname{dist}(\omega(t),Q_*)^2\,dx\geq  r^2
 \quad\text{for every }t\in\mathbb R.
\end{equation}
On the other hand, enstrophy conservation and its prescribed initial value yield
\[
 \|\omega(t)-\omega_*\|_{L^2}^2 \leavevmode{\ignorespaces = \int_\Omega \left(\omega^2(t) - 2 \omega(t)\omega_* + \omega_*^2 \right) \dd x} 
 =2Z_*-2\int_\Omega\omega(t)\omega_*\,dx,
\]
which proves~\eqref{dyn:correlation-gap}.

The linear functional $f\mapsto\int_\Omega f\omega_*\,dx$ is
weak-$*$ continuous on $L^\infty$, because $\omega_*\in L^1$.
The entire orbit therefore lies in the weak-$*$ closed half-space
where this functional is at most $Z_*-r^2/2< Z_*$, whereas its value
at $\omega_*$ is $Z_*$. This excludes weak return, i.e.
\[
 \omega_*\notin
 \overline{\{\omega(t):t\in\mathbb R\}}^{\,\mathrm{w}^*}.
\]

For the velocity conclusion, put
$L=\|\omega_0\|_{L^\infty}+\|\omega_*\|_{L^\infty}$ and choose
$\varphi\in C_c^\infty(\Omega)$ such that
\[
 \|\varphi-\omega_*\|_{L^1}<\frac{r^2}{4L}.
\]
The conserved $L^\infty$ bound and~\eqref{dyn:correlation-gap}
imply, for every real $t$, 
\[
 \int_\Omega(\omega_*-\omega(t))\varphi\,dx = \int_\Omega(\omega_*-\omega(t))(\varphi - \omega_*) \dd x - \int_\Omega \omega(t)\omega_* \dd x + Z_*
 \geq\frac{r^2}{2}-L\|\varphi-\omega_*\|_{L^1}
 >\frac{r^2}{4}.
\]
By distributional integration by parts,
\[
 \left|\int_\Omega(\omega_*-\omega(t))\varphi\,dx\right|
 =\left|\int_\Omega(v(t)-u_*)\cdot\nabla^\perp\varphi\,dx\right|
 \leq\|v(t)-u_*\|_{L^2}\|\nabla\varphi\|_{L^2}.
\]
In particular $\|\nabla\varphi\|_{L^2}>0$, and
\begin{equation}\label{dyn:velocity-gap}
 \inf_{t\in\mathbb R}\|v(t)-u_*\|_{L^2}
 \geq\frac{r^2}{4\|\nabla\varphi\|_{L^2}}>0.
\end{equation}
\end{proof}

\begin{remark}
Notice that the proposition does not require equality of energy
or circulations. Those constraints enter Theorem~\ref{intro:theorem-a}
only to ensure that $u_0\ne u_*$ implies $\omega_0\ne\omega_*$,
and to normalize the stationary multiple when its first assertion is applied at fixed energy.
\end{remark}

\begin{proof}[Proof of Theorem~\ref{intro:theorem-a}]
Corollary~\ref{con:constructed} provides the smooth annular domain and
the nonzero reference Euler flow satisfying
all the hypotheses of Lemma~\ref{rig:rigidity}. That lemma proves
both its stationary isolation and the trajectory statement in
Theorem~\ref{intro:theorem-a}\textup{(i)}, with the radius $\delta_*$ determined
above.

{Proposition~\ref{con:seed} gives the stability condition
\eqref{eq:u0-L2-stability} and distinct boundary travel times. The choice of
$\eta$ in Corollary~\ref{con:constructed} also meets the threshold in
Proposition~\ref{prop:wandering}, which applies to the same
$u_*=\nabla^\perp\psi_{F_\eta}$ and proves~\textup{(ii)}.}

Finally, suppose that
$u_0\in\mathcal M_Z(u_*)\setminus\{u_*\}$ and
$\|\omega_0-\omega_*\|_{L^\infty}<\delta_*$. If
$\omega_0=\omega_*$, the equality of boundary circulations and
Lemma~\ref{rig:streamfunction} would imply $u_0=u_*$, a
contradiction. Therefore $\omega_0\ne\omega_*$, and
Proposition~\ref{dyn:nonreturn} proves both conclusions in
\eqref{intro:nonreturn}.

\end{proof}

\section{Approximation by isolated steady Euler flows. Proof of Theorem~\ref{intro:approximation}.}
\label{sec:approximation}

We prove Theorem~\ref{intro:approximation} and
Corollary~\ref{intro:harmonic-reference}. The construction of
Proposition~\ref{con:rough-approximation} applies to a smooth reference
Euler flow whose speed varies along every level curve. We first show
that this property can be obtained by arbitrarily small perturbations
of the vorticity law on a fixed noncircular annular domain. We then apply the
construction of Section~\ref{sec:construction} and prove that, on a
concentric circular annular domain, a neighborhood of a non-stagnant circular
Euler flow contains no isolated Yudovich equilibria.

We write $u=\nabla^\perp\psi$ for the Euler flow in the theorem
and denote the boundary values of $\psi$ by $a<b$. These values are distinct,
since otherwise $\psi$ would have an interior extremum. The same
argument shows that its range is $I=[a,b]$. We also set
$
L:=\Delta-F'(\psi)\operatorname{Id}.
$
The triviality of its Dirichlet kernel, elliptic regularity and the
Fredholm alternative imply that
$L:W^{2,p}(\Omega)\cap W^{1,p}_0(\Omega)\to L^p(\Omega)$ is an
isomorphism for every $1<p<\infty$. We fix one exponent $p>2$.

\begin{proposition}\label{app:den-density}
{Let $\Omega$ be a bounded smooth annular domain which is not the region
between two concentric circles. Suppose that
$\psi\in C^\infty(\overline\Omega)$, $F\in C^\infty(I)$, and}
\[
 \Delta\psi=F(\psi),\qquad
 \psi|_{\Gamma_a}=a<b=\psi|_{\Gamma_b},\qquad
 \min_{\overline\Omega}|\nabla\psi|>0.
\]
Assume also that $L=\Delta-F'(\psi)\operatorname{Id}$ is an
isomorphism from $W^{2,p}(\Omega)\cap W^{1,p}_0(\Omega)$ onto
$L^p(\Omega)$. Then, there are
$F_j\in C^\infty(I)$ and $\psi_j\in C^\infty(\overline\Omega)$
satisfying
\begin{equation}\label{app:den-approximation-equations}
 \Delta\psi_j=F_j(\psi_j),\qquad
 \psi_j|_{\Gamma_a}=a,\qquad \psi_j|_{\Gamma_b}=b,
\end{equation}
such that, for every integer $r\geq0$,
\begin{equation}\label{app:den-smooth-convergence}
 \|F_j-F\|_{C^r(I)}
 +\|\psi_j-\psi\|_{C^r}\longrightarrow0,
\end{equation}
and
\begin{equation}\label{app:den-positive-oscillation}
 \min_{s\in I}\osc_{\{\psi_j=s\}}|\nabla\psi_j|^2>0.
\end{equation}
The functions $\psi_j$ have no critical points on $\overline\Omega$,
and their Dirichlet linearizations are invertible.
\end{proposition}

\subsection{Geometry of the level curves}

We first describe the level curves of a pair $(F,\psi)$ satisfying
the hypotheses of Proposition~\ref{app:den-density}. We write
$m=\min_{\overline\Omega}|\nabla\psi|>0$. Since an interior
extremum would be a critical point, we have $a<\psi<b$ in $\Omega$.
Moreover, $\nabla\psi$ is normal to the boundary, points into
$\Omega$ on $\Gamma_a$, and points out of $\Omega$ on $\Gamma_b$.

We consider the smooth vector field
\[
 V=\frac{\nabla\psi}{|\nabla\psi|^2},
\]
which we extend smoothly to a neighborhood of $\overline\Omega$.
Let $\mathcal F_t$ denote its flow, i.e.,
\[
 \frac{d}{dt}\mathcal F_t(x)=V(\mathcal F_t(x)),
 \qquad \mathcal F_0(x)=x,
\]
and let $\gamma:\Sone\to\Gamma_a$ be a smooth parametrization of $\Gamma_a$ {with $\gamma'\neq 0$}. Along each
trajectory, $\frac{d}{dt}\psi(\mathcal F_t(x))=1$ while it remains
in $\overline\Omega$. A trajectory starting on $\Gamma_a$ therefore
enters $\Omega$ and cannot meet its boundary again before time $b-a$,
when it reaches $\Gamma_b$. Conversely, a trajectory starting at
$x\in\overline\Omega$ reaches $\Gamma_a$ after flowing backwards
for time $\psi(x)-a$. Existence on these time intervals follows
from compactness and the boundary signs above.

Thus, the map $\Phi:I\times\Sone\longrightarrow\overline\Omega$ given by
\begin{equation}\label{app:parametrization}
 \Phi(s,\theta):=\mathcal F_{s-a}(\gamma(\theta)),
\end{equation}
satisfies $\psi(\Phi(s,\theta))=s$, and it is bijective by uniqueness of trajectories. Note that for $s\in [a,b]$ the closed streamline $\lbrace \psi = s \rbrace$ is given by $\lbrace \Phi(s,\theta) : \theta\in \Sone \rbrace$.
Since {$\partial_\theta\Phi=D\mathcal F_{s-a}(\gamma(\theta))\gamma'(\theta)$} is nonzero and tangent to the level curve,
whereas $\nabla\psi\cdot\partial_s\Phi=1$, its differential is
invertible. Consequently, $\Phi$ is a smooth diffeomorphism, also
up to the boundary, and each level curve is connected. We orient
$\theta$ so that
\[
 D_\Phi(s,\theta):=\det D\Phi(s,\theta)>0.
\]
By continuity and compactness, $D_\Phi$ has a positive minimum,
and the first derivatives of $\Phi$ and $\Phi^{-1}$ are bounded.
We also have 
\begin{equation}\label{app:gen-chart-derivative}
 \partial_s\Phi= \frac{\dd}{\dd s} \mathcal F_{s-a}(\gamma(\theta)) = V(\mathcal F_{s-a}(\gamma(\theta))) =  V(\Phi)
\end{equation}
and thus $\|\partial_s\Phi\|_{L^\infty(I\times\Sone)}=1/m$.

We next describe the smooth perturbations of the vorticity law.
Given $G\in C^\infty(I)$, we set
\[
 F_\eps=F+\eps G,
 \qquad \eps\in\R,
\]
and extend $F$ and $G$ smoothly beyond $I$. For a fixed
$0<\alpha<1$, we consider the smooth map
\[
 \mathcal G:C_0^{2,\alpha}(\overline\Omega)\times\R
 \longrightarrow C^{0,\alpha}(\overline\Omega),\qquad
 \mathcal G(h,\eps)=\Delta(\psi+h)-F_\eps(\psi+h).
\]
Its derivative in $h$ at $(0,0)$ is $L$. The assumed Sobolev
inverse and the global Dirichlet Schauder estimates imply that
$L:C_0^{2,\alpha}(\overline\Omega)\to C^{0,\alpha}(\overline\Omega)$ is an
isomorphism. The implicit function theorem therefore gives a locally
unique smooth family $\psi_\eps$ with the prescribed boundary
values. Applying the same argument in higher H\"older spaces and
using local uniqueness, we obtain smooth dependence in every fixed
$C^r(\overline\Omega)$ space, see
\cite[Chapters~6, 9 and Section~17.2]{GT}.

For sufficiently small $|\eps|$, the functions $\psi_\eps$ have no critical
points, so $a<\psi_\eps<b$ in $\Omega$. Thus, the extensions of
the profiles beyond $I$ do not affect the solution. We denote the
first variation by $v=\partial_\eps\psi_\eps|_{\eps=0}$. Differentiating
the equation and the boundary values, we obtain
\begin{equation}\label{app:den-linear-response}
 Lv=G(\psi),\qquad v|_{\partial\Omega}=0.
\end{equation}

We now compute the first variation of the speed on the perturbed
level curves. We define, on $\overline\Omega$,
\begin{equation}\label{app:den-operators}
 S=|\nabla\psi|^2,\qquad
 \leavevmode{\ignorespaces D_\psi}=\nabla\psi\cdot\nabla,\qquad
 H=\frac{D_\psi S}{S},\qquad
 \mathcal Rv=2D_\psi v-Hv.
\end{equation}
We also abuse notation and write $S(s,\theta):=S(\Phi(s,\theta))$, and similarly
$H(s,\theta) := H(\Phi(s,\theta))$, so that
\begin{align*}
    \partial_s S(s,\theta) & = \partial_s \left( |\nabla\psi|^2\circ \Phi(s,\theta) \right) \\
    &= (\partial_s \Phi)(s,\theta) \cdot [\nabla ( |\nabla\psi|^2)](\Phi(s,\theta)) = V(\Phi(s,\theta)) \cdot   [\nabla ( |\nabla\psi|^2)](\Phi(s,\theta)) = H(s,\theta).
\end{align*}
As in
\eqref{con:perturbed-chart}, we parametrize the perturbed level
curves using the identity
\begin{equation*}\label{eq:level_sets_psi_eps}
 \psi_\eps(\Phi(\sigma_\eps(s,\theta),\theta))=s.
\end{equation*}
Now, let $X_\eps(s,\theta) = \Phi(\sigma_\eps(s,\theta),\theta)$. For $\eps=0$, we recover $\sigma_0(s,\theta)=s$ and $X_0(s,\theta) = \Phi(s,\theta)$. Moreover, differentiating the identity above in $\eps$ we reach
\begin{align*}
    0 = (\partial_\eps\psi_\eps)(X_\eps) + (\partial_\eps \sigma_\eps) (\partial_s \Phi)\cdot \nabla \psi_\eps(X_\eps).
\end{align*}
At $\eps=0$ we note that $(\partial_\eps\psi_\eps)(X_\eps)|_{\eps=0}= v(X_0)$ while $(\partial_s \Phi)\cdot \nabla \psi_\eps(X_\eps)|_{\eps=0}=1$. Hence, 
\begin{align*}
    \partial_\eps \sigma_\eps(s,\theta)|_{\eps=0} = - v(\Phi(s,\theta)),
\end{align*}
and for $S_\eps(s,\theta) :=|\nabla\psi_\eps|^2(\Phi(\sigma_\eps(s,\theta),\theta))$, we obtain
\begin{align*}
    \partial_\eps S_\eps = 2 (\nabla \psi_\eps \cdot\nabla (\partial_\eps \psi_\eps))(X_\eps) + (\partial_\eps \sigma_\eps) \partial_s \Phi \cdot |\nabla \psi_\eps|^2(X_\eps),
\end{align*}
which gives
\begin{equation}\label{app:den-moving-speed}
 \left.\partial_\eps S_\eps\right|_{\eps=0} = 2(\nabla \psi \cdot \nabla v)(\Phi) - v(\Phi) H
 =\mathcal Rv\circ\Phi.
\end{equation}
This identity accounts for the displacement of the level curve.
At $s=a,b$, the curve is fixed and $v=0$, so the first variation
reduces to $2D_{\psi}v$.

\begin{lemma}\label{app:den-circle-lemma}
Let $\Omega$ be a bounded smooth annular domain, and let $F,\psi$ satisfy
the equation, boundary conditions, regularity and non-stagnation
assumptions of Proposition~\ref{app:den-density}. If one level curve,
possibly a boundary component, is a circle on which $|\nabla\psi|$
is constant, then $\Omega$ is a concentric circular annulus and
$\psi$ is radial about its center.
\end{lemma}

\begin{proof}
Let $x_*$ be the center of the circle on which $|\nabla\psi|$
is constant. We consider
\[
 w=(x-x_*)^\perp\cdot\nabla\psi.
\]
Since rotations commute with the Laplacian, differentiating
$\Delta\psi=F(\psi)$ along the rotation field gives
$\Delta w=F'(\psi)w$ in $\Omega$. On the circle, $\psi$ is
constant, so $w=0$. In polar coordinates about $x_*$, we have
$w=\partial_\theta\psi$. Moreover, the normal derivative of
$\psi$ is constant, since its magnitude is the nonzero constant
speed and its sign is fixed. Thus,
$\partial_nw=\partial_\theta(\partial_n\psi)=0$ on the circle.

We now explain why these two vanishing traces imply $w=0$, including
when the circle is a boundary component. In a sufficiently small ball $B$
centered at a point of the circle, we retain $w$ on one side contained in
$\Omega$ and extend it by zero to the other side. Integration by
parts, using $w=\partial_nw=0$, shows that the extension satisfies
\[
 \Delta\widetilde w=\widetilde V\widetilde w
 \quad\text{in the sense of distributions in }B,
\]
where $\widetilde V$ is a bounded extension of $F'(\psi)$.
The extension belongs to $W^{2,p}_{\mathrm{loc}}(B)$ for every
$1<p<\infty$ and vanishes on an open subset of $B$. Unique
continuation for this equation gives $\widetilde w=0$ in $B$, see
\cite[pp.~235--236 and Remark~2, pp.~247--248]{Aronszajn}.
A further application of unique continuation in the connected
domain $\Omega$ gives $w=0$ throughout $\Omega$.

For $x\ne x_*$, this identity implies that $\nabla\psi$ is
parallel to $x-x_*$. Since $\nabla\psi$ is nonzero and normal
to each boundary component, we obtain
$\partial_\tau|x-x_*|^2=0$ along the boundary, where $\tau$ denotes
a unit tangent vector. This derivative also
vanishes at $x=x_*$. Hence $|x-x_*|$ is constant on each connected
boundary component, and each component is an entire circle about
$x_*$. The domain is therefore the annulus between two concentric
circles. On this domain, $w=0$ means precisely that $\psi$ is radial.
\end{proof}

We next show that only finitely many level curves can have constant
speed. Let $\mathsf{C}=\{\psi=c\}$ be such a curve, and write
\[
 q=|\nabla\psi||_{\mathsf{C}}>0,\qquad
 n=\frac{\nabla\psi}{|\nabla\psi|},\qquad \kappa=\diver n|_{\mathsf{C}}.
\]
Here, $\kappa$ is the signed curvature determined by $n$. Since
\begin{align*}
    \Delta\psi = \diver (\nabla \psi) = \diver (|\nabla \psi| n) = n\cdot \nabla |\nabla \psi| + |\nabla \psi| \diver n,
\end{align*}
we have $\Delta\psi = \partial_n |\nabla \psi| + q \kappa$ on $\mathsf{C}$. Moreover, $\nabla \psi = |\nabla \psi| n$ yields 
\begin{align*}
    D_\psi S = |\nabla\psi|\partial_n (|\nabla\psi|^2) = 2|\nabla\psi|^2 \partial_n |\nabla \psi|,
\end{align*}
so that $H = \frac{D_\psi S}{S}= 2 \partial_n |\nabla \psi|$. Further using $\Delta\psi = F(c)$ on $\mathsf{C}$, we conclude that
\begin{equation}\label{app:den-H-curvature}
 H|_{\mathsf{C}}=2\bigl(F(c)-q\kappa\bigr).
\end{equation}
A connected regular closed plane curve of constant signed curvature
is a circle. Indeed, the zero-curvature case gives a straight line,
which cannot close, and the nonzero case follows by integrating the
tangent-angle equation. Lemma~\ref{app:den-circle-lemma} therefore
implies that $\kappa$ is nonconstant on every constant-speed level curve
when $\Omega$ is not a concentric circular annulus. The next lemma shows that only finitely many level curves may have constant speed.

\begin{lemma}\label{app:den-finite-levels}
Let $\Omega$ be a bounded smooth annular domain, and let $F,\psi$ satisfy
the equation, boundary conditions, regularity and non-stagnation
assumptions of Proposition~\ref{app:den-density}. If $\Omega$ is not
a concentric circular annulus, only finitely many level curves have
constant speed, including the boundary curves.
\end{lemma}

\begin{proof}
We consider the continuous function
\[
 \mathcal O(s)=\max_\theta S(s,\theta)-\min_\theta S(s,\theta),
 \qquad s\in I.
\]
If $\mathcal O(c)=0$ for some $c\in I$, the curvature $\kappa(c,\theta)$ cannot be constant in view of Lemma \ref{app:den-circle-lemma} so there exist $\theta_0,\theta_1\in \Sone$ such that $\kappa(\theta_1)\neq \kappa(\theta_0)$. The identity $\partial_sS=H$ and
\eqref{app:den-H-curvature} shows
\[
 \left.\frac{d}{ds}
 [S(s,\theta_1)-S(s,\theta_0)]\right|_{s=c} = -2 q\left( \kappa(\theta_1) - \kappa(\theta_0)\right) \ne0.
\]
The difference $S(s,\theta_1) - S(s,\theta_0)$ vanishes at $c$ and is strictly monotone near $c$,
so $c$ is an isolated zero of $\mathcal O$. The same argument
applies from one side when $c=a,b$. Thus, the zero set of
$\mathcal O$ is closed in the compact interval $I$ and all its
points are isolated, so it is finite.
\end{proof}

\subsection{Perturbations of the vorticity law}

We next construct perturbations of the vorticity law whose effect
on the speed cannot be absorbed by changing the level value.
We first consider an interior constant-speed level curve.

\begin{lemma}\label{app:den-interior-response}
Let $\mathsf{C}=\{\psi=c\}$, $a<c<b$, be a constant-speed level curve with
nonconstant curvature $\kappa$. There is a real profile
$G\in C_c^\infty((a,b))$ such that the solution to
\eqref{app:den-linear-response} satisfies
\begin{equation}\label{app:den-transverse-response}
(\mathcal Rv)|_{\mathsf{C}}\notin\operatorname{span}_{\R}\{1,\kappa\}.
\end{equation}
\end{lemma}

\begin{proof}
We first choose a family of profiles $G_\tau$ for which we can
expand $(\mathcal Rv_\tau)|_{\mathsf{C}}$ as $\tau\to\infty$. We then argue
by contradiction and show that, if the conclusion of the lemma
failed, this expansion would force the curvature $\kappa$ to be constant.

We choose $\chi\in C_c^\infty((a,b))$ equal to one near $c$.
For $\tau\geq1$, we consider the real profile
\[
 G_\tau(s)=\chi(s)\cos\left(\tau(s-c)+\frac\pi4\right)
\]
and let $v_\tau$ be the solution to
\[
 Lv_\tau=G_\tau(\psi),\qquad v_\tau|_{\partial\Omega}=0.
\]
We write $\Theta_\tau=\tau(\psi-c)+\pi/4$ 
and note that $\nabla (\cos\Theta_\tau) = -\tau \sin \Theta_\tau \nabla \psi$, so that
\begin{align*}
    \Delta(\cos\Theta_\tau) &= -\tau^2 \cos\Theta_\tau |\nabla\psi|^2 - \tau \sin\Theta_\tau \Delta\psi = -\tau^2 S \cos\Theta_\tau - \tau F(\psi)\sin\Theta_\tau.
\end{align*}
Thus, for any function $A:\overline{\Omega}\rightarrow\R$, we have
\begin{align*}
    \Delta(A\cos\Theta_\tau) &= \Delta A \cos\Theta_\tau + 2\nabla A \cdot \nabla (\cos \Theta_\tau) + A \Delta(\cos\Theta_\tau) \\
    &= (\Delta A -\tau^2 S A) \cos \Theta_\tau - \tau(2\nabla\psi \cdot \nabla A + F(\psi) A) \sin \Theta_\tau,
\end{align*}
and, setting $\mathcal{T}=2D_\psi+F(\psi)\operatorname{Id}$, we obtain
\begin{equation}\label{app:den-conjugation}
 L(A\cos\Theta_\tau)
 =(LA-\tau^2SA)\cos\Theta_\tau-\tau(\mathcal{T}A)\sin\Theta_\tau.
\end{equation}

Now, we set the approximate solution
\begin{equation}\label{eq:def_v_tau_app}
 v_\tau^{\rm app}
 =\sum_{j=2}^7\tau^{-j}a_j
   \cos\left(\Theta_\tau-\frac{j\pi}{2}\right),
\end{equation}
for some $a_j$, $2\leq j \leq 7$, to be determined so that $Lv_\tau^{\rm app} - G_\tau(\psi)$ is small. For $c_j:= \cos\left(\Theta_\tau - \frac{j\pi}{2}\right)$, it follows that 
\begin{align*}
    c_{j+1} = \sin\left(\Theta_\tau - \frac{j\pi}{2}\right), \quad c_{j-2} = -c_j,
\end{align*}
and thus
\begin{align*}
    L(\tau^{-j}a_jc_j) = -\tau^{2-j}S a_j c_j - \tau^{1-j} (\mathcal T a_j) c_{j+1} + \tau^{-j} (La_j) c_j.
\end{align*}
Since $c_2 = -\cos\Theta_\tau$ and $j=2$ corresponds to the zero order, to match the forcing $G_\tau(\psi) = \chi(\psi) \cos\Theta_\tau$ we necessarily set 
$
a_2:= \chi(\psi)/S,
$
which is well-defined since $S>0$. To determine the other coefficients so that they cancel previous errors, we collect all terms with power $\tau^{2-j}$ and phase $c_j$, setting $a_1=0$. The equation for $L(\tau^{-j}a_j c_j)$ gives the term $-\tau^{2-j} S a_j c_j$, the equation for $L(\tau^{-(j-1)}a_{j-1}c_{j-1})$ yields the terms $-\tau^{2-j}(\mathcal T a_{j-1}) c_j$ and the equation for $L(\tau^{-(j-2)}a_{j-2}c_{j-2})$ contributes with $\tau^{2-j}(L a_{j-2}) c_{j-2} = - \tau^{2-j}(La_{j-2}) c_j$. Gathering the three terms, we impose
\begin{align*}
    0 = -\tau^{2-j}\left( Sa_j + \mathcal T a_{j-1} + La_{j-2} \right)c_j
\end{align*}
which prescribes 
\begin{align}\label{app:den-amplitude-recurrence}
    a_j = -\frac{\mathcal T a_{j-1} + L a_{j-2}}{S},
\end{align}
for all $j\geq 3$. Note that each $a_j$ is supported in $\psi^{-1}(\operatorname{supp}\chi)$,
which is a compact subset of $\Omega$. Moreover, we obtain the exact residual formula
\begin{equation}\label{app:den-residual}
 Lv_\tau^{\rm app}-G_\tau(\psi)
 =-\tau^{-6}(\mathcal{T}a_7+La_6)\cos\Theta_\tau
  -\tau^{-7}La_7\sin\Theta_\tau.
\end{equation}
Note that $-\tau^{-6}(\mathcal{T}a_7+La_6)\cos\Theta_\tau$ is the contribution that would be canceled if one further introduced $a_8$. In any case, since 
\begin{align*}
    \Vert \tau^{-6}(\mathcal{T}a_7+La_6)\cos\Theta_\tau
  +\tau^{-7}La_7\sin\Theta_\tau \Vert_{L^p(\Omega)}\leq C\tau^{-6},
\end{align*}
and both $v_\tau$ and $v_\tau^{\rm app}$ vanish on
$\partial\Omega$, using the bounded inverse of $L$ and the embedding
$W^{2,p}(\Omega)\hookrightarrow C^1(\overline\Omega)$, we infer that
\begin{equation}\label{app:den-parametrix-error}
\|v_\tau-v_\tau^{\rm app}\|_{C^1(\overline\Omega)}
\leq C\tau^{-6}.
\end{equation}

Since $\Theta_\tau=\pi/4$ on $\mathsf{C}$, we have the following expansion in $C^0(\mathsf{C})$:
\begin{equation}\label{app:den-response-expansion}
 (\mathcal Rv_\tau)|_{\mathsf{C}} = (\mathcal Rv_\tau^{\rm app})|_{\mathsf{C}} + O(\tau^{-6}) 
 =\frac1{\sqrt2}\left(
 \frac{b_1}{\tau}-\frac{b_2}{\tau^2}-\frac{b_3}{\tau^3}
 +\frac{b_4}{\tau^4}+\frac{b_5}{\tau^5}\right)+O(\tau^{-6}).
\end{equation}
To compute its coefficients, we use the identity
\[
 \mathcal R(A\cos\Theta_\tau)
 =(2D_\psi A-HA)\cos\Theta_\tau-2\tau SA\sin\Theta_\tau.
\]
Applying it to the shifted phases and using
\eqref{app:den-amplitude-recurrence}, we obtain
\begin{equation}\label{app:den-b-recurrence}
b_1=2,\qquad
b_j=-2D_{\psi}a_j-(2F(\psi)+H)a_j-2La_{j-1},\qquad\textup{for }2\leq j\leq5,
\end{equation}
where the right-hand side is restricted to $\mathsf{C}$.
Writing $f=F(c)$, $g=F'(c)$, $h=F''(c)$, and using primes
for arclength derivatives of $\kappa$, we find
\begin{equation}\label{app:den-five-coefficients}
\begin{aligned}
b_1={}&2,\qquad b_2=-\frac{2\kappa}{q},\qquad
b_3=-\frac{2f\kappa}{q^3},\\
b_4={}&\frac{16\kappa''+8\kappa^3+2g\kappa}{q^3}
       -\frac{6f^2\kappa}{q^5},\\
b_5={}&\frac{176f\kappa''+88f\kappa^3+18fg\kappa}{q^5}
       -\frac{96(\kappa')^2+96\kappa\kappa''+32\kappa^4}{q^4}
       -\frac{30f^3\kappa}{q^7}-\frac{2h\kappa}{q^3}\\
       ={}&\leavevmode{\ignorespaces \frac{88(2\kappa'' + \kappa^3) +18g\kappa}{q^5}f - 32\frac{3(\kappa')^2 + 3\kappa \kappa'' + \kappa^4}{q^4} - \frac{30 f^3 \kappa}{q^7} - \frac{2h \kappa}{q^3}}
\end{aligned}
\end{equation}
We verify these formulas in Appendix~\ref{app:coefficient-calculation},
using the first five normal derivatives of $\psi$ along $\mathsf{C}$.

Having this expansion at hand, we now argue by contradiction and assume
that, for every real profile
$G\in C_c^\infty((a,b))$, the corresponding solution $v$ satisfies
\[
 (\mathcal Rv)|_{\mathsf{C}}\in\mathcal V:=\operatorname{span}_{\R}\{1,\kappa\}.
\]
Applying our assumption to $G_\tau$, we obtain
\[
 (\mathcal Rv_\tau)|_{\mathsf{C}}\in\mathcal V,
 \qquad\textup{for all }\tau\geq1.
\]
Since this subspace is closed in $C^0(\mathsf{C})$, we can recover each
coefficient in \eqref{app:den-response-expansion} by successive limits:
we first multiply by $\sqrt2\tau$ and let $\tau\to\infty$, and then
subtract the terms already recovered and multiply by the next power
of $\tau$. The remainder is $O(\tau^{-6})$, so we obtain
$b_1,\ldots,b_5\in\mathcal V$.
Using the formula for $b_4$, we find constants $\alpha,\beta\in\R$
such that
\begin{equation}\label{app:den-curvature-ode}
2\kappa''+\kappa^3=\alpha+\beta \kappa.
\end{equation}
Multiplying by $\kappa'$ and integrating in arclength, we obtain a constant
$C_0\in\R$ such that
\begin{equation}\label{app:den-curvature-first-integral}
(\kappa')^2=\alpha \kappa+\frac\beta2\kappa^2-\frac14\kappa^4+C_0.
\end{equation}
Since $b_5\in \mathcal V$ and $2\kappa''+\kappa^3\in \mathcal V$, we observe that 
\[
3(\kappa')^2+3\kappa\kappa''+\kappa^4\in\operatorname{span}\{1,\kappa\},
\]
while \eqref{app:den-curvature-ode} and
\eqref{app:den-curvature-first-integral} yield
\[
3(\kappa')^2+3\kappa\kappa''+\kappa^4
=3C_0+\frac{9\alpha}{2}\kappa+3\beta \kappa^2-\frac54\kappa^4.
\]
As a result, $-\frac54 \kappa^4 + 3\beta \kappa^2 \in \text{span}_\R \lbrace 1 , \kappa \rbrace$. Since $\kappa$ is continuous and nonconstant on the connected curve $\mathsf{C}$,
its range contains an interval of positive length. On this interval,
the polynomial $3\beta \kappa^2-\frac54\kappa^4$ would agree with an affine
function, giving a contradiction. This proves the lemma.
\end{proof}

We next show that the speed can also be made nonconstant on the
boundary level curves.

\begin{lemma}\label{app:den-boundary-response}
Suppose that a boundary component $\mathsf{C}$ has constant speed and
nonconstant curvature. There is a real profile $G\in C^\infty([a,b])$
such that $D_\psi v$ is nonconstant on $\mathsf{C}$, where $v$ is the solution to
\eqref{app:den-linear-response}.
\end{lemma}

\begin{proof}
Let $c\in\{a,b\}$ be the boundary value on $\mathsf{C}$, and write
$q=|\nabla\psi||_{\mathsf{C}}$ and $n=\nabla\psi/q$ on $\mathsf{C}$. All normal
derivatives in this proof are taken with respect to $n$, which
may point into or out of $\Omega$.

We argue by contradiction and assume that, for every real profile $G\in C^\infty([a,b])$, the
corresponding solution $v$ to \eqref{app:den-linear-response} has
$2D_{\psi}v$ constant on $\mathsf{C}$.

We first consider the solution to
\[
 LZ=0,\qquad Z|_{\mathsf{C}}=1,\qquad
 Z|_{\partial\Omega\setminus \mathsf{C}}=0.
\]
It exists because we can lift the boundary values and apply the
Dirichlet inverse of $L$. Let $\ell$ be the affine function with
$\ell(c)=1$ and value zero at the other endpoint of $[a,b]$. Then,
\[
 L(Z-\ell(\psi))=[F'\ell-\ell'F](\psi),\qquad
 Z-\ell(\psi)=0\quad\text{on }\partial\Omega.
\]
For the real profile $G=F'\ell-\ell'F$, the corresponding solution
is $v=Z-\ell(\psi)$. Since $D_{\psi}=q\partial_n$ and $\partial_n\psi=q$
on $\mathsf{C}$, we have
\[
 2D_{\psi}v=2q\bigl(\partial_nZ-\ell'q\bigr),\qquad\textup{on }\mathsf{C}.
\]
By our contradiction assumption, this expression is constant on $\mathsf{C}$.
Since $q>0$ and $\ell'$ are constants, we infer that
\begin{equation}\label{app:den-constant-dnz}
 \partial_nZ|_{\mathsf{C}}\quad\text{is constant}.
\end{equation}

Now, we choose $\chi\in C^\infty([a,b])$ equal to one near $c$
and zero near the other endpoint. With this cutoff, we repeat the
construction in the proof of Lemma~\ref{app:den-interior-response}: let $v_\tau^{\rm app}$ be given by \eqref{eq:def_v_tau_app} with $a_j$ as in \eqref{app:den-amplitude-recurrence}. Since $a_2=\frac{\chi(\psi)}{S}$, and $S|_{\mathsf{C}}=q^2$, we have ${a_2}|_{\mathsf{C}} = \frac{1}{q^2}$ and
\[
 v_\tau^{\rm app}|_{\mathsf{C}}=-\frac{1}{\sqrt2\tau^2q^2}
                  +O_{C ^{2,\alpha}(\mathsf{C})}(\tau^{-3}).
\]
Note here that $\Theta_\tau|_{\mathsf{C}} = \pi/4$, so any tangential derivative used to estimate the $C ^{2,\alpha}(\mathsf{C})$ norm only acts on the amplitudes ${a_j}|_{\mathsf{C}}$. In particular, $v_\tau^{\rm app}|_{\mathsf{C}}\neq 0$ and is not a valid approximation of $v_\tau$. To remove the trace, let $B_\tau$ be the solution to
\[
 LB_\tau=0,\qquad B_\tau|_{\mathsf{C}}=v_\tau^{\rm app}|_{\mathsf{C}},\qquad
 B_\tau|_{\partial\Omega\setminus \mathsf{C}}=0.
\]
For $B_\tau = B_\tau^{(0)} + B_\tau^{(1)}$ with
\begin{align*}
    LB_\tau^{(0)}&=0,\qquad B_\tau^{(0)}|_{\mathsf{C}}=-\frac{1}{\sqrt{2}\tau^2 q^2},\qquad
 B_\tau^{(0)}|_{\partial\Omega\setminus \mathsf{C}}=0, \\
 LB_\tau^{(1)}&=0,\qquad B_\tau^{(1)}|_{\mathsf{C}}=O_{C^{2,\alpha}(\mathsf{C})}(\tau^{-3}),\qquad
 B_\tau^{(1)}|_{\partial\Omega\setminus \mathsf{C}}=0,
\end{align*}
noting that $B_\tau^{(0)} = -\frac{Z}{\sqrt{2}\tau^2 q^2}$, and using Schauder estimates for $B_\tau^{(1)}$, gives
\begin{equation}\label{app:den-boundary-correction}
 B_\tau=-\frac{Z}{\sqrt2\tau^2q^2}
             +O_{C^{2,\alpha}(\overline\Omega)}(\tau^{-3}).
\end{equation}
Hence, the difference $\widetilde v_\tau^{\rm app}:= v_\tau^{\rm app} - B_\tau$ has zero boundary values and satisfies \eqref{app:den-residual}. As a result, $\Vert v_\tau - \widetilde v_{\tau}^{\rm app} \Vert_{C^1(\overline{\Omega})}\lesssim \tau^{-6}$. Moreover,
\begin{align*}
    2D_\psi \widetilde v_\tau^{\rm app} = (\mathcal R + H \operatorname{Id})\widetilde v_\tau^{\rm app} &= \frac{\sqrt2}{\tau} - \frac{b_2 + H a_2}{\sqrt{2}\tau^2} + \frac{\sqrt{2} \partial_n Z}{q \tau^2} + O_{C^0(\mathsf{C})}(\tau^{-3}) \\
    &= \frac1{\sqrt2}\left[
 \frac2\tau-\frac1{\tau^2}
 \left(\frac{2F(c)}{q^2}-\frac{4\kappa}{q}
                         -\frac{2\partial_nZ}{q}\right)\right]
 +O_{C^0(\mathsf{C})}(\tau^{-3}).
\end{align*}
since the $\tau^{-2}$-coefficient
of $2 D_\psi v_\tau^{\rm app}$ is $-(b_2+Ha_2)/\sqrt2$.
Equations~\eqref{app:den-H-curvature} and
\eqref{app:den-five-coefficients} give
$b_2+Ha_2=2F(c)/q^2-4\kappa/q$. Subtracting the leading term of $B_\tau$ in
\eqref{app:den-boundary-correction} adds
$\sqrt2\partial_nZ/q$ to the coefficient. Consequently, 
\begin{equation}\label{app:den-boundary-expansion}
 2D_\psi v_\tau|_{\mathsf{C}}=\frac1{\sqrt2}\left[
 \frac2\tau-\frac1{\tau^2}
 \left(\frac{2F(c)}{q^2}-\frac{4\kappa}{q}
                         -\frac{2\partial_nZ}{q}\right)\right]
 +O_{C^0(\mathsf{C})}(\tau^{-3}).
\end{equation}

By our assumption, $2D_\psi v_\tau|_{\mathsf{C}}$ is constant on $\mathsf{C}$. Subtracting the leading constant term, multiplying by $\tau^2$, and
letting $\tau\to\infty$, we infer that the second coefficient $\left(\frac{2F(c)}{q^2}-\frac{4\kappa}{q}
                         -\frac{2\partial_nZ}{q}\right)$
is also constant. Since $\partial_n Z|_{\mathsf{C}}$ is constant, this
forces $\kappa$ to be constant, giving a contradiction.
\end{proof}

\begin{proof}[Proof of Proposition~\ref{app:den-density}]
We first perturb $F$ so that the speed is nonconstant on both
boundary components. For each boundary component $\mathsf{C}$ with constant speed, we set
\[
 \mathcal K_{\mathsf{C}}:=\left\{G\in C^\infty([a,b]):
 (D_\psi L^{-1}(G(\psi)))|_{\mathsf{C}}\in\operatorname{span}_{\R}\{1\}\right\}.
\]
By Lemma~\ref{app:den-boundary-response}, each $\mathcal K_{\mathsf{C}}$ is a
proper linear subspace of $C^\infty([a,b])$. Since their finite union
does not cover this space, there exists a real profile
\[
 G\in C^\infty([a,b])\setminus\bigcup_{\mathsf{C}}\mathcal K_{\mathsf{C}}.
\]
Let $v=L^{-1}(G(\psi))$, with $D_\psi v|_{\mathsf{C}}$ non-constant. Hence, for each such $\mathsf{C}$, there exist $x_{\mathsf{C}},y_{\mathsf{C}}\in \mathsf{C}$ with
\[
 d_{\mathsf{C}}:=2D_\psi v(x_{\mathsf{C}})-2D_\psi v(y_{\mathsf{C}})\ne0.
\]
Since $(\partial_\eps S_\eps)|_{\eps=0} = \mathcal R v \circ \Phi$, with $\mathcal R v = 2D_\psi v - Hv$ and $v|_{\partial\Omega}=0$, for $0<|\eps| \ll 1$,
\[
 |\nabla\psi_\eps(x_{\mathsf{C}})|^2-|\nabla\psi_\eps(y_{\mathsf{C}})|^2
 =\eps d_{\mathsf{C}}+o(\eps)\ne0.
\]
On each boundary component $\mathsf{C}$ where the original speed is nonconstant,
uniform convergence gives, for sufficiently small $|\eps|$,
\[
 \osc_{\mathsf{C}}|\nabla\psi_\eps|^2
 \geq\frac12\osc_{\mathsf{C}}|\nabla\psi|^2>0.
\]
We choose this perturbation sufficiently small to preserve the
non-stagnation and invertibility assumptions of
Proposition~\ref{app:den-density}, as justified below.
We have thus obtained an arbitrarily close smooth pair, which
we still denote by $(F,\psi)$, with the required boundary property.

By Lemma~\ref{app:den-finite-levels} and continuity in the
parametrization \eqref{app:parametrization}, we have
\[
 \{s\in I:\osc_\theta S(s,\theta)=0\}
 =\{c_1,\ldots,c_M\}\subset(a,b).
\]
If $M=0$, compactness gives \eqref{app:den-positive-oscillation}.
Otherwise, we set $\mathsf{C}_i:=\{\psi=c_i\}$ and write $\kappa_i$ for its
signed curvature with respect to $n=\nabla\psi/|\nabla\psi|$.
By Lemma~\ref{app:den-interior-response}, the linear subspaces
\[
 \mathcal K_i:=\left\{G\in C_c^\infty((a,b)):
 (\mathcal R L^{-1}(G(\psi)))|_{\mathsf{C}_i}
 \in\operatorname{span}_{\R}\{1,\kappa_i\}\right\}
\]
are proper. As before, there exists a real profile $G\in C_c^\infty((a,b))\setminus\bigcup_{i=1}^M\mathcal K_i$ such that for $v=L^{-1}(G(\psi))$, there holds $(\mathcal R L^{-1}(G(\psi)))|_{\mathsf{C}_i}
 \not\in\operatorname{span}_{\R}\{1,\kappa_i\}$ for all $1\leq i \leq M$.

Let $1\leq i\leq M$. By
\eqref{app:den-H-curvature}, $H|_{\mathsf{C}_i} = 2(F(c_i) - q_i \kappa_i)$ so that $\operatorname{span}_{\R}\{1,H\}
 =\operatorname{span}_{\R}\{1,\kappa_i\}$. Given our choice of $G$, we conclude 
\[
\mathcal Rv\notin\operatorname{span}_{\R}\{1,H\}.
\]
We next identify these functions with their pullbacks by $\Phi$. As $H|_{\mathsf{C}_i}$ is non-constant (otherwise we would have constant curvature $\kappa|_{\mathsf{C}_i}$), there exist $\theta_0,\theta_1\in \Sone$ such that $H(c_i,\theta_1) - H(c_i,\theta_0)\neq 0$. Moreover, $\mathcal R v\notin\operatorname{span}_{\R}\{1,H\}$ shows the existence of some $\theta_2\in \Sone$ such that
\begin{align*}
    \mathcal Rv(c_i,\theta_2) \neq \mathcal Rv (c_i,\theta_0) + \frac{\mathcal R v(c_i, \theta_1) - \mathcal{R} v(c_i,\theta_0)}{H(c_i,\theta_1) - H(c_i,\theta_0)} H(c_i,\theta_2) - \frac{\mathcal R v(c_i, \theta_1) - \mathcal{R} v(c_i,\theta_0)}{H(c_i,\theta_1) - H(c_i,\theta_0)} H(c_i,\theta_0).
\end{align*}
In particular, 
\begin{align*}
    \text{det} \begin{pmatrix}
        H(c_i,\theta_1) - H(c_i,\theta_0) & H(c_i,\theta_2) - H(c_i,\theta_0) \\
        \mathcal R v(c_i, \theta_1) - \mathcal{R} v(c_i,\theta_0) & \mathcal R v(c_i, \theta_2) - \mathcal{R} v(c_i,\theta_0)
    \end{pmatrix} \neq 0.
\end{align*}

We now consider the smooth map
\[
 \mathcal A_i(s,\eps):=
 \begin{pmatrix}
 S_\eps(s,\theta_1)-S_\eps(s,\theta_0)\\
 S_\eps(s,\theta_2)-S_\eps(s,\theta_0)
 \end{pmatrix}
\]
and we note that $\mathcal A_i(c_i,0) = 0$. Moreover, as $(\partial_\eps S_\eps)|_{\eps=0} = \mathcal{R}v \circ \Phi$ and $(\partial_s S_\eps)|_{\eps=0}=H$, it follows that
\[
 \det D_{(s,\eps)}\mathcal A_i(c_i,0) = \det \begin{pmatrix}
        H(c_i,\theta_1) - H(c_i,\theta_0) & H(c_i,\theta_2) - H(c_i,\theta_0) \\
        \mathcal R v(c_i, \theta_1) - \mathcal{R} v(c_i,\theta_0) & \mathcal R v(c_i, \theta_2) - \mathcal{R} v(c_i,\theta_0)
    \end{pmatrix} \neq 0.
\]
The inverse function theorem then gives $\eta_i,\eps_i>0$ such that,
for $J_i=(c_i-\eta_i,c_i+\eta_i)\Subset(a,b)$, we have
\[
 \mathcal A_i^{-1}(\{0\})
 \cap\bigl(J_i\times(-\eps_i,\eps_i)\bigr)=\{(c_i,0)\}.
\]
Therefore, $\osc_\theta S_\eps(s,\theta)>0$ for all $s\in J_i$ and all $0<|\eps|<\eps_i$. On the remaining compact set, we have
\[
 \mu:=\min_{s\in I\setminus\bigcup_{i=1}^M J_i}
              \osc_\theta S(s,\theta)>0.
\]
Since $S_\eps\to S$ uniformly on $I\times\Sone$, for sufficiently
small $|\eps|$ and all $s\in I\setminus\bigcup_{i=1}^M J_i$,
\[
 \osc_\theta S_\eps(s,\theta)
 \geq\mu-2\|S_\eps-S\|_{L^\infty(I\times\Sone)}
 \geq\frac\mu2>0.
\]
We now take $0<|\eps|<\min_{1\leq i\leq M}\eps_i$ sufficiently
small that this estimate holds. Continuity and compactness then give
\[
 \min_{s\in I}\osc_{\{\psi_\eps=s\}}|\nabla\psi_\eps|^2
 =\min_{s\in I}\osc_\theta S_\eps(s,\theta)>0,
\]
which is \eqref{app:den-positive-oscillation}.

For each perturbation, we preserve the absence of critical points by
choosing its size so that
\[
 \begin{aligned}
 \min_{\overline\Omega}|\nabla\psi_\eps|
 &\geq\min_{\overline\Omega}|\nabla\psi|
       -\|\nabla\psi_\eps-\nabla\psi\|_{L^\infty} \geq\frac12\min_{\overline\Omega}|\nabla\psi|>0.
 \end{aligned}
\]
To check invertibility, we write $L_\eps=L-M_\eps$ and $M_\eps=\bigl(F_\eps'(\psi_\eps)-F'(\psi)\bigr)\operatorname{Id}$. For every $w\in W^{2,p}(\Omega)\cap W^{1,p}_0(\Omega)$ and
sufficiently small $|\eps|$, we have
\[
 \begin{aligned}
 \|L^{-1}M_\eps w\|_{W^{2,p}}
 &\leq C\|F_\eps'(\psi_\eps)-F'(\psi)\|_{L^\infty}
             \|w\|_{W^{2,p}} \leq\frac12\|w\|_{W^{2,p}}.
 \end{aligned}
\]
The Neumann series and the identity $ L_\eps=L(\operatorname{Id}-L^{-1}M_\eps)$ therefore give the required invertibility. Both properties persist
in a fixed neighborhood of the original pair, where we keep both
perturbations.

Finally, we return to the original notation $(F,\psi)$. For each
$j$, we choose each of the two changes smaller than $1/(2j)$ in
the sum of the corresponding $C^r$ norms, for $0\leq r\leq j$.
By the triangle inequality, the resulting pair $(F_j,\psi_j)$ satisfies
\[
 \|F_j-F\|_{C^r(I)}+\|\psi_j-\psi\|_{C^r}
 <\frac1{2j}+\frac1{2j}=\frac1j,
 \qquad\textup{for all }0\leq r\leq j.
\]
This proves \eqref{app:den-smooth-convergence} and completes the proof.
\end{proof}

\subsection{The circular case}

We now prove the converse. The main point is to show that every steady
Euler flow sufficiently close to a smooth circular Euler flow without
stagnation points is also circular. We use polar coordinates
$r=|x|$ and $\theta\in \mathbb S^1$.

\begin{proposition}
\label{app:cir-local-rigidity}
Let
\[
 \Omega=\{x\in\R^2:r_-<|x|<r_+\},\qquad 0<r_-<r_+,
\]
and let $u_0=V_0(r)e_\theta$ be a smooth circular Euler flow with
$\min_{[r_-,r_+]}|V_0|>0$. We set $\omega_0=\curl u_0$.
There exists $\eta>0$ such that every steady Euler flow
$v\in\mathcal Y(\Omega)$ satisfying
\begin{equation}\label{app:cir-neighborhood}
 \|v-u_0\|_{\mathcal Y}<\eta,
\end{equation}
is circular.
\end{proposition}

\begin{proof}
We first normalize the boundary values of the stream function of $v$.
Since $V_0$ never vanishes, we can assume that $V_0>0$, replacing
$u_0$ and $v$ by their negatives if necessary. This preserves the
steady Euler equations and the distance~\eqref{app:cir-neighborhood}.
We choose a radial stream function $\psi_0=\Psi_0(r)$ with
$\Psi_0'=V_0$ and set
\[
 a=\Psi_0(r_-),\qquad b=\Psi_0(r_+),\qquad a<b.
\]
Then $\Delta\psi_0=F_0(\psi_0)$, where
\[
 F_0(s)=\left(V_0'(r)+\frac{V_0(r)}r\right)_{r=\Psi_0^{-1}(s)}
 \in C^\infty([a,b]).
\]
By Lemma~\ref{rig:streamfunction}, we write $v=\nabla^\perp\phi$,
with $\phi=a$ on the inner boundary, and denote its outer boundary
value by $\widetilde b$. For a fixed $p>2$,
Lemma~\ref{rig:streamfunction}, applied to $v-u_0$, gives
\[
 \|v-u_0\|_{W^{1,p}}\leq C\|v-u_0\|_{\mathcal Y}.
\]
Since $\phi-\psi_0$ vanishes on the inner boundary, Poincar\'e's
inequality and the embedding
$W^{2,p}(\Omega)\hookrightarrow C^1(\overline\Omega)$ imply
\begin{equation}\label{app:cir-stream-close}
 \|\phi-\psi_0\|_{W^{2,p}}+
 \|\phi-\psi_0\|_{C^1}
 \leq C\|v-u_0\|_{\mathcal Y},
 \qquad |\widetilde b-b|\leq C\|v-u_0\|_{\mathcal Y}.
\end{equation}
For $\eta$ sufficiently small, we have $\widetilde b>a$ and can define
\begin{equation}\label{app:cir-normalization}
 \lambda=\frac{b-a}{\widetilde b-a},\qquad
 \psi=a+\lambda(\phi-a),\qquad u=\nabla^\perp\psi=\lambda v.
\end{equation}
Thus $\psi=a,b$ on the inner and outer circles, respectively, and
\begin{equation}\label{app:cir-normalized-close}
 |\lambda-1|+\|\psi-\psi_0\|_{C^1}
 +\|\lambda\curl v-\omega_0\|_{L^\infty}
 \leq C\|v-u_0\|_{\mathcal Y}.
\end{equation}
The velocity $u$ is a steady Euler flow, with pressure equal to
$\lambda^2$ times that of $v$. Moreover, after decreasing $\eta$,
we have
\begin{equation}\label{app:cir-radial-monotonicity}
 \partial_r\psi\geq\tfrac12\min_{[r_-,r_+]}V_0>0
 \qquad\text{on }\overline\Omega.
\end{equation}

We next express the vorticity as a bounded measurable function of
$\psi$. For every $\theta$, the function $r\mapsto\psi(r,\theta)$
is strictly increasing from $a$ to $b$. We denote its inverse by
$R(s,\theta)$, so that $\psi(R(s,\theta),\theta)=s$ and
\begin{equation}\label{app:cir-R-boundary}
 R(a,\theta)=r_-,\qquad R(b,\theta)=r_+,\qquad \partial_s R>0.
\end{equation}
We set $\mathsf{D}=(a,b)\times\Sone$. By the inverse function theorem,
$R\in C^1(\overline{\mathsf{D}})$, and
\[
 \partial_s R=\frac1{\partial_r\psi(R,\theta)},\qquad
 \partial_\theta R=-\frac{\partial_\theta\psi(R,\theta)}{\partial_r\psi(R,\theta)}.
\]
These identities and the weak chain rule give $R\in W^{2,p}(\mathsf{D})$
for every $1<p<\infty$: each second derivative is a linear
combination of second derivatives of $\psi$, composed with the
parametrization, with bounded coefficients. The change of variables
maps $L^p(\Omega)$ continuously into $L^p(\mathsf{D})$ because its Jacobian
and reciprocal are bounded. The same identities give
\begin{equation}\label{app:cir-R-close}
 \|R-R_0\|_{C^1}\leq C\|v-u_0\|_{\mathcal Y},
 \qquad R_0=\Psi_0^{-1}.
\end{equation}
More explicitly, the parametrization is the $C^1$ diffeomorphism
\begin{equation}\label{app:cir-parametrization}
 \begin{aligned}
 X:[a,b]\times\Sone&\longrightarrow\overline\Omega,\\
 X(s,\theta)&=R(s,\theta)e_r(\theta),
 \qquad\psi(X(s,\theta))=s.
 \end{aligned}
\end{equation}
Its area Jacobian is $R\partial_s R>0$. We write $\omega=\curl u$ and,
for $\xi\in C_c^1(\mathsf{D})$, take $\zeta=\xi\circ X^{-1}$ in
\eqref{rig:transport}. Since $u\cdot\nabla\psi=0$ and
$u\cdot\nabla\theta=\partial_r\psi/r$, we obtain $\partial_\theta \xi = \partial_\theta X\cdot \nabla \zeta(X)$. Likewise, since $\psi(X(s,\theta))=s$ we deduce that $0 = \partial_\theta X \cdot \nabla \psi(X)$ which yields $\partial_\theta X = \lambda(X) u(X)$ for some $\lambda(X)$ to be determined. Moreover, as $\theta(X(s,\theta))=\theta$, we observe that 
\begin{align*}
    \lambda(X)\frac{1}{R\partial_s R}=\lambda(X) \frac{\partial_r\psi}{r}=\lambda(X) u(X)\cdot \nabla \theta(X) = \partial_\theta X \cdot \nabla \theta(X) = 1,
\end{align*}
which in turn gives $\lambda(X) = R\partial_s R>0$. Hence, it follows that $\partial_\theta\xi =\lambda(X) u(X) \cdot \nabla \zeta(X)$ and thus $(u\cdot\nabla\zeta)(X)=\frac{\partial_\theta\xi}{R\partial_s R}$. In particular, we have that
\begin{align*}
    0 = \int_\Omega \omega(x) (u\cdot \nabla \zeta)(x) \dd x = \int_{\mathsf{D}}\omega(X(s,\theta))\partial_\theta\xi(s,\theta)\,\dd s\,\dd\theta.
\end{align*}
Lemma~\ref{rig:zero-derivative} therefore yields a bounded measurable
function $F:(a,b)\to\R$ such that
\begin{equation}\label{app:cir-measurable-profile}
 (\omega\circ X)(s,\theta)=F(s)\quad\text{a.e. in }\mathsf{D},
 \qquad\Delta\psi=F(\psi)\quad\text{a.e. in }\Omega.
\end{equation}
Here we have used that $X$ and $X^{-1}$ preserve null sets.
Moreover, for a.e. $(s,\theta)\in \mathsf{D}$, we have
\[
 |F(s)-F_0(s)|
 \leq\|\lambda\curl v-\omega_0\|_{L^\infty}
       +\|F_0'\|_{L^\infty}\|\psi-\psi_0\|_{L^\infty}.
\]
Consequently, \eqref{app:cir-normalized-close} gives
\begin{equation}\label{app:cir-F-close}
 \|F-F_0\|_{L^\infty}\leq C\|v-u_0\|_{\mathcal Y}.
\end{equation}

At this point, let us remark that it suffices to prove that $\partial_\theta R=0$, since then $\psi$
is radial and $v$ is circular by~\eqref{app:cir-normalization}.
Having this representation at hand, we derive an equation for
$\partial_\theta R$ that does not require differentiating $F$. 

Let $\mathsf{F}(s) = \int_a^s F(t) \dd t$ be a primitive of $F$, we define the energy functional
\begin{align*}
    \mathsf E[\varphi] = \int_\Omega\frac12 |\nabla\varphi|^2 \dd x  +\int_\Omega \mathsf F(\varphi) \dd x
\end{align*}
for all $\varphi\in H^1(\Omega)$ with $\varphi(r_-,\theta)=a$ and $\varphi(r_+,\theta) = b$. Let $\phi\in C^\infty_0(\Omega)$ and $\psi_\eps = \psi + \eps \phi$. For $0<|\eps| \ll 1$, $\partial_r\psi_\eps >0$ on $\overline{\Omega}$ and there exists $R_\eps(s,\theta)$ such that 
\begin{align*}
    \psi_\eps(R_\eps(s,\theta),\theta) =s
\end{align*}
for all $(s,\theta)\in [a,b]\times\Sone$. In particular,
\[
 \partial_s R_\eps=\frac1{\partial_r\psi_\eps(R_\eps,\theta)},\qquad
 \partial_\theta R_\eps=-\frac{\partial_\theta\psi_\eps(R_\eps,\theta)}{\partial_r\psi_\eps(R_\eps,\theta)},
\]
and the area Jacobian is $R_\eps \partial_s R_\eps >0$. We next note that
\begin{align*}
     \frac12 \int_\Omega |\nabla \psi_\eps|^2 \dd x & = \frac12 \int_\Omega \left( |\partial_r\psi_\eps|^2 + \frac{|\partial_\theta\psi_\eps|^2}{r^2}\right) \dd x  \\& = \frac12 \int_{\mathsf{D}} \left( \frac{1}{(\partial_s R_\eps)^2} + \frac{1}{R_\eps^2}\frac{(\partial_\theta R_\eps)^2}{(\partial_s R_\eps)^2} \right) R_\eps \partial_s R_\eps \dd s \dd \theta = \frac{1}{2} \int_{\mathsf{D}} \frac{R_\eps + \frac{(\partial_\theta R_\eps)^2}{R_\eps}}{\partial_s R_\eps} \dd s \dd \theta,
\end{align*}
while we also have
\begin{align*}
    \int_\Omega \mathsf F(\psi_\eps) \dd x = \int_{\mathsf{D}} \mathsf{F}(s) R_\eps \partial_s R_\eps \dd s \dd \theta = \pi \mathsf F(b) (r_+)^2 - \frac12 \int_{\mathsf{D}} F(s) R_\eps^2 \dd s \dd \theta.
\end{align*}
Therefore, for $T(s,\sigma,\rho,z) := \frac{\sigma + z^2/\sigma}{2\rho} - \frac12 F(s) \sigma^2$, with $\sigma,\rho>0$, it follows that
\begin{align}\label{eq:equiv_mathsfE_T}
    \mathsf E[\psi_\eps] = \pi \mathsf F(b) (r_+)^2 + \int_{\mathsf{D}} T(s, R_\eps, \partial_s R_\eps, \partial_\theta R_\eps) \dd s \dd \theta.
\end{align}
Since $\Delta\psi = F(\psi)$ a.e. in $\Omega$, there holds
\begin{align}
    \frac{\dd }{\dd \eps}\Big|_{\eps=0} \mathsf E[\psi_\eps] = \int_\Omega \left( \nabla \psi \cdot \nabla \phi + F(\psi)\phi \right) \dd x = 0.
\end{align}
We next define
\begin{align*}
    \mathcal J_F[R] := \int_{\mathsf{D}} T(s,R, \partial_s R, \partial_\theta R) \dd s \dd \theta,
\end{align*}
and note that, in view of \eqref{eq:equiv_mathsfE_T}, we reach
\begin{align*}
    0 =  \frac{\dd }{\dd \eps}\Big|_{\eps=0} \mathcal{J}_F[R_\eps] &=  \int_{\mathsf{D}} \left( \partial_\sigma T \vartheta + \partial_\rho T \partial_s\vartheta + \partial_z T \partial_\theta\vartheta\right) \dd s \dd \theta =: \mathcal{J}_F'[R]\vartheta,
\end{align*}
for $\vartheta:=(\partial_\eps R_\eps)|_{\eps=0} = -\partial_s R(s,\theta) \phi(X(s,\theta))\in H^1_0(\mathsf{D})$ and for all $\phi\in C_0^\infty(\Omega)$. 

We claim that $\mathcal J'_F[R]w = 0$ for all $w\in H_0^1(\mathsf{D})$. Indeed, for $w\in H_0^1(\mathsf{D})$ we define
\begin{align}
    \phi(x) = - \left( \frac{ w}{\partial_s R}\right) \circ X^{-1}(x)
\end{align}
and we note that $\phi\in H_0^1(\Omega)$ because $X$ is bi-Lipschitz, $\partial_s R$ is bounded away from zero and $\partial_s R, (\partial_s R)^{-1}\in W^{1,p}(\mathsf{D})$ for $p>2$. By density of $C^\infty_0(\Omega)$ in $H_0^1(\Omega)$, there exist $\phi_n\in C^\infty_0(\Omega)$ with $\phi_n \rightarrow \phi$ in $H_0^1(\Omega)$ such that $\mathcal J'_F[R]w_n = 0$, for $w_n := - \partial_s R \phi_n(X)\rightarrow -\partial_s R\phi(X)=w$ in $H^1_0(\mathsf{D})$. By continuity, $\mathcal{J}_F'[R]w = 0$.

Next, we consider the $\theta$-average of $R(s,\theta)$ and define
\begin{align*}
    w (s,\theta):= R(s,\theta)- \overline R(s), \quad \overline R(s) := \frac{1}{2\pi}\int_{\Sone} R(s,\theta) \dd \theta.
\end{align*}
It is mean-free and since $R(a,\theta) = r_-$ and $R(b,\theta) = r_+$ for all $\theta\in \Sone$ we have $w\in H_0^1(\mathsf{D})$ so that $\mathcal J'_F[R]w = 0$. Additionally, since $w$ is mean-free, we have that
\begin{align}\label{eq:stationarity_overline_R}
    \mathcal J'_F[\overline R]w = 0.
\end{align}
Indeed, $\partial_s w$ is also mean-free and $\partial_\sigma T(s, \overline{R}, \partial_s \overline{R}, 0)$, $\partial_\rho T(s, \overline{R}, \partial_s \overline{R}, 0)$ depend only on $s\in[a,b]$, while $\partial_{z} T(s, \overline{R}, \partial_s \overline{R}, 0) = 0$. {In fact, \eqref{eq:stationarity_overline_R} holds because $\overline R$ describes a radial flow and mean free variations $w$ preserve every sub-level area to first order. In other words, $\mathcal J_F'[\overline R]w=0$ is equivalent to the radial stream-function defined by $\overline R$ being stationary for the kinetic energy under area-preserving rearrangements. This is the case since radial stream-functions are steady solutions of Euler and thus stationary for the kinetic energy under such variations. Note moreover that the potential term depends only on the stream-function distribution, which is preserved under area-preserving rearrangements.}

Next, we consider the linear interpolation $\widetilde R_t = \overline R + t w$ for $0 \leq t \leq 1$, with $\Vert \widetilde R_t - R_0 \Vert_{C^1}\leq C \eta$. Then,
\begin{align*}
    0 = \mathcal J_F'[R]w - \mathcal J_F'[\overline R]w = \int_0^1 \mathcal{J}_F''[\widetilde R_t](w,w) \dd t.
\end{align*}
To compute the second variation, we note that 
\begin{align*}
    \mathcal J''_F[\widetilde R_t](w,w) &= \frac{\dd^2}{\dd \epsilon^2}\Big|_{\epsilon=0} \int_{\mathsf{D}} T(s,\widetilde R_t + \epsilon w, \partial_s(\widetilde R_t + \epsilon w), \partial_\theta(\widetilde R_t + \epsilon w)) \dd s \dd \theta \\
    &= \int_{\mathsf{D}} (w, \partial_s w, \partial_\theta w) D^2_{(\sigma, \rho, z)}T(s, \widetilde R_t, \partial_s \widetilde R_t, \partial_\theta\widetilde R_t) (w, \partial_s w, \partial_\theta w)^T \dd s \dd \theta, 
\end{align*}
where
\begin{align*}
    D^2_{(\sigma, \rho, z)}T(s,\sigma,\rho, z) = \begin{pmatrix}
        \frac{z^2}{\sigma^3 \rho} - F(s) & -\frac{1-z^2/\sigma^2}{2\rho^2} & -\frac{z}{\sigma^2\rho} \\
         -\frac{1-z^2/\sigma^2}{2\rho^2} &\frac{\sigma + z^2/\sigma}{\rho^3} & -\frac{z}{\sigma\rho^2} \\
         -\frac{z}{\sigma^2\rho} & -\frac{z}{\sigma\rho^2}& \frac{1}{\sigma\rho}
    \end{pmatrix}.
\end{align*}
Importantly, for $R_0$ we have $\rho_0=\partial_s R_0 \neq 0$, $z_0=\partial_\theta R_0=0$ and $F_0 = -\frac{\rho_0'}{\rho_0^3} + \frac{1}{R_0 \rho_0}$. In particular, we define
\begin{align*}
   Q_0[w]:= \mathcal J_{F_0}''[R_0](w,w) = \int_{\mathsf{D}} \left( \frac{R_0}{\rho_0^3}(\partial_s w)^2 + \frac{1}{R_0 \rho_0}\left( (\partial_\theta w)^2 - w^2 \right) \right) \dd s \dd \theta.
\end{align*}
Since $w$ is mean-free, Poincar\'{e} inequality yields $\Vert w(s,\cdot) \Vert_{L^2(\Sone)}^2 \leq \Vert \partial_\theta w(s,\cdot)\Vert_{L^2(\Sone)}^2$ for all $s\in[a,b]$.  Hence, the second contribution to $Q_0[w]$ is non-negative and positive lower bounds of $R_0$, $\rho_0$ and $w \in H_0^1(\mathsf{D})$ yield
\begin{align*}
    \|\partial_s w\|_{L^2}^2\leq C Q_0[w],\qquad
 \|w\|_{L^2}^2\leq C\|\partial_s w\|_{L^2}^2\leq C Q_0[w].
\end{align*}
Additionally, for $C_0:= \frac{1}{R_0\rho_0}$,
\[
 \int_{\mathsf{D}} C_0(\partial_\theta w)^2\,\dd s\,\dd\theta
 \leq Q_0[w]+\|C_0\|_{L^\infty}\|w\|_{L^2}^2.
\]
Thus, there exists $c_*>0$, depending only on the reference Euler
flow and the annular domain, such that
\begin{equation}\label{app:cir-coercivity}
 Q_0[w]\geq c_*\bigl(\|\partial_s w\|_{L^2}^2+
 \|\partial_\theta w\|_{L^2}^2+\|w\|_{L^2}^2\bigr).
\end{equation}
By \eqref{app:cir-R-close} and \eqref{app:cir-F-close}, 
\begin{align*}
     |\mathcal J_F''[\widetilde R_t](w,w)-Q_0[w]|&\leq C\left( \Vert F - F_0\Vert_{L^\infty} + \Vert R- R_0 \Vert_{C^1}\right)
 \bigl(\|\partial_s w\|_{L^2}^2+\|\partial_\theta w\|_{L^2}^2+
 \|w\|_{L^2}^2\bigr) \\
 &\leq C\|v-u_0\|_{\mathcal Y}
 \bigl(\|\partial_s w\|_{L^2}^2+\|\partial_\theta w\|_{L^2}^2+
 \|w\|_{L^2}^2\bigr) \\
 &\leq C \eta 
 \bigl(\|\partial_s w\|_{L^2}^2+\|\partial_\theta w\|_{L^2}^2+
 \|w\|_{L^2}^2\bigr)
\end{align*}
uniformly for all $t\in[0,1]$. As a result,
\begin{align*}
  &  c_*\bigl(\|\partial_s w\|_{L^2}^2+
 \|\partial_\theta w\|_{L^2}^2+\|w\|_{L^2}^2\bigr)  \leq  Q_0[w] \\
 & \qquad = \int_0^1 Q_0[w] \dd t = -\int_0^1 (\mathcal{J}_F''[\widetilde R_t](w,w) - Q_0[w]) \dd t \leq C \eta \bigl(\|\partial_s w\|_{L^2}^2+
 \|\partial_\theta w\|_{L^2}^2+\|w\|_{L^2}^2\bigr),
\end{align*}
reaching a contradiction for $\eta$ small enough unless $w=0$. Hence $R(s,\theta) = \overline R(s)$ is radial, $\psi$ is radial and $u$ is circular. By \eqref{app:cir-normalization}, $v$ is circular as well.
\end{proof}

\begin{lemma}
\label{app:cir-nonisolation}
Let $\Omega=\{x\in\R^2:r_-<|x|<r_+\}$ be an annulus, with
$0<r_-<r_+$.
Every nonzero circular steady Euler flow $v=V(r)e_\theta$ in
$\mathcal Y(\Omega)$ is a limit of distinct circular steady Euler
flows of the same kinetic energy, with convergence in $\mathcal Y(\Omega)$
and in $W^{1,p}(\Omega;\R^2)$ for every $1<p<\infty$.
In particular, $v$ is not
isolated on its energy surface.
\end{lemma}

\begin{proof}
We use the classical family of circular Euler flows discussed by
Hamel and Nadirashvili in~\cite[Section~1.1]{HamelNadirashvili2023}.
By Lemma~\ref{rig:streamfunction}, $V$ is bounded. Writing
$\omega=\curl v$, the identity $V'=\omega-V/r$ then gives
$V\in W^{1,\infty}(r_-,r_+)$.
We choose $h\in C_c^\infty(r_-,r_+)$ not proportional to $V$ and,
for $|t|$ sufficiently small, set
\[
 c_t=\frac{\|v\|_{L^2}}{\|(V+th)e_\theta\|_{L^2}},
 \qquad v_t=c_t(V+th)e_\theta.
\]
Each $v_t$ is a circular steady Euler flow, with pressure satisfying
$\partial_rP_t=c_t^2(V+th)^2/r$. By construction, $E(v_t)=E(v)$
and $v_t\notin\R v$ for $t\ne0$, since $h$ is not proportional to $V$.
Finally, $c_t\to1$ gives $v_t\to v$ in
$W^{1,\infty}(\Omega;\R^2)$, which yields both asserted convergences.
\end{proof}

\begin{proof}[Proof of Theorem~\ref{intro:approximation}]
We first suppose that $\Omega$ is not a circular annulus.
By Proposition~\ref{app:den-density}, there are smooth pairs
$(F_j,\Psi_j)$ with the original boundary values such that
{
\[
 F_j\longrightarrow F\quad\text{in }C^\infty(I),\qquad
 \Psi_j\longrightarrow\psi\quad\text{in }C^\infty(\overline\Omega),
 \qquad \Delta\Psi_j=F_j(\Psi_j),
\]
}
and
\[
 \mathbf{o}_j:=\min_{s\in I}\osc_{\{\Psi_j=s\}}|\nabla\Psi_j|^2>0.
\]
Each $\Psi_j$ has no critical points and has an invertible
Dirichlet linearization. We do not need a uniform lower bound
for $\mathbf{o}_j$.

We apply Proposition~\ref{con:rough-approximation} to each pair,
using the finitely valued Borel law obtained from $F_j$ by
\eqref{con:profile}. We choose $\eta_j>0$ sufficiently small
that the resulting stream function $\widetilde\psi_j$, with
$\widetilde u_j=\nabla^\perp\widetilde\psi_j$, satisfies
\begin{equation}\label{app:diagonal-choice}
 \begin{aligned}
 &\|\widetilde\psi_j-\Psi_j\|_{W^{2,k}}<\frac1j
 \qquad\textup{for every integer }2\leq k\leq j+2,\\
 &\|\widetilde u_j-\nabla^\perp\Psi_j\|_{\mathcal Y}<\frac1j.
 \end{aligned}
\end{equation}
For each $j$, there are only finitely many estimates to impose.
The velocity $\widetilde u_j$ is a steady Euler flow in
$\mathcal Y(\Omega)$ with finitely valued vorticity.
Proposition~\ref{con:rough-approximation} also gives a $C^1$
parametrization of its level curves and a positive lower bound
for the speed oscillation. Since its Borel vorticity law is not
essentially constant on any open interval, Lemma~\ref{rig:rigidity}
shows that $\widetilde u_j$ is isolated modulo scaling.

For any fixed $1<r<\infty$, we choose an integer
$k\geq\max\{2,r\}$. By \eqref{app:diagonal-choice} and the
inclusion $W^{2,k}(\Omega)\subset W^{2,r}(\Omega)$, we obtain
\[
 \widetilde u_j\longrightarrow u\quad\text{in }W^{1,r}(\Omega),
 \qquad
 \|\widetilde u_j-u\|_{\mathcal Y}\longrightarrow0.
\]
Here we have also used the uniform convergence
$F_j(\Psi_j)\to F(\psi)$. Finally, since $u\ne0$, we set
\[
 c_j=\frac{\|u\|_{L^2}}{\|\widetilde u_j\|_{L^2}},\qquad
 u_j=c_j\widetilde u_j.
\]
Then $c_j\to1$, $E(u_j)=E(u)$, and the same convergence holds
for $u_j$. Multiplication by $c_j$ preserves stationarity, with
pressure multiplied by $c_j^2$, and preserves isolation modulo
scaling: we can divide any nearby steady Euler flow by $c_j$
and apply the isolation of $\widetilde u_j$.
This proves the approximation assertion.

We now assume that
$\Omega=\{r_-<|x-x_*|<r_+\}$ and consider
\[
 w=(x-x_*)^\perp\cdot\nabla\psi.
\]
This function vanishes on both boundary circles. Since the Laplacian
commutes with rotations, differentiation of the equation for $\psi$
gives $Lw=0$. The Dirichlet kernel is trivial, so $w=0$ and $u$
is circular. Moreover, $u$ has no stagnation points by assumption.

Proposition~\ref{app:cir-local-rigidity} gives a neighborhood of $u$
in the norm $\|\cdot\|_{\mathcal Y}$ whose steady Euler flows are circular.
We can shrink this
neighborhood so that they are all nonzero.
Lemma~\ref{app:cir-nonisolation} then shows that none is isolated
on its energy surface. Consequently, the sequence of isolated Euler
flows in the statement cannot exist on a concentric circular annulus.
\end{proof}

We conclude the section by proving Corollary~\ref{intro:harmonic-reference}.

\begin{proof}[Proof of Corollary~\ref{intro:harmonic-reference}]
Let $\Omega$ be a smooth bounded annular domain whose boundary
components are not concentric circles. We consider the harmonic
function $h$ equal to $0$ on the inner boundary and $1$ on the outer
boundary. By elliptic regularity, $h\in C^\infty(\overline\Omega)$.
By \cite[Lemma~2.9]{AdamowiczVeronelli2022}, we have $\nabla h\ne0$
in $\Omega$. The maximum principle and Hopf lemma give the same
conclusion on $\partial\Omega$. Thus, by continuity and compactness,
\[
 \min_{\overline\Omega}|\nabla h|>0.
\]
The vorticity law is $F=0$,
and its Dirichlet linearization has trivial kernel.
Thus, $u=\nabla^\perp h$ satisfies the hypotheses of
Theorem~\ref{intro:approximation}, with $E(u)>0$. Applying the theorem,
we obtain nonzero isolated steady Euler flows with the same energy.
The construction in its proof gives finitely valued vorticity,
and this property is preserved by the final energy normalization.
\end{proof}

\section{Failure of attraction at fixed energy and circulations. Proof of Theorem~\ref{intro:theorem-c}.}
\label{sec:no-attraction}

In this section we prove Theorem~\ref{intro:theorem-c}. We first
explain how to prescribe the energy and the boundary circulations.
We then treat the case of constant vorticity by constructing nearby
steady Euler flows with the same constraints. When the vorticity is
nonconstant, we construct perturbations that decrease either the
enstrophy or the fourth vorticity moment while preserving these
constraints. Conservation of this quantity will then exclude weak
return, and we will use the equality of energies to obtain a positive
distance from the whole scaling line. Neither argument requires
$u_*$ to be isolated or uses the particular structure of the Euler flow
constructed in Section~\ref{sec:construction}. Throughout this section,
$\Omega$ is as in Theorem~\ref{intro:theorem-c}.

We begin by describing how to recover the velocity from its vorticity
and prescribed inner circulations. We write the boundary of $\Omega$ as
$\Sigma_0\cup\Sigma_1\cup\cdots\cup\Sigma_m$, where $\Sigma_0$ is
the outer component. All components are oriented by $\tau=n^\perp$.
We fix the inner circulations
$\gamma_j=\Gamma_j(u_*)$, $1\leq j\leq m$, and, for
$f\in L^\infty(\Omega)$, let $u_f=\nabla^\perp\psi_f$ be the unique
velocity satisfying
\begin{equation}\label{na:recovery}
 \Delta\psi_f=f,\qquad
 \psi_f|_{\Sigma_0}=0,\qquad
 \psi_f|_{\Sigma_j}\text{ is constant},\qquad
 \int_{\Sigma_j}\partial_n\psi_f\dd s=\gamma_j.
\end{equation}
The recovery map $f\mapsto\psi_f$ is affine and continuous from
$L^p(\Omega)$ to $W^{2,p}(\Omega)$ for $1<p<\infty$.

The existence and uniqueness of $\psi_f$, together with its affine
continuous dependence on $f$, follow from the classical planar Hodge
decomposition and standard Dirichlet estimates, see for instance
\cite[Section~4]{Sueur2011} and \cite[Chapter~9]{GT}.
If $\Omega$ is simply connected, then $m=0$ and \eqref{na:recovery}
reduces to the Dirichlet problem $\Delta\psi_f=f$ in $\Omega$, with
$\psi_f=0$ on $\partial\Omega$. The only boundary circulation is then
$\Gamma_0(u_f)=\int_\Omega f\dd x$.

We now define $\mathcal E(f)=\frac12\|u_f\|_{L^2}^2$ and describe
its variation. For a vorticity increment $h$, we set
\[
 \xi_h:=\psi_{f+h}-\psi_f,\qquad Kh:=\nabla^\perp\xi_h.
\]
Subtracting the equations in \eqref{na:recovery}, we see that
$\Delta\xi_h=h$, that $\xi_h$ vanishes on $\Sigma_0$, it is constant
on each $\Sigma_j$, and that
$\int_{\Sigma_j}\partial_n\xi_h\dd s=0$ since $\gamma_j$ is fixed, for $1\leq j\leq m$.
These conditions determine $\xi_h$ uniquely: indeed, the difference of any two solutions would be harmonic, {constant and with zero circulation on each $\Sigma_i$}, and with zero trace on $\Sigma_0$, thus zero on $\Omega$. Hence,
$K:L^p(\Omega)\to W^{1,p}(\Omega;\R^2)$ is a bounded linear map
for $1<p<\infty$, and $u_{f+h}=u_f+Kh$. In particular, the velocity
increment $Kh$ has zero inner circulations. Since $\curl Kh=h$, these
estimates also give $\|Kh\|_{\mathcal Y}\leq C_\Omega\|h\|_{L^\infty}$
for bounded $h$. Integration by parts then gives 
\begin{equation}\label{na:energy-variation}
 \mathcal E(f+h)-\mathcal E(f) =\int_\Omega \nabla \psi_f\cdot \nabla \xi_h \dd x+\frac12\|Kh\|_{L^2}^2
   =-\int_\Omega\psi_f h\dd x+\frac12\|Kh\|_{L^2}^2,
\end{equation}
and we note that
\begin{equation*}
    D\mathcal E(f)[h]=\frac{\dd}{\dd \varepsilon}\Big|_{\varepsilon=0} \mathcal E(f+\varepsilon h)=-\int_\Omega\psi_f h\dd x.
\end{equation*}
Indeed, 
\begin{align*}
    \int_\Omega \nabla \psi_f\cdot \nabla \xi_h \dd x = - \int_\Omega \psi_f h \dd x  + \int_{\Sigma_0} \psi_f \partial_n\xi_h \dd s+ \sum_{j=1}^m\int_{\Sigma_j} \psi_f \partial_n \xi_h \dd s = -\int_\Omega \psi_f h \dd x
\end{align*}
since the boundary term on $\Sigma_0$ vanishes because $\psi_f=0$
there and on each inner component it vanishes because $\psi_f$ is
constant and $\int_{\Sigma_j}\partial_n\xi_h\dd s=0$.
On the other hand, Stokes' formula gives
\begin{equation}\label{na:outer-circulation}
 \Gamma_0(u_f)=\int_\Omega f\dd x-\sum_{j=1}^m\gamma_j,
\end{equation}
and thus
\begin{align*}
    \Gamma_0(u_{f+h}) = \Gamma_0(u_f) + \int_\Omega h \dd x,
\end{align*}
Consequently, it suffices to
impose $\int_\Omega h\dd x=0$ to preserve all boundary circulations. Notice that $\mathcal E$ is the
kinetic energy of the full velocity $u_f$, including its harmonic
part.

From now on, we write $\psi_*:=\psi_{\omega_*}$ for the normalized
stream function of $u_*$. Since $u_*\ne0$, the function $\psi_*$ is
nonconstant. In particular,
the functions $1$ and $\psi_*$ are linearly independent. We write
$\Pi$ for the $L^2(\Omega)$ orthogonal projection onto 
\begin{equation}\label{na:projection-space}
 \mathcal{V}=\operatorname{span}\{1,\psi_*\}.
\end{equation}
Since $\psi_*$ is bounded, $\Pi$ maps $L^\infty(\Omega)$ into itself.

We first consider constant vorticity. In this case, we construct nearby
steady Euler flows with the same energy and boundary circulations that will serve as the initial data for Theorem~\ref{intro:theorem-c}.

\begin{lemma}\label{na:constant-vorticity}
Let $u_*\in\mathcal Y(\Omega)\setminus\{0\}$ have constant
vorticity. There is a family of steady Euler flows $u_c\in\mathcal Y(\Omega)$,
defined for sufficiently small $c\ne0$, with nonconstant vorticity
and with the same energy and all the same boundary circulations as
$u_*$. Moreover, $\|u_c-u_*\|_{\mathcal Y}\to0$ and
$u_c\to u_*$ in $W^{1,p}(\Omega)$ and $C^{0,\alpha}(\overline\Omega)$
for every $1<p<\infty$ and $0<\alpha<1$.
In particular, $u_c$ is not a multiple of $u_*$, and $u_*$ is not
isolated modulo scaling. If $\partial\Omega$ is smooth, then
$u_c\to u_*$ also in $C^\infty(\overline\Omega)$.
\end{lemma}

\begin{proof}
We write $\omega_*=q\in\R$ and keep the inner circulations in
\eqref{na:recovery} fixed. We seek solutions $(f,a,b,c)\in L^\infty(\Omega)\times \R^3$ to
\begin{equation}\label{na:constant-branch}
 f=q+a+b\psi_f+c\psi_f^2,\qquad
 \int_\Omega f\dd x=q|\Omega|,\qquad
 \mathcal E(f)=\mathcal E(q).
\end{equation}
We regard $c$ as a parameter and solve for
$(f,a,b)\in L^\infty(\Omega)\times\R^2$ near $(q,0,0)$. We define the smooth map
\begin{align*}
    \mathcal G:L^\infty(\Omega)\times \R^3 \rightarrow L^\infty(\Omega)\times \R^2
\end{align*}
given by
\begin{align*}
    \mathcal G(f,a,b,c) := \left(f-q-a-b\psi_f - c\psi_f^2, \int_\Omega f \dd x - q|\Omega|, \mathcal E(f) - \mathcal E(q) \right).
\end{align*}
We note that $\mathcal G(q,0,0,0) = 0$ and we define
\begin{align}\label{eq:def_lin_op_G}
    L[h,\alpha,\beta]:=D_{(f,a,b)}\mathcal G(q,0,0,0)[h,\alpha,\beta] = \left( h - \alpha - \beta\psi_*, \int_\Omega h \dd x, -\int_\Omega \psi_* h \dd x \right).
\end{align}
To show that $L$ is an isomorphism, let $(g,\Gamma, \mathrm E)\in L^\infty(\Omega)\times \R^2$. We shall find $(h,\alpha,\beta)\in L^\infty(\Omega)\times \R^2$ such that $L(h,\alpha,\beta)=(g,\Gamma,\mathrm E)$. Setting $h:= g + \alpha + \beta\psi_*$ we obtain
\begin{align*}
 \begin{pmatrix}
 |\Omega|&\displaystyle\int_\Omega\psi_*\dd x\\[2mm]
 \displaystyle\int_\Omega\psi_*\dd x&
 \displaystyle\int_\Omega\psi_*^2\dd x
 \end{pmatrix} \begin{pmatrix}
     \alpha \\[2mm] \beta
 \end{pmatrix} = \begin{pmatrix}
     \Gamma - \int_\Omega g \dd x \\[2mm]
     - \mathrm E - \int_\Omega \psi_* g \dd x
 \end{pmatrix}.
\end{align*}
Since $1$ and $\psi_*$ are linearly
independent, Cauchy-Schwarz does not saturate and $\left(\int_\Omega \psi_* \dd x\right)^2 < |\Omega|\int_\Omega \psi_*^2 \dd x$. Namely, the matrix is invertible and there exists a unique solution $(\alpha,\beta)$. Hence, $L$ is an isomorphism. We then apply the implicit function theorem to
obtain a smooth branch
$c\mapsto(f_c,a(c),b(c))$ satisfying
\eqref{na:constant-branch} for $c$ sufficiently small.

We next establish the convergence. Smoothness of the branch in
$L^\infty(\Omega)\times\R^2$ gives
\[
 \|f_c-q\|_{L^\infty}+|a(c)|+|b(c)|\leq C|c|.
\]
Using the estimates for \eqref{na:recovery}, we infer that
$\psi_{f_c}-\psi_*=O(c)$ in $W^{2,p}(\Omega)$ for every $1<p<\infty$,
and hence in $C^{1,\alpha}(\overline\Omega)$ for each $0<\alpha<1$ after choosing
$p$ sufficiently large. Thus, $u_c:=u_{f_c}\in\mathcal Y(\Omega)$
converges to $u_*$ in $\mathcal Y(\Omega)$, $W^{1,p}(\Omega)$ and
$C^{0,\alpha}(\overline\Omega)$ for the stated exponents. The identity
\[
 f_c-q=a(c)+b(c)\psi_{f_c}+c\psi_{f_c}^2
\]
then gives $f_c-q=O(c)$ in $C^{1,\alpha}(\overline\Omega)$.

If $\partial\Omega$ is smooth, then $\psi_*$ is smooth, since its
Laplacian is constant and its boundary values are constant on each
component. By the Schauder estimates for \eqref{na:recovery}, we
obtain the corresponding $C^{3,\alpha}(\overline\Omega)$ bound for
$\psi_{f_c}-\psi_*$. Repeating this argument, we obtain,
for every integer $r\geq0$,
\[
 \|f_c-q\|_{C^{r,\alpha}}
 +\|\psi_{f_c}-\psi_*\|_{C^{r+2,\alpha}}
 \leq C_{r,\alpha}|c|.
\]
Thus, when $\partial\Omega$ is smooth, $u_c$ is smooth and converges
to $u_*$ in $C^\infty(\overline\Omega)$.

To verify that $u_c$ is a steady Euler flow, we define the vorticity law and pressure  
\[
 F_c(s)=q+a(c)+b(c)s+cs^2,
 \qquad \mathcal{P}_c(s)=\int_0^s F_c(r)\dd r,
\]
and use the identity
\[
 (u_c\cdot\nabla)u_c
   =\nabla\left(\frac12|u_c|^2-\mathcal P_c(\psi_{f_c})\right)
   \qquad\text{a.e. in }\Omega,
\]
which holds under the $C^{1,1}$ boundary assumption, to obtain the
pressure $p_c=\mathcal{P}_c(\psi_{f_c})-\frac12|u_c|^2$,
which is defined throughout $\Omega$.

The equality of energies implies that $u_c\ne0$, so that
$\psi_{f_c}$ is nonconstant. Its continuous range contains an
interval of positive length. For $c\ne0$, the polynomial $F_c$ cannot be
constant on this interval, and therefore $f_c$ is nonconstant.
Thus, $u_c$ cannot be a multiple of $u_*$. 

Finally, the energy and the total
vorticity are fixed by \eqref{na:constant-branch}. Together with
the prescribed inner circulations, \eqref{na:outer-circulation}
therefore gives $\Gamma(u_c)=\Gamma(u_*)$. Notice that the argument
also applies when $q=0$, since the nonzero harmonic velocity $u_*$
still has a nonconstant stream function.
\end{proof}

We next consider nonconstant vorticity and construct perturbations that
decrease a convex Casimir while preserving the energy and boundary
circulations. These perturbations will then be shown to remain far from the ray generated by $u_*$ for all time.

\begin{proposition}\label{na:lower-casimir}
Let $u_*\in\mathcal Y(\Omega)\setminus\{0\}$ have nonconstant
vorticity. For one of the functions $H(s)=s^2$ and $H(s)=s^4$,
there exist $\eps_0>0$ and a family
$u_\eps\in\mathcal Y(\Omega)$, $0<\eps<\eps_0$, such that
\begin{gather}
 E(u_\eps)=E(u_*),\qquad
 \Gamma(u_\eps)=\Gamma(u_*),\label{na:exact-constraints}\\
 0<\|u_\eps-u_*\|_{\mathcal Y}\longrightarrow0,
 \qquad
 C_H(\omega_\eps)<C_H(\omega_*),\label{na:casimir-decrease}
\end{gather}
where $\omega_\eps=\curl u_\eps$ and
$C_H(f)=\int_\Omega H(f)\dd x$.
\end{proposition}

\begin{proof}
We first show that at least one of $\omega_*$ and $\omega_*^3$ does
not belong to the space $\mathcal{V}= =\operatorname{span}\{1,\psi_*\}$.
We argue by contradiction and assume that both belong to
$\mathcal{V}$. Then, we can write
\[
 \omega_*=a\psi_*+b,\qquad \omega_*^3=d\psi_*+e
 \qquad\text{a.e. in }\Omega.
\]
Since $\omega_*$ is nonconstant, we have $a\ne0$. Thus, it follows that
\[
 (a\psi_*+b)^3-d\psi_*-e=0
\]
a.e. in $\Omega$, and hence everywhere by continuity. Since the
range of $\psi_*$ contains an interval of positive length, the
polynomial $(as+b)^3-ds-e$ vanishes on that interval. This is
impossible because its leading coefficient is $a^3\ne0$.

We can therefore choose $\ell\in\{2,4\}$ such that
$G:=\omega_*^{\ell-1}\notin\mathcal{V}$. We set $H(s):=s^\ell$ and
define
\begin{equation}\label{na:directions}
 h=-(\operatorname{Id}-\Pi)G,\qquad
 k=\psi_*-\frac1{|\Omega|}\int_\Omega\psi_*\dd x.
\end{equation}
Both functions are bounded. By the definition of $\Pi$, we have $(\operatorname{Id}-\Pi)G \perp \mathcal V$ , so that $(\operatorname{Id}-\Pi)G \perp \Pi G$ and
\begin{equation}\label{na:first-order-constraints}
 \int_\Omega h\dd x=\int_\Omega\psi_*h\dd x=0,
 \qquad
 \int_\Omega Gh\dd x=-\|(\operatorname{Id}-\Pi)G\|_{L^2}^2<0.
\end{equation}
Thus, the perturbation $\omega_*+\eps h$ preserves the total
vorticity and has zero first variation of the energy, while the
first variation of $C_H$ is strictly negative. We will use $k$ to
correct the quadratic error in the energy. Indeed, $k$ has zero mean and
\begin{equation}\label{na:energy-correction-direction}
 D\mathcal E(\omega_*)[k]
   =-\int_\Omega\psi_*k\dd x = -\int_\Omega \left( k + \frac{1}{|\Omega|}\int_\Omega \psi_*(y) \dd y \right) k \dd x =-\|k\|_{L^2}^2<0.
\end{equation}
We now apply the scalar implicit function theorem to the map $\mathcal G:\R^2 \rightarrow \R$ given by
\[
 \mathcal G(\eps,a):= \mathcal E(\omega_*+\eps h+ak)-\mathcal E(\omega_*).
\]
By \eqref{na:energy-correction-direction}, its derivative with
respect to $a$ at $(0,0)$ is nonzero,
\begin{align*}
    D_a\mathcal{G}(0,0) = D\mathcal E(\omega_*)[k]= -\Vert k \Vert_{L^2}^2 < 0.   
\end{align*}
We therefore obtain $\eps_0>0$ and a smooth function $a(\eps)$ with $a(0)=0$ such that 
\begin{align}\label{na:perturbed-vorticity}
    \omega_\eps:=\omega_*+\eps h+a(\eps)k,
\end{align}
satisfies $\mathcal E(\omega_\eps)=\mathcal E(\omega_*)$ for all $|\eps|<\eps_0$. Moreover, $a'(0)=0$ and $a(\eps)=O(\eps^2)$. Indeed, $\mathcal G(\eps,a(\eps)) = 0$ gives 
\begin{align*}
    0 = \frac{\dd}{\dd \eps}\Big|_{\eps=0} \mathcal{G}(\eps,a(\eps)) = \int_\Omega \psi_* h \dd x - a'(0) \int_\Omega \psi_* k \dd x = - a'(0) \int_\Omega \psi_* k \dd x
\end{align*}
readily gives $a'(0)=0$. We define $u_\eps:=u_{\omega_\eps}$ using \eqref{na:recovery}, so
the inner circulations are fixed. Since both $h$ and $k$ have zero
integral, \eqref{na:outer-circulation} also gives equality of the
outer circulations. We have therefore proved \eqref{na:exact-constraints}. Finally, expanding $C_H$ and using $a(\eps)=O(\eps^2)$, we obtain
\begin{align*}
 C_H(\omega_\eps)-C_H(\omega_*)
  &=\ell\eps\int_\Omega\omega_*^{\ell-1}h\dd x
       +O(\eps^2)=-\ell\eps\|(\operatorname{Id}-\Pi)G\|_{L^2}^2+O(\eps^2).
\end{align*}
The right-hand side is strictly negative for all sufficiently small
$\eps>0$. On the other hand, \eqref{na:perturbed-vorticity} and the
bound for $K$ give $\|u_\eps-u_*\|_{\mathcal Y}=O(\eps)$.
Since the values of $C_H$ are different, $\omega_\eps\ne\omega_*$,
and hence $u_\eps\ne u_*$. This proves
\eqref{na:casimir-decrease}.
\end{proof}

We now complete the proof by using these perturbations to exclude weak
return and obtain a positive distance from the scaling line.

\begin{proof}[Proof of Theorem~\ref{intro:theorem-c}]
We first suppose that $\omega_*=q$ is constant. We choose a
sufficiently small nonzero $c$ in Lemma~\ref{na:constant-vorticity}
and set $u_0:=u_c$. Since $u_c$ is a steady Euler flow, uniqueness
of Yudovich solutions gives $u(t)=u_c$ for all $t\in\R$. Thus,
the weak-$*$ closure of the vorticity orbit is $\{f_c\}$, which does
not contain $\omega_*$ because $f_c$ is nonconstant. Moreover,
$u_c\notin\R u_*$ and $\R u_*$ is closed in $L^2(\Omega)$, so
\[
 \inf_{t,\lambda\in\R}
       \|u(t)-\lambda u_*\|_{L^2}
 =\operatorname{dist}_{L^2}(u_c,\R u_*)>0.
\]
Lemma~\ref{na:constant-vorticity} gives equality of energy and boundary circulations of $u_c$, as well as convergence to $u_*$ as $c\to0$. This proves all the assertions when $\omega_*$ is constant.

Suppose now that $\omega_*$ is nonconstant. We choose a sufficiently
small perturbation $u_0=u_\eps$ from
Proposition~\ref{na:lower-casimir} and let $u(t)$ be the global
Yudovich solution with initial datum $u_0$. We set
\[
 \delta=C_H(\omega_*)-C_H(\omega_0)>0.
\]
Since the vorticity is transported by a flow that preserves volume,
we have
\begin{equation}\label{na:conserved-casimir}
 C_H(\omega(t))=C_H(\omega_0)=C_H(\omega_*)-\delta
 \qquad\text{for all }t\in\R.
\end{equation}
We first show that neither $\omega_*$ nor $-\omega_*$ belongs to
the weak-$*$ closure of the vorticity orbit. Let the sign
$\sigma\in\{\pm 1\}$ be fixed but arbitrary. Since $H$ is convex
and even, we infer from \eqref{na:conserved-casimir} that
\begin{equation}\label{na:weak-separation}
 \int_\Omega H'(\sigma\omega_*)
       \big(\omega(t)-\sigma\omega_*\big)\dd x
 \leq C_H(\omega(t))-C_H(\sigma\omega_*)=-\delta.
\end{equation}
The function $H'(\sigma\omega_*)$ belongs to $L^1(\Omega)$, so
the left-hand side is weak-$*$ continuous as a function of
$\omega(t)$. Therefore, \eqref{na:weak-separation} holds on the
weak-$*$ closure of the entire orbit. Since its left-hand side
vanishes at $\sigma\omega_*$ and $\delta>0$, we conclude that
\begin{equation}\label{na:no-weak-return}
 \{\omega_*,-\omega_*\}\cap
 \overline{\{\omega(t):t\in\R\}}^{\,\mathrm{w}^*}=\varnothing.
\end{equation}

We next prove the corresponding statement for the velocities.
For $\sigma\in\{\pm 1\}$, we set
\[
 d_\sigma=\inf_{t\in\R}\|u(t)-\sigma u_*\|_{L^2}.
\]
We argue by contradiction and assume that $d_\sigma=0$. Then,
there is a sequence $t_j\in\R$ such that
$u(t_j)\to\sigma u_*$ in $L^2(\Omega)$. Taking the curl, we obtain
$\omega(t_j)\to\sigma\omega_*$ in the sense of distributions in $\Omega$.
Moreover, $\|\omega(t_j)\|_{L^\infty}=\|\omega_0\|_{L^\infty}$.
Using this uniform bound and the density of $C_c^\infty(\Omega)$
in $L^1(\Omega)$, we deduce that the convergence also holds in the
weak-$*$ topology of $L^\infty(\Omega)$. This contradicts
\eqref{na:no-weak-return}. We have therefore proved that
\begin{equation}\label{na:two-sign-distance}
 d:=\min\{d_+,d_-\}>0.
\end{equation}
In particular, both assertions in \eqref{intro:nonreturn} hold.

Finally, conservation of energy and
\eqref{na:exact-constraints} give
$\|u(t)\|_{L^2}=\|u_*\|_{L^2}$ for all $t\in\R$. Let $c\in\R$
be fixed but arbitrary. We choose $\sigma=1$ if $c\geq0$ and
$\sigma=-1$ otherwise. By the reverse triangle inequality, we have
\[
 \big||c|-1\big|\,\|u_*\|_{L^2}\leq\|u(t)-cu_*\|_{L^2}.
\]
Combining this estimate with \eqref{na:two-sign-distance}, we obtain
\[
 d\leq\|u(t)-\sigma u_*\|_{L^2}
 \leq\|u(t)-cu_*\|_{L^2}+\big||c|-1\big|\,\|u_*\|_{L^2}
 \leq2\|u(t)-cu_*\|_{L^2}.
\]
Taking the infimum over $t,c\in\R$, we conclude that
\[
 \inf_{t,c\in\R}
       \|u(t)-cu_*\|_{L^2}\geq\frac d2>0.
\]
Since the initial perturbation can be chosen arbitrarily small,
the proof is complete.
\end{proof}

\section{Further remarks and extensions}
\label{app:constraints}

Let us first prove the assertions in Remarks~\ref{intro:immediate-departures}
and~\ref{intro:regularity}. We first construct nontrivial perturbations
preserving the energy, enstrophy and boundary circulations.
Let $u_*=\nabla^\perp\psi_*$ be the steady Euler flow constructed in
Corollary~\ref{con:constructed}, and write
\[
 \omega_*=F_*(\psi_*),\qquad
 Q_*:=\operatorname{ess\,ran}F_*=\{q_1,\ldots,q_N\},\qquad
 2\delta_*:=\min_{i\ne j}|q_i-q_j|.
\]

\begin{proposition}\label{dyn:invariant-perturbations}
For every nonempty open set $U\Subset\Omega$, there is a smooth
family of functions $g_\varepsilon\in C_c^\infty(U)$, with supports
in a fixed compact subset of $U$, such that, for all sufficiently
small nonzero $\varepsilon$,
\[
 v_\varepsilon:=u_*+\nabla^\perp g_\varepsilon
 \in\mathcal M_Z(u_*)\setminus\{u_*\}.
\]
The differences $v_\varepsilon-u_*$ and
$\operatorname{curl}v_\varepsilon-\omega_*$ tend to zero in every
$C^k$ norm. In particular, there are nontrivial data in
$\mathcal M_Z(u_*)$ arbitrarily close to $u_*$ in $\|\cdot\|_{\mathcal Y}$.
\end{proposition}

\begin{proof}
Let $U\Subset\Omega$ be a nonempty open set, fixed but arbitrary.
The derivatives of the energy and enstrophy constraints in a direction
$\nabla^\perp f$, with $f\in C_c^\infty(U)$, are
\begin{equation}\label{dyn:constraint-derivatives}
 L_E(f):=-\int_\Omega\omega_*f\,dx,
 \qquad
 L_Z(f):=2\int_\Omega\omega_*\Delta f\,dx.
\end{equation}
We first show that these functionals are linearly independent on
$C_c^\infty(U)$.

To that end, let us start pointing out that the vorticity $\omega_*$
is not essentially constant on any nonempty open subset $O\subset\Omega$. Indeed, the inverse image
of $O$ under the parametrization \eqref{con:constructed-parametrization}
contains a rectangle $I_0\times J_0$.
If $\omega_*=c$ a.e. in $O$, change of variables and
Fubini would imply $F_*(s)=c$ a.e. in $I_0$, contrary
to its construction.

Suppose now that $aL_E(f)+bL_Z(f)=0$ for all $f\in C^\infty_c(U)$. Then
\[
 2b\Delta\omega_*-a\omega_*=0
 \quad\text{in }\mathcal D'(U).
\]
If $b=0$ and $a\ne0$, this contradicts the local nonconstancy just
proved. If $b\ne0$, interior elliptic regularity for
$\Delta\omega_*=(a/(2b))\omega_*$, followed by iteration, gives a
smooth representative of $\omega_*$ on $U$. A continuous function
that lies in the finite set $Q_*$ a.e. lies in $Q_*$
everywhere, and hence is constant on each connected ball. This is
again a contradiction. Therefore $a=b=0$, and thus $L_E$ and $L_Z$ are linearly independent.

We now choose $f_1,f_2\in C_c^\infty(U)$ such that
\[
 M:=\begin{pmatrix}
 L_E(f_1)&L_E(f_2)\\
 L_Z(f_1)&L_Z(f_2)
 \end{pmatrix}
\]
is invertible. Since $C_c^\infty(U)$ is infinite dimensional, we can choose
$f_3\in C_c^\infty(U)\setminus\operatorname{span}\{f_1,f_2\}$ and set
\[
 \begin{pmatrix}c_1\\c_2\end{pmatrix}
 =M^{-1}\begin{pmatrix}L_E(f_3)\\L_Z(f_3)\end{pmatrix},
 \qquad h:=f_3-c_1f_1-c_2f_2\neq 0.
\]
Then, we have $L_E(h)=L_Z(h)=0$ and
$h\notin\operatorname{span}\{f_1,f_2\}$, so $h,f_1,f_2$ are linearly independent. For
$f=\varepsilon h+\alpha f_1+\beta f_2$, we consider
\[
 \mathcal G(\varepsilon,\alpha,\beta)
 = \begin{pmatrix}
     \mathcal{G}_E(\varepsilon,\alpha,\beta)\\
     \mathcal{G}_Z(\varepsilon,\alpha,\beta)
 \end{pmatrix} := \begin{pmatrix}
 E(u_*+\nabla^\perp f)-E(u_*)\\
 Z(\omega_*+\Delta f)-Z(\omega_*)
 \end{pmatrix},
\]
and write $\mathcal G=(\mathcal G_E,\mathcal G_Z)$ for its two components.
Using the compact support of $f$ and $\Delta\psi_*=\omega_*$, we obtain the exact formulas
\begin{align}
 \mathcal G_E
 &=-\int_\Omega\omega_*f\,dx
       +\frac12\int_\Omega|\nabla f|^2\,dx,
       \label{dyn:energy-polynomial}\\
 \mathcal G_Z
 &=2\int_\Omega\omega_*\Delta f\,dx
       +\int_\Omega(\Delta f)^2\,dx.
       \label{dyn:enstrophy-polynomial}
\end{align}
Thus, we see that $\mathcal G(0,0,0)=0$, and its derivative with respect to
$(\alpha,\beta)$ at zero is $M$, which is invertible. We then apply the finite-dimensional implicit
function theorem to obtain smooth functions
$\alpha=a(\varepsilon)$ and $\beta=b(\varepsilon)$ such that
\[
 \mathcal G(\varepsilon,a(\varepsilon),b(\varepsilon))=0,
 \qquad a(0)=b(0)=0.
\]
Since $L_E(h)=L_Z(h)=0$, we differentiate the preceding identity at zero and obtain
$a'(0)=b'(0)=0$. In particular, we have $a(\varepsilon),b(\varepsilon)
=O(\varepsilon^2)$.

Finally, we define
\[
 g_\varepsilon
 :=\varepsilon h+a(\varepsilon)f_1+b(\varepsilon)f_2\in C^\infty_c(U).
\]
The energy and enstrophy are exactly preserved for $u_*+\nabla^\perp g_\eps$. The velocity
perturbation $\nabla^\perp g_\eps$ vanishes near the boundary, so every boundary circulation
is preserved as well. The fixed support and the preceding estimates
give convergence of $\nabla^\perp g_\eps$ to 0 in every $C^k$ norm.
For $\varepsilon\ne0$, linear independence gives
$g_\varepsilon\ne0$. Moreover, $\Delta g_\varepsilon\ne0$:
otherwise integration by parts would give
$\int|\nabla g_\varepsilon|^2=0$, forcing the compactly supported
function $g_\varepsilon$ to vanish. Thus the velocity and vorticity
perturbations are both nonzero, as claimed.
\end{proof}

\begin{proof}[Proof of the immediate departure statement {following} Remark~\ref{intro:immediate-departures}]
Proposition~\ref{dyn:invariant-perturbations} supplies the
nearby nontrivial data asserted in Remark~\ref{intro:immediate-departures}. Every such datum
is nonstationary: stationary rigidity and equality of energy would
otherwise give $u_0=u_*$. For any $0<\rho<\delta_*$ and any $T>0$,
its orbit must contain a time $t\in(0,T)$ with
$\|\omega(t)-\omega_*\|_{L^\infty}>\rho$. If not, the orbit would
satisfy the rigidity hypothesis throughout $(0,T)$ and hence would
equal $u_*$ globally, which is impossible. Applying the same
argument to $(-T,0)$ proves the asserted immediate departures in
both time directions.
\end{proof}

We next make precise the extension of these results to $W^{k,\infty}$ and the remark that our equilibria are not isolated in weak-$*$ $L^\infty$. This follows essentially from estimates we have already established, so we only provide a sketch of the proof of this result. As is customary, we use the convention $W^{0,p}=L^p$. 

\begin{proposition}\label{app:finite-order-isolation}
Let $\Omega$ and $u$ satisfy the hypotheses of
Theorem~\ref{intro:approximation}, and suppose that $\Omega$ is not
a concentric circular annulus. For each integer $k\ge0$, there are
steady Euler flows $u_n$, with vorticities $\omega_n$, such that
\[
 E(u_n)=E(u),\qquad
 \|u_n-u\|_{L^2}+\|\omega_n-\curl u\|_{W^{k,\infty}}\longrightarrow0,
 \qquad \omega_n\in W^{k,\infty}(\Omega)\setminus C^k(\overline\Omega).
\]
For each $u_*=u_n$ there is a $\delta_*>0$ such that any Yudovich
solution $v$ satisfying
\[
 \operatorname*{sup}_{t\in J}
 \|\curl v(t)-\curl u_*\|_{W^{k,\infty}}<\delta_*
\]
on a nonempty open interval $J$ is a constant positive multiple of
$u_*$ for all time. Nevertheless, there are smooth steady flows
$v_j\notin\R u_*$ of the same energy as $u_*$ such that
\[
 \curl v_j\longrightarrow\curl u_*
 \quad\text{in }W^{k,p}(\Omega)\quad\text{for every }1\le p<\infty.
\]
At order zero their vorticities also converge weakly-$*$ in $L^\infty(\Omega)$.
\end{proposition}

\begin{proof}
Fix a smooth pair $(F_0,\psi_0)$ given by
Proposition~\ref{app:den-density}. For $k=0$ take $F=F_\eps$ as
in~\eqref{con:profile}. For $k\ge1$ take $F=F_0+\eps P_k$, where
$P_k$ is the $k$th primitive of $q:=\mathbf1_B$ with zero initial
data at the left endpoint of $I$.

To obtain the smooth approximants at the same time, extend and mollify
$F-F_0$ and add back $F_0$. The resulting smooth laws $F_j$ converge
to $F$ in $W^{k,p}(I)$ for every finite $p$ and satisfy
$\|F_j-F_0\|_{W^{k,\infty}(I)}\le C_k\eps$. Apply
Proposition~\ref{con:rough-approximation} to obtain stream functions
$\psi_j$ with the required level geometry. Starting from
\eqref{con:rough-estimates}, successive finite-$p$ Dirichlet estimates
for $\Delta(\psi_j-\psi_0)=F_j(\psi_j)-F_0(\psi_0)$ give, uniformly
in $j$,
\[
 \|\psi_j-\psi_0\|_{W^{k+2,p}}
 +\|F_j(\psi_j)-F_0(\psi_0)\|_{W^{k,\infty}}
 \le C_{k,p}\eps \qquad (1<p<\infty).
\]
Pass to a subsequence with $\psi_j\to\psi$ in $C^{k+1}$. At the
highest order, Lemmas~\ref{con:pushforward} and~\ref{con:composition}
give $F_j^{(k)}(\psi_j)\to F^{(k)}(\psi)$ in every finite $L^p$.
The other factors in the chain rule converge uniformly. Thus
$F_j(\psi_j)\to\omega_*:=F(\psi)=\Delta\psi$ in $W^{k,p}$.
At $k=0$ the uniform vorticity bound also gives weak-$*$ convergence.

For $k=0$, weak-$*$ continuity extends the essential bound to every
time in $J$, so Lemma~\ref{rig:rigidity} gives rigidity. For $k\ge1$
we recover $q(\psi)$ by differentiating across the levels. Set
$D=(\nabla\psi/|\nabla\psi|^2)\cdot\nabla$. Since $D\psi=1$,
\[
 \eps^{-1}D^k\bigl(\omega_*-F_0(\psi)\bigr)=q(\psi)
 \qquad\text{a.e.}
\]
For $\varpi=\curl v$, set $h:=\eps^{-1}D^k(\varpi-F_0(\psi))$.
The coefficients of $D^k$ are bounded, so the assumed closeness gives
$|h-q(\psi)|<1/3$ almost everywhere for sufficiently small $\delta_*$.

Let {$\Xi_{t,t_0}$} be the particle map from a time $t_0$ at which the bound holds. For almost every $x$, the function $t\mapsto h(t,{\Xi_{t,t_0}(x)})$ has
a continuous representative. Indeed, div--curl estimates and energy
conservation bound $v$ in $W^{k+1,p}$ for $p>2$, so compactness and
$L^2$ continuity give $v\in C(J;C^k)$. The spatial chain rule for
$\varpi(t)=\varpi(t_0)\circ {\Xi_{t,t_0}^{-1}}$ then expresses $h(t,{\Xi_{t,t_0}(x)})$
in terms of fixed derivatives of $\varpi(t_0)$ at $x$ and coefficients
continuous in time. For almost every $t$, this function lies in
$[-1/3,1/3]$ if $q(\psi({\Xi_{t,t_0}(x)}))=0$ and in $[2/3,4/3]$ if
$q(\psi({\Xi_{t,t_0}(x)}))=1$. By continuity it stays in one interval
throughout $J$. Hence $q(\psi({\Xi_{t,t_0}(x)}))=q(\psi(x))$ for almost every
$(t,x)$.

This gives~\eqref{rig:frozen-labels} with $q$ in place of $F$.
The tangency and multiplier arguments in the proof of
Lemma~\ref{rig:rigidity} now give $v=c\nabla^\perp\psi$ on $J$.
Taking $\delta_*<\|\omega_*\|_\infty$ gives $c>0$, and Yudovich
uniqueness extends the identity to all time.

Finally, both values of $q$ occur on sets of positive measure in
every subinterval of $I$. In level coordinates the formula for
$D^k(\omega_*-F_0(\psi))$ therefore excludes a $C^k$ representative
of $\omega_*$. At $k=0$ this follows from
Lemma~\ref{con:profile-properties}. Thus the smooth approximants lie
outside the scaling line. Choosing $\eps$ separately for each smooth
pair and rescaling to the prescribed energy, as in the proof of
Theorem~\ref{intro:approximation}, completes the argument.
\end{proof}

To conclude, for the regularity statement, let $\Omega$ and $u_*$ be as in
Theorem~\ref{intro:theorem-c}. We use $\omega_*=\curl u_*$, the
normalized stream function $\psi_*$, and the operator $K$ introduced
in Section~\ref{sec:no-attraction}.

\begin{proof}[Proof of the regularity statement in Remark~\ref{intro:regularity}]
We assume that $\partial\Omega$ is smooth. For constant $\omega_*$,
we use Lemma~\ref{na:constant-vorticity}. Otherwise, we choose
$G=\omega_*^{\ell-1}$ as in Proposition~\ref{na:lower-casimir}.
Since $1,\psi_*,G$ are linearly independent in $L^2(\Omega)$, the map
$f\mapsto(\int_\Omega f\dd x,\int_\Omega\psi_*f\dd x,\int_\Omega Gf\dd x)$
from $C_c^\infty(\Omega)$ to $\R^3$ is surjective: otherwise, a
nontrivial linear combination of $1,\psi_*,G$ would vanish in the
sense of distributions. We can therefore choose
$h,k\in C_c^\infty(\Omega)$ with images $(0,0,-1)$ and $(0,1,0)$,
respectively. The argument of Proposition~\ref{na:lower-casimir}
then gives the same constraints and a strict decrease of $C_H$, with
$a(\eps)=O(\eps^2)$. Finally, elliptic regularity for
\eqref{na:recovery} gives $Kh,Kk\in C^\infty(\overline\Omega;\R^2)$.
Since $u_\eps-u_*=\eps Kh+a(\eps)Kk$, these differences converge to
zero in every $C^m$ norm. 
\end{proof}

To conclude, we show that the non-attraction property of Theorem~\ref{intro:theorem-c} applies not only to equilibria, but to any solution with precompact vorticity orbit. 

\begin{remark}\label{app:precompact-motions}
Let $\Omega$ be as in Theorem~\ref{intro:theorem-c}, and let $v(t)$
be a prescribed nonstationary Yudovich solution whose vorticity
$\varpi(t)$ has precompact orbit in $L^2$. Write
$v(0)=\nabla^\perp\psi$ and $f=\varpi(0)$. By uniqueness,
$f\notin\operatorname{span}\{1,\psi\}$, since an affine relation
would make $v(0)$ steady. The proof of
Proposition~\ref{na:lower-casimir} for the enstrophy casimir therefore gives arbitrarily close
data $u_0$ in $\mathcal Y$, with the same energy and circulations,
and $Z(\omega_0)<Z(f)$, where $\omega_0=\curl u_0$.

Every element of
$\mathcal K_v:=\overline{\{\varpi(t):t\in\R\}}^{\,L^2}$
has enstrophy $Z(f)$. Conservation and weak lower semicontinuity
therefore exclude $\mathcal K_v\cup(-\mathcal K_v)$ from the
weak-$*$ closure of the perturbed vorticity orbit. Moreover,
\[
 \inf_{c,t,t'\in\R}
       \|u(t)-c\,v(t')\|_{L^2}>0.
\]
Indeed, if this infimum vanished, there would be sequences with
$\|u(t_j)-c_jv(t'_j)\|_{L^2}\to0$. Equality of the positive energies
would then give $|c_j|\to1$. By compactness, after passing to a subsequence if necessary,
$c_j\to\sigma\in\{-1,1\}$ and $\varpi(t'_j)\to\eta\in\mathcal K_v$
in $L^2$. Taking curls and using the uniform vorticity bound gives
$\omega(t_j)\rightharpoonup\sigma\eta$ in $L^2$, yielding
\[
 Z(f)=Z(\eta)\leq \liminf Z(\omega(t_j)) = Z(\omega_0)<Z(f),
\]
a contradiction. This covers periodic motions and quasiperiodic motions with
continuous $L^2$ dependence on their phases.
\end{remark}

\section*{Acknowledgements}
{The authors would like to thank Theodore D. Drivas and Vladimir \v{S}ver\'{a}k for valuable comments} {on a preliminary version of the paper.} This work was supported by the European Research Council (ERC) under
the European Union's Horizon 2020 research and innovation programme
through grant agreement 862342 (A.E. and M.N.). A.E. also acknowledges support from the Spanish State Research Agency through grants
PID2022-136795NB-I00 and PID2025-168410NB-I00, and from the
ICMAT--Severo Ochoa grant CEX2023-001347-S. A.J.F. acknowledges support by the grants PID2023-149451NA-I00 of MCIN/AEI/\allowbreak 10.13039/\allowbreak 501100011033/FEDER,UE; and RYC2024-049142-I, MICINN (Spain).

\section*{Declaration on the use of AI tools}
The authors used ChatGPT to assist with mathematical arguments
and calculations, the identification of errors and imprecisions, and
the organization, drafting and revision of the manuscript. The authors
take full responsibility for the mathematical results, arguments and
final content of this paper.

\appendix

\section{Computation of the coefficients in Lemma~\ref{app:den-interior-response}}
\label{app:coefficient-calculation}

We prove \eqref{app:den-five-coefficients}, keeping the notation of
Lemma~\ref{app:den-interior-response}. Let $\mathsf{C}=\{\psi=c\}$ be the
constant-speed level curve, and write
\[
 q=|\nabla\psi||_{\mathsf{C}}>0,\qquad n=\frac{\nabla\psi}{q}\quad\text{on }\mathsf{C}.
\]
We choose an arclength parametrization $t\mapsto\gamma(t)$ of $\mathsf{C}$
and write $n(t)=n(\gamma(t))$. We use the normal coordinates
\[
 x=\gamma(t)+zn(t),\qquad |z|<\delta,
\]
where $\delta>0$ is sufficiently small and $z$ is the signed distance
to $\mathsf{C}$ in the direction of $n$. Since $\partial_tn=\kappa\partial_t\gamma$, the Laplacian takes the form
\begin{equation}\label{app:eq-density-fermi-operator}
 \Delta=\partial_z^2+\frac{\kappa}{1+\kappa z}\partial_z
 +\frac{1}{(1+\kappa z)^2}\partial_t^2
 -\frac{z\kappa'}{(1+\kappa z)^3}\partial_t. 
\end{equation}
As in Section~\ref{sec:approximation}, we use primes for arclength
derivatives of $\kappa$.

We first reduce the computation to the case $c=0$ and $q=1$.
To that end, we set
\[
 \widehat\psi=\frac{\psi-c}{q},\qquad
 \widehat F(s)=\frac{F(c+qs)}q.
\]
Then, we have $\widehat L=L$ and
$\widehat{\mathcal R}=\mathcal R/q$. With
$\widehat\chi(s)=\chi(c+qs)$ and $\widehat\tau=q\tau$, the corresponding
solutions satisfy $\widehat v_{\widehat\tau}=v_\tau$, since their
phases and right-hand sides coincide. Comparing their
expansions in \eqref{app:den-response-expansion}, we obtain
\begin{equation}\label{app:eq-density-coefficient-rescaling}
 b_j=q^{1-j}\widehat b_j,
 \qquad
 \widehat F(0)=\frac{F(c)}q,\quad
 \widehat F'(0)=F'(c),\quad
 \widehat F''(0)=qF''(c).
\end{equation}
Here $\widehat b_j$ denotes the coefficient for the normalized
stream function. For the rest of the computation, we omit the hats
and assume that $c=0$ and $q=1$. We write
\[
 f=F(0),\qquad g=F'(0),\qquad h=F''(0),\qquad j=F'''(0).
\]

We now compute the normal derivatives of $\psi$ along $\mathsf{C}$. We set
$p_i=(\partial_z^i\psi)|_{\mathsf{C}}$ for $1\leq i\leq5$. Using
$\Delta\psi=F(\psi)$ and the identities
$\psi(t,0)=0$, $\partial_z\psi(t,0)=1$, we obtain
\begin{equation}\label{app:eq-density-normal-jets}
\begin{aligned}
 p_1&=1,\qquad p_2=f-\kappa,\\
 p_3&=g-\kappa f+2\kappa^2,\\
 p_4&=h+gf-2g\kappa+3\kappa^2f-6\kappa^3+\kappa'',\\
 p_5&=f\kappa''-12f\kappa^3-2fg\kappa-7(\kappa')^2-11\kappa\kappa''+24\kappa^4+7g\kappa^2+g^2+3fh-4h\kappa+j.
\end{aligned}
\end{equation}
For instance, setting $z=0$ in \eqref{app:eq-density-fermi-operator},
we find $p_2=f-\kappa$. Differentiating once, we obtain
$p_3-\kappa^2p_1+\kappa p_2=gp_1$, which gives the displayed value of $p_3$.
We obtain the remaining derivatives in the same way, by comparing
powers of $z$ using
\begin{equation}\label{app:eq-density-jet-recursion}
 \partial_z^2\psi=F(\psi)-\frac{\kappa}{1+\kappa z}\partial_z\psi
 -\frac{1}{(1+\kappa z)^2}\partial_t^2\psi
 +\frac{z\kappa'}{(1+\kappa z)^3}\partial_t\psi.
\end{equation}
For $\psi(t,z)=\sum_{i=1}^5\frac{p_i(t)}{i!}z^i+O(z^6)$, we have that
\begin{align*}
    \partial_z^2 \psi(t,z) = \sum_{i=0}^3 \frac{p_{i+2}(t)}{i!}z^i + O(z^4)
\end{align*}
so that $p_{i+2}$ is determined in terms of $p_1, ..., p_{i+1}$ by \eqref{app:eq-density-jet-recursion} for all $0\leq i\leq3$. Having these derivatives at hand, we compute the amplitudes in
\eqref{app:den-amplitude-recurrence}. In the same coordinates, we have
\[
 D_\psi=(\partial_z\psi)\partial_z+(1+\kappa z)^{-2}(\partial_t\psi)\partial_t,
 \qquad S=(\partial_z\psi)^2+(1+\kappa z)^{-2}(\partial_t\psi)^2.
\]
Since $\partial_t\psi(t,0)=0$ and $p_2=f-\kappa$, we obtain
\[
 \partial_t\psi=-\frac{\kappa'}2z^2+O(z^3),\qquad
 (1+\kappa z)^{-2}(\partial_t\psi)^2=\frac{(\kappa')^2}{4}z^4+O(z^5).
\]
We must therefore keep this term when computing $b_5$.

Since $\chi(\psi)=1$ near $\mathsf{C}$, the recurrence
\eqref{app:den-amplitude-recurrence} becomes
\begin{equation}\label{app:eq-density-arithmetic-recursion}
 a_1=0,\qquad a_2=S^{-1},\qquad
 a_i=-S^{-1}\bigl(\mathcal{T}a_{i-1}+La_{i-2}\bigr),
 \qquad\textup{for all }i\geq3.
\end{equation}
We write
\[
 A_i=a_i|_{\mathsf{C}},\qquad B_i=\partial_za_i|_{\mathsf{C}},\qquad
 E_i=(La_i)|_{\mathsf{C}}.
\]
We compute these quantities from
\eqref{app:eq-density-normal-jets} and
\eqref{app:eq-density-arithmetic-recursion}. To make the calculation
explicit, let $U_i$ and $V_i$ denote the coefficients of $z^i$ in
the Taylor expansions of $U$ and $V$. We use the identities
\[
 (UV)_i=\sum_{\ell=0}^iU_\ell V_{i-\ell},\qquad
 (U^{-1})_0=U_0^{-1},\qquad
 (U^{-1})_i=-U_0^{-1}\sum_{\ell=1}^iU_\ell(U^{-1})_{i-\ell},
\]
where $U_0\ne0$ for the reciprocal formulas and $i\geq1$ in the last
identity. For the derivatives, we use
$(\partial_zU)_i=(i+1)U_{i+1}$ and
$(\partial_tU)_i=\partial_tU_i$. Moreover, by
\eqref{app:eq-density-fermi-operator}, we have
\[
 E_i=(\partial_z^2a_i)|_{\mathsf{C}}+\kappa B_i+\partial_t^2A_i-gA_i.
\]
Substituting in the recurrence, we obtain
\begin{align*}
 A_2={}&1,\\
 A_3={}&3f-4\kappa,\\
 A_4={}&15f^2-32f\kappa+12\kappa^2-3g,\\
 A_5={}&105f^3-300f^2\kappa+220f\kappa^2-45fg
          +16\kappa''-32\kappa^3+48g\kappa+3h;
\end{align*}
\begin{align*}
 B_2={}&-2f+2\kappa,\\
 B_3={}&-12f^2+24f\kappa-8\kappa^2+3g,\\
 B_4={}&-90f^3+250f^2\kappa-176f\kappa^2+42fg
          -8\kappa''+24\kappa^3-44g\kappa-3h,\\
 B_5={}&-840f^4+2940f^3\kappa-3120f^2\kappa^2+585f^2g
          -200f\kappa''+1040f\kappa^3-1110fg\kappa\\
       &\quad+32(\kappa')^2+96\kappa\kappa''-64\kappa^4+412g\kappa^2-45g^2
          -57fh+60h\kappa+3j;
\end{align*}
\begin{align*}
 E_2={}&6f^2-12f\kappa+4\kappa^2-3g,\\
 E_3={}&60f^3-168f^2\kappa+120f\kappa^2-36fg-16\kappa^3+40g\kappa+3h,\\
 E_4={}&630f^4-2220f^3\kappa+2380f^2\kappa^2-495f^2g
          +80f\kappa''-800f\kappa^3+960fg\kappa\\
       &\quad+16(\kappa')^2-32\kappa\kappa''+48\kappa^4-364g\kappa^2+45g^2
          +51fh-56h\kappa-3j.
\end{align*}
For instance, we have $A_3=-(2B_2+fA_2)=3f-4\kappa$ and
$\partial_z^2a_2|_{\mathsf{C}}=6f^2-10f\kappa+2\kappa^2-2g$. Substituting the latter
identity in $E_2=\partial_z^2a_2|_{\mathsf{C}}+\kappa B_2-gA_2$, we obtain the
formula for $E_2$.

Finally, we recover the coefficients $b_1,\ldots,b_5$. Since $H|_{\mathsf{C}}=2(f-\kappa)$ and $D_\psi|_{\mathsf{C}}=\partial_z$, we rewrite
\eqref{app:den-b-recurrence} as
\begin{equation}\label{app:eq-density-final-coefficient-rule}
 b_1=2,\qquad
 b_i=-2B_i-(4f-2\kappa)A_i-2E_{i-1},\qquad\textup{for }2\le i\le5,
 \qquad E_1=0.
\end{equation}
The first two substitutions give $b_2=-2\kappa$ and $b_3=-2f\kappa$.
Using the values of $B_4,A_4,E_3$, we obtain
\[
 b_4=16\kappa''+8\kappa^3+(2g-6f^2)\kappa.
\]
In the expression for $b_5$, the terms in
$f^4,f^2\kappa^2,f^2g,g\kappa^2,g^2,fh,j$ cancel, and we obtain
\[
 b_5=176f\kappa''+88f\kappa^3+18fg\kappa
      -96(\kappa')^2-96\kappa\kappa''-32\kappa^4-30f^3\kappa-2h\kappa.
\]
We now apply \eqref{app:eq-density-coefficient-rescaling} to recover
\eqref{app:den-five-coefficients} for arbitrary $c$ and $q>0$.
Thus, the first five normal derivatives of $\psi$ suffice to compute
all five coefficients. We use the recurrence to eliminate $a_6$
from $b_5$; the amplitude $a_7$ is only needed for the remainder
estimate \eqref{app:den-parametrix-error}.

\end{document}